\documentclass[preprint,sort&compress,12pt]{elsarticle}
\usepackage{bbm}
\usepackage{amsmath}
\usepackage{mathrsfs}
\usepackage{stmaryrd}
\usepackage{mathrsfs}
\usepackage{bm}
\usepackage{amsmath}
\usepackage{amsfonts}
\usepackage{amssymb}
\usepackage{amsthm}
\usepackage{lineno}
\usepackage[RGB]{xcolor}
\newcommand{\mm}{\mathrm}
\newcommand{\ml}{\mathcal}

\newcommand{\be}{\begin{equation}}
\newcommand{\bea}{\begin{equation}\begin{aligned}}
\newcommand{\beas}{\begin{equation*}\begin{aligned}}
\newcommand{\eeas}{\end{aligned}\end{equation*}}
\newcommand{\eea}{\end{aligned}\end{equation}}
\newcommand{\ee}{\end{equation}}

\renewcommand{\div}{{\rm div }}

\begin{document}
\begin{frontmatter}
\title{On the Inhibition of
 Thermal Convection by \\
 a Magnetic Field under Zero Resistivity}
%%%%%

%%% algebraic decay

\author[FJ]{Fei Jiang}\ead{jiangfei0591@163.com}
%\cortext[cor1]{We  certify that the general content of this article, in whole or in part, is not submitted, accepted, or published elsewhere, including conference proceedings. }
\author[sJ]{Song Jiang}
\ead{jiang@iapcm.ac.cn
}
%\author[ww]{Guochun Wu}
%\ead{guochunwu@126.com}
%\author[ww]{Xin Zhong}
%\ead{xzhong1014@amss.ac.cn}
%\author[ww]{Tong Yang}
\address[FJ]{College of Mathematics and
Computer Science, Fuzhou University, Fuzhou, 350108, China.}
\address[sJ]{Institute of Applied Physics and Computational Mathematics, P.O. Box 8009,
 Beijing 100088, China.}
%\address[ww]{Institute of Applied Mathematics, AMSS, Chinese Academy of %Sciences, Beijing 100190, China.}

\begin{abstract}
We investigate the stability and instability of the magnetic Rayleigh--B\'enard problem with zero resistivity. An stability criterion is established, under which the magnetic B\'enard problem is stable.  The proof mainly is based on a three-layers energy method and an idea of magnetic inhibition mechanism.  The stable result first mathematically verifies Chandrasekhar's assertion in 1955 that the thermal instability can be inhibited by strong magnetic field in  magnetohydrodynamic (MHD) fluid with zero resistivity (based on a linearized steady magnetic B\'enard equations). In addition, we also provide an instability criterion, under which  the  magnetic Rayleigh--B\'enard problem is unstable. The proof mainly is based on the bootstrap instability method by further developing new analysis technique. Our instability result presents that the thermal instability occurs for a small magnetic field.
\end{abstract}
\begin{keyword}
Inhibition effect; Rayleigh--B\'enard problem; Magnetohydrodynamic fluid; Stability; Thermal instability;
%\MSC[2000] 35Q35\sep  76D03.
%(2000 is the default)
\end{keyword}
\end{frontmatter}
% \linenumbers

\newtheorem{thm}{Theorem}[section]
\newtheorem{lem}{Lemma}[section]
\newtheorem{pro}{Proposition}[section]
\newtheorem{cor}{Corollary}[section]
\newproof{pf}{Proof}
\newdefinition{rem}{Remark}[section]
\newtheorem{definition}{Definition}[section]
% \linenumbers

\section{Introduction}\label{introud}
\numberwithin{equation}{section}

Thermal instability often arises when a fluid is heated from below.
The classic example of this is a horizontal layer of fluid with its lower side hotter than its upper. The basic state is then one of rest states
with light and hot fluid below heavy and cool fluid. When the temperature difference across the layer is great enough, the stabilizing effects
of viscosity and thermal conductivity are overcome by the destabilizing (thermal) buoyancy, and an overturning instability ensues
 as thermal convection: hotter part of fluid is lighter and tends to rise as colder part tends to sink according
to the action of the gravity force \cite{DPGRWHHC}.

The effect of an impressed magnetic field on the onset of thermal instability in MHD fluids is first considered by Thompson \cite{TWBTCIA} in 1951.
Then Chandrasekhar theoretically further discovered the inhibiting effect of the magnetic field on the thermal instability based on the linear magnetic Boussinesq equations in 1952 \cite{CSTPRSOTICB,CSTPRSOTICBII}. Later Nakagawa experimentally verified Chandrasekhar's linear magnetic inhibition theory in 1955 \cite{YNAEN,YNAEN2}. Since the inhibition of thermal instability by magnetic fields was discovered, many prominent writers tried to mathematically verify it, see \cite{JOSEPHDD,RSSSRM} for examples. Until 1985, Galdi successfully proved the magnetic inhibition phenomenon  by (nonlinear) magnetic Boussinesq equations  with non-zero resistivity in 1985 \cite{GGPNA}, also see \cite{GGPPMFR,GGPPMNA,MGRSNSC} for further improved mathematical results. However, the mathematical verification of magnetic inhibition phenomenon by  magnetic Boussinesq equations  with zero resistivity is still a long-standing open problem. The main aim of this paper to mathematically solve this problem. Since our result is closely relevant to Galdi's conclusion, next let recall his result.

The three-dimensional (3D) magnetic Boussinesq equations in a layer domain with height $h$ (see \cite[Chapter IV]{CSHHS} for a derivation) read as follows:
\begin{equation}\label{0101ooo}\left\{{\begin{array}{ll}
 v_t+   v\cdot \nabla v+\nabla \left(p+ \lambda|M|^2/8\pi\right)/  \rho \\
  =g (\alpha(\Theta-{\Theta}_{\mm{b}})-1)e_3
  + \nu \Delta v+\lambda {M}\cdot \nabla M/4\pi\rho , \\
  \Theta_t +v\cdot\nabla \Theta =\kappa \Delta \Theta , \\[1mm]
M_t+v\cdot \nabla M=M\cdot \nabla v+\sigma\Delta M,\\
\mm{div}v= \mm{div }M=0 . \end{array}}  \right.
\end{equation}
Next we shall explain the notations in the equations above.

The unknowns  ${v}= {v}(x,t)$, $\Theta=\Theta(x,t)$,
$M= {M}(x,t)$ and $p=p(x,t)$ denote velocity, temperature,
magnetic field and pressure of the incompressible MHD fluid, resp..
The parameters $ \rho$, $\alpha$, $\lambda $ and $g>0$ denote the density constant at some properly chosen temperature parameter $\Theta_{\mm{b}}$, the coefficient of volume expansion, the permeability of vacuum and the gravitational constant, resp.. $\kappa=k/\rho c_{\mm{v}}$ and $\nu=\mu/\rho$ represents the coefficient of thermometric conductivity and the kinematic viscosity, reps., where $c_{\mm{v}}$ is the specific heat at constant volume, $k$
the coefficient of heat conductivity and $\mu$ the coefficient of shear viscosity.
$e_3=(0,0,1)^{\mm{T}}$ stand for the vertical unit vector, $g \rho\alpha(\Theta-{\Theta}_{\mm{b}})e_3$ the buoyancy (caused by expanding with heat and contracting with cold) and $-\rho g {e}_3$ for the gravitational force.
In this paper, we consider the horizontally periodic motion of MHD fluids, and thus shall introduce a horizontally periodic domain with finite height $h$, i.e.,
$$
\Omega_h:=\{ (x_{\mm{h}},x_3)\in \mathbb{R}^3~|~x_{\mm{h}}:=(x_1,x_2)\in \mathcal{T},\ 0<x_3<h\},
$$
where $\mathcal{T} :=(2\pi L_1\mathbb{T})\times(2\pi L_2\mathbb{T}) $, $\mathbb{T}=\mathbb{R}/\mathbb{Z}$, and $2\pi L_1$, $2\pi L_2>0$
 are the periodicity lengths. In addition, we denote the boundary of $\Omega_h$  by $\partial\Omega_h:=\mathbb{T}^2\times \{y_3=0,\ h\}$.

The rest state of the above  magnetic  Boussinesq equations can be given by ${r_\mm{B}}:=(0,\bar{\Theta},\bar{M})$ with an associated
pressure profile $\bar{p}$, where $\bar{M}$ is a constant vector and called an impressive magnetic field, and the temperature profile $\bar{\Theta}$ and $\bar{p}$  depend on $x_3$ only, and satisfy
the equilibrium (or rest) state
\begin{align}
\label{201903111042}
 &\nabla \bar{p} =g \rho (\alpha(\bar{\Theta}-\Theta_{\mm{b}})-1)e_3 , \\
 & \Delta \bar{\Theta} =0.  \label{201902161355}
 \end{align}
For the sake of simplicity, we consider that
$$\bar{\Theta}= \Theta_{\mm{b}}-\varpi x_3\mbox{ for }0 \leqslant x_3\leqslant h,$$
where $\varpi>0$ is a constant of adverse temperature gradient.

Denoting the perturbation to the magnetic Boussinesq equilibrium state by
$${v}= {v}- {0},\ \theta=\Theta-\bar{\Theta},\  N= M -\bar{M},\ \beta=(p+ \lambda|M|^2/8\pi-\bar{p}- \lambda|\bar{M}|^2/8\pi )/\rho ,$$
then, $( {v},\theta, N,\beta)$ satisfies the perturbation  equations in $\Omega_h$:
\begin{equation}\label{0101oooxx}\left\{{\begin{array}{ll}
 v_t+   v\cdot \nabla v +\nabla \beta  =g \alpha \theta  e_3
  +\nu \Delta v +\lambda (\bar{M}+{N} )\cdot \nabla N /4\pi \rho, \\
  \theta_t  +v\cdot\nabla (\theta+\bar{\Theta}) =\kappa \Delta\theta  , \\[1mm]
N_t +v \cdot \nabla N =(\bar{M}+N )\cdot \nabla v +\sigma \Delta N ,\\
\mm{div}v = \mm{div }N =0 . \end{array}}  \right.
\end{equation}
We call \eqref{0101oooxx} the perturbed magnetic Rayleigh--B\'enard equations (MB equations for short). In particular, we call \eqref{0101oooxx} with $\sigma=0$ the zero-MB equations for simplicity. To investigate the stability of the rest state ${r_\mm{B}}$ for $\sigma\neq 0$,
Galdi considered the following initial-boundary value conditions
\begin{align}
&(v,\theta,N)|_{t=0}=(v^0,\theta^0,N^0),
\label{201903081125}\\
& (v_3, \partial_3 v_{\mm{h}},\theta,N)|_{\partial\Omega}=0 .\label{201903081125xx}
 \end{align}
The initial-boundary value problem of \eqref{0101oooxx}--\eqref{201903081125xx}  with $\sigma\neq0$  is called magnetic Rayleigh--B\'enard problem (MB problem for short).
We mention that, if $\sigma=0$, we can pose the following initial-boundary value conditions for the well-posedness of zero-MB equations:
\begin{align}
&(v,\theta,N)|_{t=0}=(v^0,\theta^0,N^0),
\label{201903081125xxxx}\\
& (v,\vartheta)|_{\partial\Omega}=0.\label{201903081125xxxxxx}
 \end{align}
The initial-boundary value problem of the zero-MB equations,  \eqref{201903081125xxxx} and \eqref{201903081125xxxxxx}  is called  magnetic Rayleigh--B\'enard problem  with zero resistivity (\emph{zero-MB problem} for short). In addition,  the following initial-boundary value problem is called the Rayleigh--B\'enard problem with the first type of boundary value condition (\emph{RB problem} for short)
\begin{equation}\label{0101oooxxxxx}\left\{{\begin{array}{ll}
 v_t+   v\cdot \nabla v +\nabla \beta  =g \alpha \theta  e_3
  +\nu \Delta v , \\
  \theta_t  +v\cdot\nabla (\theta+\bar{\Theta}) =\kappa \Delta\theta  , \\[1mm]
  \mm{div}v=0,\\
  (v,\theta)|_{t=0}=(v^0,\theta^0),\\
  ( v,\theta)|_{\partial\Omega}=0. \end{array}}  \right.
\end{equation}

Motivated by the stability method in \cite{JDDSSTC},
Galdi proved that the MB problem   is exponentially stable in time by (generalized) energy method under the stability conditions \cite{GGPNA}:
\begin{alignat}{2}
&R_\sigma <\frac{1}{2}R_{\mm{s}}+\left(\frac{R_{\mm{s}}^2}{4}+\pi^2R_{\mm{s}} Q_\sigma
\right)^{\frac{1}{2}} & &\
\mbox{ when }P_\mm{m}\leqslant P_\vartheta,\nonumber \\
& R_\sigma <\frac{1}{2}R_{\mm{s}}+\left(\frac{R_{\mm{s}}^2}{4}+\frac{\pi^2R_{\mm{s}}Q_\sigma
}{(P_\mm{m}/P_\vartheta)^2}\right)^{\frac{1}{2}} & &\
\mbox{ when } P_\mm{m}>P_\vartheta,\label{201903072309}
\end{alignat}
where $R_{\mm{s}}\approx 657.5$ is the critical Rayleigh number in the absence of a magnetic field,  and the Rayleigh number $R:=R^2_\sigma$, Chandrasekhar number $Q_\sigma$,  Prandtl number  $P_\mm{\vartheta}$ and magnetic Prandtl numbers $\mm{P_m}$ are given by
$$
R =\frac{g \varpi \alpha h^4}{\nu \kappa},\ Q_\sigma=\frac{\lambda \bar{M}^2_3 h^2}{4\pi \rho \sigma \nu},\ P_\vartheta=\frac{\nu}{\kappa} , \ P_\mm{m}=\frac{\nu}{\sigma}.
$$
 Galdi's result present that the magnetic fields can enlarge the stability condition, and thus mathematically prove the magnetic inhibition phenomenon by magnetic Boussinesq equations with resistivity.

Formally let $\sigma\to 0$ in \eqref{201903072309}, we easily see that the limit state is that
$$R_\sigma<R_{\mm{s}},$$
which is just the stability condition  in the absence of a magnetic field.
Hence, for the case of zero resistivity, we cannot directly observe the conclusion of magnetic inhibition phenomenon from Galdi's result.
However, based on a linearized steady magnetic B\'enard equations with zero resistivity, Chandrasekhar ever asserted  that
\emph{a MHD fluid of zero resistivity should, therefore, be thermally stable for all adverse temperature gradients; moreover, it is clear that a (strong) magnetic field will generally have the effect of inhibiting the onset of convection} in 1956 (see p. 160 in Chandrasekhar's monograph).

An alternative question arises whether we can mathematically verify the inhibition of other flow instabilities by a magnetic field in a non-resistive MHD fluid. The answer is positive, for example,  the inhibition of Rayleigh--Taylor (RT) instability by magnetic fields in non-resistive MHD fluids was proved based on a three-layers energy method \cite{JFJSJMFMOSERT,WYJSARMA}, motivated by Guo--Tice's work in a viscus surface-wave problem \cite{GYTIAED}. Later, the authors of this paper further gave the mechanism of magnetic inhibition phenomenon of non-resistive MHD fluids \cite{JFJSOMARMA}. More precisely, a non-resistive MHD fluid in a rest state in a layer domain with a no-slip boundary condition of velocity can be regarded to be made up by infinite elastic strings (due to magnetic tension) with fixed points along the impressed field. Thus all bent element lines vibrate around their initial locations. In particular, the perturbed MHD fluid will revert the original rest state, if the fluid has viscosity.

In this paper, motivated by the three-layers energy method and the magnetic inhibition mechanism, we establish a stability criterion (see Theorem \ref{thm3origxxx}), under which the solution of zero-MB problem is algebraically stable in time. Moreover, we find out that the stability criterion can be satisfied, if
\begin{equation}
\label{201903111923}
Q > {R^2(1+4{P_\vartheta}^{-1})^2 }/8
\end{equation}
(see Section \ref{201903102308} for derivation), where the Chandrasekhar number (with zero resistivity) $Q$ is defined by
$Q =\lambda \bar{M}_3^2 h^2/ 4\pi\rho \nu^2$.

It should be noted that there exists $R_0$ such that if  the convection condition $R>R_0$ is satisfied, the  RB problem is unstable (referring to \cite{GYHYQCR} or the proofs of Theorem 1.1 in \cite{JFZYYOAR} and Theorem \ref{thmsdfa3}); if  $R<R_0$, the RB problem is stable (referring to \cite{JDDNARAMA}), see   \eqref{201903101912} for the definition of $R_0$. Thus, by \eqref{201903111923}, we see that the strong magnetic filed can inhibit the thermal instability. Therefore, our result, together with Galdi's result, successively completes the mathematical proof of the inhibition of thermal instability by magnetic fields for $\sigma\geqslant 0$. It should be noted the  mechanisms of magnetic inhibition in the both cases of $\sigma=0$ and $\sigma>0$ are different, which will be briefly discussed in Section \ref{201903111229}.

To the best of our knowledge, it is still open problem that whether the MB problem is (nonlinearly) unstable if Galdi's stability conditions fail. However, we can prove that  the zero-MB problem is unstable under some instability criterion based on bootstrap instability method by further developing new analytical skills, see Theorem \ref{201806dsaff012301}. In particular, we find out that the instability stability criterion can be satisfied, if  \begin{equation}
  \sqrt{Q}\leqslant  \min\left\{1, \frac{2(R-R_0)\xi}{ \mathcal{H} +2\sqrt{  \mathcal{H} }},\frac{2(R-R_0)\xi}{1+\sqrt{ \mathcal{H}}}\right\},  \label{201903071510xx}
 \end{equation}
 see Remark \ref{201903131454}, where $R$ satisfies the convection condition $R>R_0$, $\mathcal{H}:=  R/\sqrt{ P_\vartheta}-2(R-R_0)\xi $ and the parameter $\xi$  depends on $\Omega:=\Omega_1$ and $P_\vartheta$ (see \eqref{2019003101953} for definition). Our instability result presents that a small magnetic field cannot inhibit the thermal instability, if the zero-MB problem satisfies convection condition. In particular, we see from \eqref{201903071510xx} that $Q\to 0$ as $R\to R_0$ for given $P_\vartheta$.

We mention that the estimates \eqref{201903111923} and \eqref{201903071510xx} can improved by more careful analysis. Our results of stability and instability can be generalized to a general constant impressive magnetic field $\bar{M}$, or a general  temperature profile $\bar{\Theta}$, which depends on $x_3$, and satisfy \eqref{201903111042}, \eqref{201902161355} and the B\'enard condition
\begin{equation} \label{201705021821nn}
\bar{\Theta}'|_{x_3=x_3^0}<0 \mbox{ for some }x^0_{3}\in (0,h).  \end{equation}
Of course, the stability and  instability criteria should be modified for the both cases of the general $\bar{M}$ and the general $\bar{\Theta}$.

The rest of this paper are organized as follows. In Section \ref{201903102255}, we rewrite the zero-MB problem in Lagrangian coordinates, then state the results of  stability and instability in Lagrangian coordinates, see Theorems \ref{thm3} and \ref{thmsdfa3}, resp., and finally we introduce  the corresponding  results in  Eulerian coordinates, see Theorems \ref{thm3origxxx} and \ref{201806dsaff012301}. In Sections \ref{lowerorder} and \ref{sec:201903102311}, we give the proofs of  Theorems \ref{thm3} and  \ref{thmsdfa3}, resp..
%%%%%%%%%%%%%%%%%%%%%%%%%%%%%%%%%%%%%%%%%%%
\section{Main results}\label{201903102255}
\subsection{Reformulation}
In general, it is difficult to directly investigate the stability and instability of zero-MB problem. Similarly to \cite{JFJSJMFMOSERT,WYJSARMA}, we switch our analysis to Lagrangian
coordinates. To this end, we assume that there is an invertible mapping
$\zeta^0:=\zeta^0(y):\Omega_h\to\Omega_h$, such that
\begin{equation}
\label{zeta0inta}
\partial\Omega=\zeta^0(\partial\Omega_h)\mbox{ and } \det\nabla\zeta^0= 1,
\end{equation}
where $\zeta_3^0$ denotes the third component of $\zeta^0$ and $\det$ the determinant of matrix.
Then we define a flow map $\zeta$ as the solution to
\begin{equation}\label{0114}
\left\{
            \begin{array}{l}
 \zeta_t (y,t)=v(\zeta(y,t),t)
\\
\zeta (y,0)= {\zeta_0}.
                  \end{array}    \right.
\end{equation}
Denote the Eulerian coordinates by $(x,t)$ with $x=\zeta(y,t)$, whereas  $(y,t)\in \Omega_h\times\mathbb{R}^+$ stand for the
Lagrangian coordinates. In order to switch back and forth from Lagrangian to Eulerian coordinates, we assume that $\zeta(\cdot,t)$ is invertible and
$\Omega_h=\zeta(\Omega_h, t)$. In other words, the Eulerian domain of the fluid is the image of $\Omega_h $
under the mapping $\zeta$. In view of the non-slip boundary condition $v|_{\partial\Omega_h}=0$,
we have $\partial\Omega_h=\zeta(\partial\Omega_h, t)$. In addition, since $\det \nabla \zeta_0=1$, the $\zeta$ satisfies the volume-preserving condition
\begin{equation}
\label{zetanabla}
\det \nabla \zeta =1
\end{equation}
due to $\mm{div}v=0$, see \cite[Proposition 1.4]{MAJBAL}.

Now, we further define the Lagrangian unknowns by
\begin{equation}\label{0115}
( u,\vartheta,B,q)(y,t)=( v,\theta,M,\beta)(\zeta(y,t),t) \mbox{ for } (y,t)\in \Omega_h\times\mathbb{R}_+ ,
\end{equation}
where $\mathbb{R}_+:=(0,\infty)$, and $(v,\theta,M,\beta)$ is the solution of the zero-MB problem.
Thus in Lagrangian
coordinates the evolution equations for  $u$, $\vartheta$, $B$ and $q$ read as
\begin{equation}\label{0116} \left\{
                              \begin{array}{ll}
\zeta_t=u, \\[1mm]
  u_t-\nu \Delta_{\ml{A}}u+\nabla_{\ml{A}}q=
  g\alpha \vartheta e_3+
 \lambda B\cdot \nabla_{\ml{A}} B/4\pi\rho, \\[1mm]
\vartheta_t =\kappa \Delta_{\mathcal{A}}\vartheta+\varpi u\cdot \nabla_{\mathcal{A}} \zeta_3,
 \\[1mm]
B_t-B\cdot\nabla_{\ml{A}}u=0, \\[1mm]
\div_\ml{A}B=0, \\[1mm]
\div_\ml{A}u=0
\end{array}
                            \right.
\end{equation}
with initial-boundary value conditions
$$(u,\zeta-y)|_{\partial\Omega_h}=0 \mbox{ and } (\zeta , u, B)|_{t=0}=(\zeta^0, u^0, B^0).
$$Moreover, $\div_\ml{A}B=0$ automatically holds, if the initial data $\zeta^0$ and $B^0$
satisfy
\begin{equation}
\label{diveA00}
\div_{\ml{A}^0}B^0=0.
\end{equation}
Here $\mathcal{A}^0$ denotes the initial value of $\mathcal{A}$,  the matrix $\mathcal{A}:=(\ml{A}_{ij})_{3\times 3}$ is defined via
$\ml{A}^{\mm{T}}=(\nabla \zeta)^{-1}:=(\partial_j \zeta_i)^{-1}_{3\times 3}$, and  the subscript $\mm{T}$ represents the transposition.
The differential operators $\nabla_{\ml{A}}$, $\mm{div}\ml{A}$ and $\Delta_{\ml{A}}$ are defined by $\nabla_{\ml{A}}f:=(\ml{A}_{1k}\partial_kf,
\ml{A}_{2k}\partial_kf,\ml{A}_{3k}\partial_kf)^{\mm{T}}$, $\mm{div}_{\ml{A}}(X_1,X_2,X_3)^{\mm{T}}:=\ml{A}_{lk}\partial_k X_l$ and $\Delta_{\mathcal{A}}f:=\mm{div}_{\ml{A}}\nabla_{\ml{A}}f$
for appropriate $f$ and $X$. It   should be noted that  we have used the Einstein convention of summation over repeated indices, and $\partial_k:=\partial_{y_k}$.

We can  derive from \eqref{0116}$_4$ a differential version of magnetic flux conservation in Lagrangian coordinates \cite{JFJSOMARMA}:
\begin{equation}
\label{abjlj0}
\ml{A}_{jl}B_j=\ml{A}_{jl}^0B_j^0,
 \end{equation}
 which yields
\begin{eqnarray}
\label{0124}   B=\nabla\zeta \ml{A}^{\mm{T}}_0 B_0.
\end{eqnarray}
Here and in what follows, the notation $f_0$ also denotes the initial data of the function $f$. If we assume that
\begin{equation}
\label{201903081437}
\ml{A}^{\mm{T}}_0 B^0=\bar{M}
 \mbox{ (i.e., }B^0= \partial_{\bar{M}}\zeta^0\mbox{)},
 \end{equation} then
 \eqref{0124} reduces to
 \begin{eqnarray}
\label{0124xx}   B=\partial_{\bar{M}}\zeta  .
\end{eqnarray}
It should be noted that \eqref{0116}$_4$ and \eqref{0116}$_5$ automatically holds for $B$ given by \eqref{0124xx}.
Moreover, by \eqref{0124xx}, the magnetic tension can be represented by
\begin{equation}\label{0131}
\begin{aligned}B\cdot \nabla_{\ml{A}} B=& \partial_{\bar{M}}^2\zeta.
\end{aligned}\end{equation}

Thus we easily  see that, under  the initial condition \eqref{201903081437},
 the equations \eqref{0116} can be changed  to a Navier--Stokes system with a magnetic tension and sinking-buoyancy force by using the relation \eqref{0131}:
\begin{equation}\label{peturbationequation}
\begin{cases}
 \eta_t=  u, \\[1mm]
 u_{t}-\nu\Delta_{\ml{A}}  u+\nabla_{\ml{A}}  q= \lambda \partial_{\bar{M}}^2\eta/4\pi\rho +g \alpha  \vartheta e_3
,\\[1mm]
\vartheta_t   =\kappa \Delta_{\mathcal{A}}\vartheta +\varpi u\cdot \nabla_{\mathcal{A}} \zeta_3, \\[1mm]
\div_{\ml{A}} u  =0,
\end{cases}
\end{equation}
where $\mathcal{A}=(I+\nabla \eta)^{-1}$ and $I=(\delta_{ij})_{3\times 3}$.
The  associated  initial and boundary value conditions read as follows
\begin{equation}\label{defineditcon}
(\eta,u,\vartheta)|_{t=0}=(\eta^0,u^0,\vartheta^0),\ (\eta,u,\vartheta)|_{\partial\Omega_h}=0.
\end{equation}
It should be noted that the shift function $q$ is the sum of the perturbation fluid and magnetic pressures in Lagrangian coordinates.
However we still call $q$ the perturbation pressure for the sake of simplicity.

Now we use dimensionless variables
\begin{align}
& y^*= y/h,\ t^*=\nu t/h^2,\ \eta^*=\eta/h, \ u^*= h u/\nu,\  \theta^*=R \kappa \theta/ h \varpi \nu \mbox{ and }q^*= h^2q/\nu^2,  \nonumber
\end{align}
  to rewrite \eqref{peturbationequation}  as  a non-dimensional formed defined in $\Omega$, and then omitting the superscript $*$, we obtain a so-called (non-dimensional) \emph{transformed  MB problem} defined in $\Omega$:
\begin{equation}\label{0101xxxx}\left\{\begin{array}{ll}
 \eta_t=  u, \\[1mm]
 u_{t}- \Delta_{\ml{A}}  u+\nabla_{\ml{A}}  q= Q \partial_{3}^2\eta +R \vartheta e_3
,\\[1mm]
P_{\vartheta} \vartheta_t -\Delta_{\mathcal{A}}\vartheta =R u\cdot \nabla_\mathcal{A} \zeta_3, \\[1mm]
\div_{\ml{A}} u  =0,
 \end{array}\right.\end{equation}
%\lambda h^2/4\pi\rho \nu^2 , & \hbox{ if }\bar{M}_3=0,
% and $$\widetilde{M}:=\left\{                \begin{array}{ll}          (\bar{M}_1,\bar{M}_2,\bar{M}_3)/\bar{M}_3 & \hbox{ for }\bar{M}_3\neq 0;  \\     \bar{M} & \hbox{ for }\bar{M}_3= 0.                \end{array}              \right.$$
The  associated  initial and boundary conditions read as follows
\begin{equation}\label{defineditconxx}
(\eta,u,\vartheta)|_{t=0}=(\eta^0,u^0,\vartheta^0),\ (\eta,u,\vartheta)|_{\partial\Omega}=0.
\end{equation}

\subsection{Simplified notations}

Before stating our main results for the transformed MB problem in details, we shall introduce some simplified notations used throughout this paper.

(1) Basic notations:
 \begin{align}& \mathbb{R}^+_0:=\overline{\mathbb{R}_+},\ I_T:=(0,T),\
L^p:=L^p (\Omega):=W^{0,p}(\Omega)\mbox{ for }1< p\leqslant \infty,\nonumber \\
& {H}^1_0:=W^{1,2}_0(\Omega),\ {H}^k:=W^{k,2}(\Omega ),\ {H}^k_0:=H^{1}_0\cap H^k,\ H^k_\sigma:=H^k_0 \cap\{\mm{div}u=0\},\nonumber \\
& H^{k,1}_0:=H^k_0\cap \{\det(I+\nabla u)=1\},\ \underline{H}^{k}:=H^k\cap \left\{\int u\mm{d}y=0\right\}\mbox{ with }\int:=\int_\Omega,\nonumber \\
& H^{\infty}_0:=\cap_{k=1}^\infty H^k_0,\ H^{\infty}_\sigma:=\cap_{k=1}^\infty H^k_\sigma,\ \underline{H}^{\infty}:=\cap_{k=1}^\infty \underline{H}^k,\ \|\cdot \|_k :=\|\cdot \|_{H^k(\Omega)}\mbox{ for }k\geqslant 1,\nonumber \\
&\Delta_{\mm{h}}:=\partial_1^2+\partial_2^2,\ \partial_{\mm{h}}^k\mbox{ deontes } \partial_{ 1}^{o_1}\partial_{ 2}^{o_2}\mbox{ for any } o_1+o_2=k,\nonumber \\
 & H_0^{1,i}:=\{ u\in  H_0^{1}~|~ \partial_{ 1}^{o_1}\partial_{ 2}^{o_2} u\in L^2\mbox{ for any } o_1+o_2\leqslant i \} ,\ H_\sigma^{1,i}:=H_\sigma^{1}\cap H_0^{1,i}, \nonumber \\
&\|\cdot\|_{m,k}:=\sum_{o_1+o_2=m}
\|\partial_{1}^{o_1}\partial_{2}^{o_2}\cdot\|_{k}^2,\  a\lesssim b\;\mbox{ means that }\; a\leqslant cb.\nonumber  \end{align}
Here and in what follows the letter $c$ denotes a generic constant which may depend on the domain $\Omega$ and the known physical parameters (in the transformed MB problem).

(2) Auxiliary parameters:
\begin{align}
&\mathcal{R}:=R/\sqrt{P_\vartheta}-\Lambda_0 ,\  \frac{1}{R_0}:=\sup_{(\varpi,\phi)\in H_\sigma^1\times H_0^1}\frac{2  \int \varpi_3\phi\mathrm{d}y}{ \|\nabla (\varpi, \phi)\|_0^2},  \label{201903101912}  \\
 & \Lambda_0:=\sup_{(\varpi,\phi)\in \mathcal{A} }
 D_R(\varpi,\phi,1) , \quad\xi:=\sup_{(\varpi,\phi)\in\mathcal{B} }\left\{ \int \varpi_3\phi \mm{d}y\right\} ,\label{2019003101953}
\end{align}
where $ \mathcal{A}:=\{(\varpi,\phi)\in H_\sigma^1\times H^1_0~|~ J(\varpi,\phi) =1\}$, $ J(\varpi,\phi) :=\|\varpi\|_0^2  +P_\vartheta\|\phi \|_0^2$ and
$$\mathcal{B}:=\left\{(\varpi,\phi)\in\mathcal{A}~\left|~ \frac{1}{R_0}=\frac{2  \int \varpi_3\phi\mathrm{d}y}{ \|\nabla (\varpi, \phi)\|_0^2}\right.\right\}.$$

(3) Functionals for  the derivation of  stability result:
\begin{equation*}
  \begin{aligned}
&\mathcal{E}^{L}:=\| \nabla \eta\|_{3,0}^2+\|(\eta,u,\theta)\|_{3}^2+\|(u_t,\vartheta_t,\nabla q)\|_{1}^2,\\
&\mathcal{D}^{L}:=\|(\partial_{3} \eta, \nabla u)\|^2_{3,0}+
\|(\eta,u)\|_{3}^2+\|\theta\|_4^2+\sum_{k=1}^2\|\partial_t^k(u,\theta)\|_{4-2k}^2+\|\nabla q\|_{1}^2+\|\nabla q_t\|_{0}^2,\\
&\mathcal{E}^{H}:=\|  \nabla \eta\|^2_{6,0}+\|\eta\|_{6}^2 +\sum_{k=0}^3 \|\partial_t^k (u,\theta)\|_{6-2k}^2+\sum_{k=0}^2
\|\nabla \partial_t^k q\|_{4-2k}^2, \\
&
\mathcal{D}^{H}:=   \|(\partial_{3} \eta,\nabla u)\|^2_{6,0}+\|(\eta, u )\|_{6}^2
+\|\theta\|_7^2\\
&\qquad \quad +\sum_{k=1}^3\|\partial_t^k (u,\theta)\|_{7-2k}^2+\sum_{k=1}^2\|\nabla \partial_t^kq\|_{5-2k}^2+\|\nabla q\|_{4}^2,
 \\
 &\mathcal{G}_1(t)=\sup_{0\leqslant \tau< t}\|\eta(\tau)\|_7^2,\qquad
 \mathcal{G}_2(t)=\int_0^t\frac{\|(\eta,u )(\tau)\|_7^2}{(1+\tau)^{3/2}}
\mm{d}\tau,\\
& \mathcal{G}_3(t)=\sup_{0\leqslant \tau< t}\mathcal{E}^H(\tau)+\int_0^t
\mathcal{D}^H(\tau)\mm{d}\tau,\qquad \mathcal{G}_4(t)=\sup_{0\leqslant \tau< t}(1+\tau)^{3}\mathcal{E}^L(\tau),
\end{aligned}
\end{equation*}

(4)  Functionals for the derivation of  instability result:
\begin{align}
&D_{R}(u,\vartheta,\tau) =2R\int u_3\vartheta\mathrm{d}y - \tau\|\nabla u\|_{0}^2-\|\nabla \vartheta\|_{0}^2,\nonumber \\
 &F( u,\vartheta,s,\tau):=D_R(u,\vartheta,\tau)-Q \|\partial_3 u\|_0^2 /s,\nonumber \\
& \mathfrak{E} :=\|\eta \|_5^2+\|(u,\vartheta) \|_4^2+ \|\nabla q\|_{2}^2+\|
(u,\vartheta)_t\|_2^2+\|\nabla q_t\|_{0}^2+\|(u,\vartheta)_{tt}\|_0^2 ,\nonumber
\\ &  \mathfrak{D}:= \|(\eta,u,\vartheta)\|_5^2+\|\nabla q\|_{3}^2+\|(u,\vartheta)_t\|_3^2+\|\nabla q_t\|_{1}^2+\|(u,\vartheta)_{tt}\|_1^2,  \nonumber\end{align}

(5) Functionals for stability criterion:
\begin{equation*}
\begin{aligned}
& E_{1}(\varphi,\phi):=Q \|\partial_{3}
\varphi\|^2_{0}+ \frac{2 R^2}{P_\vartheta} \|  \varphi_3\|^2_0+P_{\vartheta} \|   \phi \|_{0}^2 -  4R  \int \varphi_3  \phi \mm{d}y,\\
&\tilde{\mathcal{E}}_{1,1}(\varphi,\psi,\phi):=\frac{1}{2}\|\nabla \varphi\|_{\underline{3},0}^2+ \sum_{0\leqslant |\alpha|\leqslant 2}  E_1(\partial^{\alpha}_{\mm{h}} \varphi,  \partial_{\mm{h}}^j \phi)+
\|   \psi\|^2_{\underline{2},0}
 ,\\
 &\tilde{\mathcal{E}}_1(\varphi,\psi,\phi ):=\tilde{\mathcal{E}}_{1,1}(\varphi,\psi,\phi)+
\sum_{0\leqslant |\alpha|\leqslant 3}
 \int  \partial^{\alpha}_{\mm{h}}\varphi\cdot \partial^{\alpha}_{\mm{h}} \psi\mm{d}y
, \\
 &  \tilde{{\mathcal{D}}}_{1,1}(\varphi,\psi,\phi ):=
Q\|\partial_{3} \varphi\|_{\underline{3},0}^2
 +2  \sum_{ 0\leqslant |\alpha|\leqslant 2}\| \nabla \partial_{\mm{h}}^\alpha (\psi,\phi)\|^2_{0} ,
\\
&\tilde{{\mathcal{D}}}_{1,2}(\varphi,\psi,\phi):=
\|\psi\|_{\underline{3},0}^2+\frac{4 R}{P_\vartheta}\sum_{0\leqslant |\alpha|\leqslant 2}\int\nabla \partial_{\mm{h}}^\alpha\varphi_3\cdot
\nabla\partial_{\mm{h}}^\alpha \phi\mm{d}y \\
&\qquad\qquad\qquad\quad  +R\sum_{0\leqslant |\alpha|\leqslant 3}\int \partial_{\mm{h}}^\alpha\varphi_3\cdot
\partial_{\mm{h}}^\alpha \phi\mm{d}y ,\\
&
 \tilde{{\mathcal{D}}}_1(\varphi,\psi,\phi):= \tilde{{\mathcal{D}}}_{1,1}(\varphi,\psi,\phi )
-\tilde{{\mathcal{D}}}_{1,2}(\varphi,\psi,\phi).
\end{aligned}
\end{equation*}

(6)  Variation forms for stability criterion:
\begin{align}
& \Upsilon_1:= \sup_{(\varphi,\phi)\in H^1_\sigma\times L^2}
 \frac{ 4R  \int  \varphi_3  \phi \mm{d}y }{
Q \|\partial_{3}
\varphi\|^2_{0} +P_{\vartheta} \| \phi\|_{0}^2+ {2 R^2}\|  \varphi_3\|^2_0/P_\vartheta
}, \label{201902282142} \\
  &\Upsilon_2:= \sup_{(\varphi,\psi,\phi)\in H^4_\sigma\times H^3_\sigma\times H_0^3  }\frac{\tilde{\mathcal{D}}_{1,2}(\varphi,\psi,\phi )
 }{\tilde{\mathcal{D}}_{1,1}(\varphi,\psi,\phi )} \nonumber .
\end{align}

\subsection{Stability results}\label{201903111229}
Next  we state our first main result, which presents the inhibiting effect of a magnetic field in the transformed MB problem.
\begin{thm}\label{thm3}
Let   $R$, $Q$
and $P_\vartheta$ satisfy
\begin{equation}
\label{201903032131}
\Upsilon_1,\ \Upsilon_2<1.
\end{equation}
Then there is a sufficiently small constant $\delta >0$, such that for any
$(\eta^0, u^0, \vartheta^0)\in H^{7,1}_0\times H^6_0\times H^6_0$ satisfying the following conditions:
\begin{enumerate}
  \item[(1)] $ \| \eta^0 \|_7^2+\|(u^0,\vartheta^0)\|_{6}^2 \leqslant \delta $;
  \item[(2)] $\zeta^0:=y+\eta^0$ satisfies \eqref{zeta0inta};
  \item[(3)]  $(\eta^0, u^0, \vartheta^0)$ satisfies necessary compatibility conditions (i.e., $\partial_t^i u(y,0)=0$ and $\partial_t^i \vartheta(y,0)=0$ on $\partial\Omega$ for $i=1$ and $2$),
\end{enumerate}
 there exists a unique global classical solution $(\eta,u,\vartheta)\in C^0(\mathbb{R}^+_0,H^{7,1}_0\times H^6\times H^6)$
to the transformed MB problem \eqref{0101xxxx}--\eqref{defineditconxx} with an associated  perturbation pressure  $q$.
Moreover, $(\eta,u,\vartheta,q)$ enjoys the following stability estimate:
\begin{equation}\label{1.19}
  \mathcal{G}(\infty):=\sum_{k=1}^4\mathcal{G}_k(\infty) \leqslant  c (\| \eta^0 \|_7^2+\|(u^0,\vartheta^0)\|_{6}^2).
\end{equation}
   Here the positive constants $\delta$ and  $c$  depend on the domain $\Omega$ and other known physical parameters.
\end{thm}
\begin{rem}\label{201903110946}
For $k\geqslant 2$, let
 \begin{align}& {H}_{0,*}^{k,1}:=\{  \varpi\in H^{k,1}_0 ~|~  \phi (y): =\varpi+y : \overline{\Omega } \mapsto  \overline{\Omega  } \mbox{ is a homeomorphism mapping},\nonumber \\
 &\qquad \qquad\qquad \qquad \quad  \phi  (y): \Omega \mapsto   \Omega  \mbox{ is a } C^{k-2}\mbox{-diffeomorphic  mapping}\}.\nonumber
\end{align}
  Since $\eta=0$ on boundary, then, if $\delta$ in Theorem 2.1 is sufficiently small, we can further have $\eta\in {H}_{0,*}^{k,1}$ (referring to Lemma 4.2 in \cite{JFJSOMARMA} for derivation).  This means the obtained  stability result in Lagrangian coordinates can be recovered to the one in Eulerian coordinates.
\end{rem}
\begin{rem}
We easily see from the proof of Theorem \ref{thm3} that the transformed MB problem always stable for any $R<R_0$ and any $Q>0$. In the other words, Theorem \ref{thm3} holds with the both conditions of $R<R_0$ and $Q>0$ in place of the stability conditions in \eqref{201903032131}.
\end{rem}

The  proof of Theorem \ref{thm3} is based on the ideas of multi-layer energy method and magnetic inhibition mechanism. Next we briefly introduce the proof idea.

 Let's first observe the basic energy identity of the transformed MB problem:
\begin{align}
&\frac{1}{2}\frac{\mm{d}}{\mm{d}t}\left(Q \|\partial_3 \eta\|_{0}^2
+\| u\|^2_0+P_{\vartheta}  \|   \vartheta \|_{ 0}^2\right)+  \|\nabla (u,\vartheta)\|_{0}^2\nonumber \\
&=2R \int u_3 \vartheta\mm{d}y + \mbox{ int (i.e., integrals involving nonlinear terms)},\label{201903111644}
\end{align}
where the first integral term on the right hand of \eqref{201903111644} is called thermal instability integral.
By magnetic inhibition theory in non-resistive MHD fluids in \cite{JFJSOMARMA}, the RT instability can inhibited by the magnetic tension for strong magnetic fields. The mathematical relation to describe this physical phenomenon is that the (perturbation) RT instability integral (or perturbation gravity potential energy) can be controlled by the (perturbation) magnetic energy. Motivated by this idea, we naturally use \eqref{0101xxxx}$_1$ to write \eqref{201903111644} as follows
\begin{align}
&\frac{1}{2}\frac{\mm{d}}{\mm{d}t}\left(Q \|\partial_3 \eta\|_{0}^2
+P_{\vartheta} \|   \vartheta \|_{ 0}^2+\|  u\|^2_0-4R \int \eta_3 \vartheta \mm{d}y \right)+  \|\nabla (u,\vartheta)\|_{0}^2\nonumber \\
& =-2R \int \eta_3 \vartheta_t\mm{d}y + \mbox{ int}.
\end{align}

To deal with $\vartheta_t$, the unique method is to put the perturbation temperature equation \eqref{0101xxxx}$_2$ into the above identity, and thus we get
\begin{align}
 \frac{1}{2}\frac{\mm{d}}{\mm{d}t}\left(E_1(\eta, \vartheta )+\| u\|^2_0\right)+  \|\nabla (u,\vartheta)\|_{0}^2  =\frac{2R}{P_\vartheta} \int   \nabla  \eta_3\cdot
\nabla  \vartheta \mm{d}y + \mbox{int}.
\label{201903112042}
\end{align}
 To control the first integral on the right hand of \eqref{201903112042}, we shall derive from \eqref{0101xxxx}$_1$ and \eqref{0101xxxx}$_2$ that, for $0\leqslant i \leqslant 1$,
\begin{align}
 \frac{\mm{d}}{\mm{d}t}\left(\int  \partial_{\mm{h}}^i \eta\cdot \partial_{\mm{h}}^i  u\mm{d}y +\frac{1}{2}\|\nabla \partial_{\mm{h}}^i \eta\|_{0}^2\right)
+ Q \|\partial_{3} \partial_{\mm{h}}^i\eta\|_{0}^2
  = \| \partial_{\mm{h}}^i   u\|^2_{0}+  R\int\partial_{\mm{h}}^j\eta_3\partial_{\mm{h}}^j \vartheta\mm{d}y  + \mbox{int}. \nonumber
\end{align}
Thus we further derive from the above two identities that
\begin{align}
 \frac{\mm{d}}{\mm{d}t}\tilde{\mathcal{E}}_2( \eta, u,\vartheta)
+ \tilde{\mathcal{D}}_2( \eta, u,\vartheta)
  =\mbox{int}.,\label{201903122102}
\end{align}
where
$$\begin{aligned}
 &\tilde{\mathcal{E}}_2( \eta, u,\vartheta):= \frac{1}{2} (2( E_1(\eta, \vartheta )+\| u\|^2_0 ) + \|\nabla  \eta\|_{\underline{1},0}^2)+ \sum_{0\leqslant |\alpha|\leqslant 1}\int \partial_{\mm{h}}^\alpha \eta\cdot \partial_{\mm{h}}^\alpha  u\mm{d}y,\\
 & \tilde{\mathcal{D}}_2( \eta, u,\vartheta)=2 \|\nabla (u,\vartheta)\|_{0}^2 +Q \|\partial_{3} \eta\|_{\underline{1}, 0}^2 - \|   u\|^2_{\underline{1},0}\nonumber \\
&\qquad \qquad \qquad- \frac{4R }{P_\vartheta} \int   \nabla  \eta_3\cdot
\nabla  \vartheta \mm{d}y
-  R  \sum_{0\leqslant |\alpha|\leqslant 1}\int\partial_{\mm{h}}^\alpha\eta_3\partial_{\mm{h}}^\alpha \vartheta\mm{d}y.
\end{aligned} $$

Exploiting some analytic tools, we can prove that, for sufficiently large $Q$, the  magnetic energy $\|\partial_3\eta\|_0^2$ can absorb the non-positivity terms such as $-4R \int \eta_3 \vartheta \mm{d}y$, $- \|   u\|^2_{\underline{1},0}$ and so on, and thus $\tilde{\mathcal{E}}_2( \eta, u,\vartheta)$ and $\tilde{\mathcal{D}}_2( \eta, u,\vartheta)$ are equivalent to $\|(\eta,u,\vartheta)\|_0^2+\|\nabla \eta  \|_{\underline{1},0}^2$ and $\|(u,\vartheta)\|_1^2+\|\partial_3 \eta\|_{\underline{1},0}^2$, resp.. Of course, to make sure the equivalence, it suffices to propose proper stability conditions, which are similar to \eqref{201903032131}.

\emph{From the above mathematical analysis, we easily see that the magnetic energy controls the thermal instability integral in non-resistive MHD fluids by an indirect manner. This mathematical mechanism is different to the ones in the known magnetic inhibition phenomenons. In fact, in the inhibition of RT instability, the magnetic energy directly controls the RT instability integral caused gravity (i.e., magnetic tension directly resists gravity) \cite{JFJSOMARMA}. In the inhibition of thermal instability in the resistive MHD fluids, due to the viscosity and heat conduction, the thermal instability integral is also directly controlled dissipative term caused by magnetic energy
\begin{equation}
\label{201903122243}
-(Q_\sigma\varsigma P_\vartheta/R P_m)\int |\nabla \partial_3 \theta|^2\mm{d}y
\end{equation}
(see variational form of stability criterion (4.2) in \cite{GGPNA}),
where $\varsigma\in (0,1)$ is a parameter.
The  dissipative term \eqref{201903122243} reflects that the magnetic energy will convert the thermal energy, when the magnetic fields  are diffusing in the resistive MHD fluids.
}

By the basic idea in \eqref{201903122102}, we can further derive from the transformed MB problem that there are two functionals $\tilde{\mathcal{E}}^{L}$ and $Q$ of $(\eta,\vartheta,u)$ satisfying the lower-order energy inequality (see  Proposition \ref{pro:0301n})
\be
\label{latdervin}
  \frac{\mm{d}}{\mm{d}t}\tilde{\mathcal{E}}^{L}+\mathcal{D}^L\leqslant \mathcal{Q}\mathcal{D}^L,\ee
where    the functional $\tilde{\mathcal{E}}^{L}$ is equivalent to $\mathcal{E}^L$ under  the stability condition \eqref{201903032131}.
Unfortunately, we can
not close the energy estimates for the global well-posedness of the transformed MB problem only based on \eqref{latdervin}, since  $\mathcal{Q}$ can not be controlled by $\tilde{\mathcal{E}}^{L}$.
Guo and Tice also faced similar trouble in the investigation of surface wave problem \cite{GYTIDAP,GYTIAED}.  However, they developed a three-layers energy method to overcome this trouble. Later Tan--Wang also used the three-layers energy method to prove the global well-posedness of the initial-boundary value problem of non-resistive MHD fluids \cite{TZWYJGw,WYJSARMA}.

Motivated by the three-layers energy method, we shall further look after a higher-order energy inequality to match the lower-order energy inequality \eqref{latdervin}.
Since $\tilde{\mathcal{E}}^L$ contains $\|\eta\|_3$,
we  find that the higher-order energy at least includes $\|\eta\|_6$. Thus, similar to \eqref{latdervin},  we shall establish the  higher-order energy inequality  (see Proposition \ref{pro:0301n}):
\be
\label{latdervins}
  \frac{\mm{d}}{\mm{d}t}\tilde{\mathcal{E}}^{H}+ {\mathcal{D}}^H\leqslant \sqrt{\mathcal{E}^{L}}\|(\eta,u)\|_7^2,\ee
where the functional $\tilde{\mathcal{E}}^{H}$ is equivalent to $\mathcal{E}^H$ under  the stability condition. Moreover,
the highest-order norm $\|(\eta,u)\|_7^2$ enjoys the   highest-order energy inequality
\be
\label{latdervsin}\frac{\mm{d}}{\mm{d}t}  \|\eta\|_{7,*}^2+ \|( \eta,u)\|_{7}^2\lesssim \mathcal{E}^H+\mathcal{D}^H,
\ee
where the norm  $\|\eta\|_{7,*} $ is equivalent to $ \|\eta\|_{7} $.
In the derivation of \emph{a priori} estimates,  we have $\mathcal{Q}\lesssim \mathcal{E}^H$, and thus
 \eqref{latdervin} implies (see Proposition \ref{pro:0301n})
 \be   \label{starlatdervin}
  \frac{\mm{d}}{\mm{d}t}\tilde{\mathcal{E}}^{L}+\mathcal{D}^L\leqslant 0.\ee
Consequently, by a standard  three-layers energy method, we easily deduce the global-in-time stability estimate \eqref{1.19} based on \eqref{latdervins}--\eqref{starlatdervin}.

With Theorem  \ref{thm3} and Remark \ref{201903110946} in hand,
we easily obtain the stability result of zero-MB problem by employing the transformation of Lagrangian coordinates (referring to
\cite[Theorem 1.2]{JFJSJMFMOSERT} and \cite[Theorem 1.2]{JFJSOMITNWZ} for the derivation).
\begin{thm}\label{thm3origxxx}
Under the condition \eqref{201903032131},
then there is a sufficiently small constant $\delta_{0}> 0$, such that, for any
$( v^0,\theta^0,N^0)\in H^6_{\sigma,h}\times H^6_{0,h}\times H^6_h$ satisfying the following conditions:
\begin{enumerate}
    \item[(1)]  there exists an invertible mapping $\zeta^0:=\zeta^0(y):\Omega_h\to \Omega_h$, such that
\eqref{zeta0inta} holds, where $\ml{A}^{\mm{T}}_0=(\nabla \zeta^0)^{-1}$;
 \item[(2)] $\mm{div}v^0=0$, $(\bar{M}+N^0)(\zeta^0)=\bar{M}_3\partial_3 \zeta^0 $;
\item[(3)] $ \| \zeta_0 -x\|_{7,h}^2+\|(v^0,\theta^0)\|_{6,h}^2 \leqslant \delta_0$;
  \item[(4)] the initial data $v^0$, $\theta^0$, $N^0$ satisfy necessary compatibility conditions (i.e., $ \partial_t^iv(x,0)=0$ and $\partial_t^i\theta(x,0)=0$ on ${\partial\Omega_h}$ for $i=1$ and $2$),
\end{enumerate}
  there exists a unique global solution $((v,\theta, N),\beta)\in C^0(\mathbb{R}^+_0,H^6_h\times \underline{H}^5_h)$
to the zero-MB problem. Moreover, $(  v,\theta,N,\beta)$ enjoys
 the following  stability estimate:
\begin{equation}\label{1.19st}
\begin{aligned}
 &\sup_{0\leqslant t< \infty}\left(\|N\|_{6,h}^2+\sum_{k=0}^3 \|\partial_t^k ({v},\theta)(t)\|_{6-2k,h}^2+\sum_{k=0}^2
\|\nabla \partial_t^k \beta(t)\|_{4-2k,h}^2\right)\\
&+\sup_{0\leqslant t< \infty}(1+t)^{3}\left( \|( v,\theta,N)\|_{3,h}^2+\|(v_t,\theta_t,\nabla \beta)\|_{1,h}^2\right) \leqslant c ( \|\zeta^0-x \|_{7,h}^2+\|(v^0,\theta^0)\|_{6,h}^2).
\end{aligned}
\end{equation}
 Here the positive constants $\delta_{0}$ and  $c$  depend on the domain $\Omega_h$ and other known physical parameters in zero-MB problem, and $H^k_{\sigma,h}$, $ H^k_{0,h}$, $H^k_h$, $\underline{H}^k_h$ and $\|\cdot\|_{k,h}$ are defined as $H^k_{\sigma}$, $ H^k_{0}$, $H^k$, $\underline{H}^k$ and $\|\cdot\|_{k}$ with $\Omega_h$ in place of $\Omega$, resp.
\end{thm}
\begin{rem}
The conclusion in Theorem \ref{thm3origxxx} can be also extended to the case of the fluid domain being $\mathbb{R}^2\times (0,h)$.
\end{rem}

\subsection{Instability results}
Now we turn to introduce the instability result of the transformed MB problem, which presents that the small magnetic field cannot inhibit the thermal instability.
\begin{thm}\label{thmsdfa3}
Let $R$ satisfy the  convection condition $R>R_0$. Under the instability condition
  \begin{equation}
 \sqrt{Q}\leqslant  \min\left\{1, \frac{\Lambda_0 }{ \mathcal{R} +2\sqrt{ \mathcal{R}}},\frac{\Lambda_0}{1+\sqrt{ \mathcal{R} }}\right\},\label{201903071510}
  \end{equation} the zero solution to the transformed MB problem \eqref{0101xxxx}--\eqref{defineditconxx} is unstable in the Hadamard sense, that is, there are positive constants $\Lambda^*$ (see Lemma \ref{thm:0201201622xx} for the definition $\Lambda^*$), $m_0$, $\epsilon$, $\delta_0$,
and  $(\tilde{\eta}^0,\tilde{u}^0,\tilde{\vartheta}^0,
\tilde{q}^0,\eta^\mm{r},(u^\mm{r},
\vartheta^\mm{r}),q^\mm{r})\in H^5_\sigma\times  H^4_\sigma\times H_0^4\times \underline{H}^3\times H^5_0\times H_0^4\times \underline{H}^3$,
such that, for any $\delta\in (0,\delta_0)$ and the initial data
\begin{equation}
\label{201901092134}
(\eta^0, u^0, \vartheta^0, q^0):=\delta(\tilde{\eta}^0,\tilde{u}^0,\tilde{\vartheta}^0, \tilde{q}^0)  + \delta^2(\eta^{\mm{r}},u^\mm{r}, \vartheta^\mm{r}, q^{\mm{r}})\in H^{5,1}_{0}\times H^4_\sigma\times H^4_0\times \underline{H}^3,
\end{equation}
there is a unique classical solution $(\eta,(u,\vartheta),q)\in C^0(\overline{I_T},H^{5,1}_{0}\times  H^4_0\times \underline{H}^3)$ with initial data $(\eta^0,u^0,\vartheta^0,q^0)$ to the  transformed MB problem  satisfying
\begin{align}
 \|\chi (T^\delta)\|_{0} \geqslant  {\varepsilon}  \label{201901092145xx}
\end{align}
for some escape time $T^\delta:=\frac{1}{\Lambda^*}\mm{ln}\frac{2\epsilon}{m_0\delta}\in I_T$, where $\chi=\eta_{\mm{h}}$, $\eta_3$, $u_{\mm{h}}$, $u_3$ and $\vartheta$ (we can further take  $\chi=\partial_3 \eta_{\mm{h}}$ and $\partial_3 \eta_3$ for $Q\neq0$). Moreover, the initial data $(\eta^0,u^0,\vartheta^0,q^0)  $ satisfies the compatibility  conditions:
\begin{alignat}{2}
&\mm{div}_{\mathcal{A}^0}u^0=0& &\ \mbox{ in }\Omega, \label{20180923221611xx123}  \\
& \mm{div}_{\mathcal{A}^0}\left( u^0 \cdot \nabla_{\mathcal{A}^0} u^0 + \nabla_{\mathcal{A}^0}\sigma^0-\Delta_{\mathcal{A}^0}u^0- Q\partial_3^2\eta^0-  R \vartheta^0 e_3\right)=0 & &\ \mbox{ in }\Omega, \\
&\nabla_{\mathcal{A}^0}\sigma^0- \Delta_{\mathcal{A}^0}u^0 -Q \partial_3^2\eta^0  =0 & &\ \mbox{ on }\partial \Omega, \label{20180923221611xx1234}
 \\
&  \Delta_{\mathcal{A}^0}\vartheta^0 =0 & &\ \mbox{ on }\partial \Omega, \label{201903111314} \end{alignat}
where   $\mathcal{A}^0$  denotes the initial data of $\mathcal{A}$, and is defined by $\eta^0$.
\end{thm}
\begin{rem}
\label{201903131454} It is should be noted that $\Lambda_0>0$ and $2R>\Lambda_0\sqrt{P_\vartheta}$ are positive under the convection condition $R>R_0$, see Remark \ref{201903102310}. In addition, using \eqref{201903092105} in Remark \ref{201903102310}, we find $Q$ must satisfy \eqref{201903071510}, as
 \eqref{201903071510xx} is satisfied.
\end{rem}

The proof of Theorem \ref{thmsdfa3} is based on a  bootstrap instability method, which has its origin in \cite{GYSWIC,GYSWICNonlinea}. Later, various versions of bootstrap approaches were established by many authors, see \cite{FSSWVMNA,GYHCSDDC,FSNPVVNC,FrrVishikM} for examples. Since the instability solutions shall satisfies the compatibility  conditions, we adapt the version of bootstrap instability method in \cite{JFJSOMITNWZ} %\cite{JFJSJMFMNTICF} or \cite[Theorem 2]{JFJSNSIRTPSCMHD}
to prove Theorem \ref{thmsdfa3} by further introducing new ideas.

In view of the bootstrap instability method in \cite{JFJSOMITNWZ}, the proof of Theorem \ref{thmsdfa3} should be divided five steps:
First, we introduce unstable solutions to the linearized MR problem (see Proposition \ref{thm:0201201622}). Second, we use energy method to derive that the local-in-time solutions of the (nonlinear) RT problem enjoy a Gronwall-type energy inequality (see Proposition \ref{pro:0401nxd}). Third we use the existence theories of the both of Stokes and elliptic problems and an iterative technique to modify initial data of solutions of the linearized MB problem, so that the obtained modified initial data  satisfy the compatibility conditions \eqref{20180923221611xx123}--\eqref{201903111314} and the volume-preserving condition, and are close to the original initial data (see Proposition \ref{lem:modfied}). Fourth, we derive the error estimates between the solutions of the linearized and nonlinear (transformed) MB problems (see Proposition \ref{lem:0401xxxxxx}), and thus get Theorem \ref{thmsdfa3}. Finally we prove the existence of escape time $T^\delta$.
Next we shall briefly how to derive the instability condition \eqref{201903071510}, which involves new idea.

First, we adapt the modified variational method of PDE  in \cite{JFJSO2014} to construct the linear unstable solutions. We mention that the modified variational method was probably first used by Guo and Tice to construct unstable solutions to a class of ordinary differential equations (ODEs) arising from a linearized RT
instability problem \cite{GYTI2}. By the modified variational method of PDE, we natural expect that there exists a positive fixed point of the variational curve $\alpha(s,\tau)$ defined on $\mathbb{R}_+$ (see \eqref{201903111430} for the definition), where the parameter $\tau \in (0,1]$ is given. By careful analyzing the properties of $\alpha(s,\tau)$, we find that if
\begin{equation}
\label{201903052059}
\sqrt{Q}\leqslant  {\Lambda_0 }/{2\sqrt{\mathcal{R}}},
\end{equation}
then $\alpha(\Lambda_0/2,\tau)>\Lambda_0/2$. This relation makes sure the existence of a positive fixed point, and thus we prove the instability of  linearized MB problem with parameter  $\tau$ and with a largest growth rate $\Lambda(\tau)$ (see Proposition  \ref{thm:0201201622}). In particular, taking $\tau=1$, we get an unstable solution of linearized MB problem.

Second, in the derivation of the error estimate, we have an error integral identity:
\begin{align}
Q \|\partial_3 \eta^{\mathrm{d}}\|_{0}^2+\|  u^{\mathrm{d}}\|^2_0+{P_\vartheta} \|\vartheta^{\mathrm{d}} \|_0^2 - 2 \int_0^t  D_{R}(u^{\mm{d}},\vartheta^{\mm{d}},1)\mm{d}\tau  = 2 \int_0^t {R}_1(\tau)\mm{d}\tau  +  R_2,\label{0314}
\end{align}
where ${R}_1(\tau)$ is defined in \eqref{201903111523} and $R_2:=Q \|\partial_3 \eta^{\mathrm{d}}_0\|_{0}^2+\|  u^{\mathrm{d}}_0\|^2_0+{P_\vartheta} \|\vartheta_0^{\mathrm{d}}\|_0^2$. If we directly apply the method of largest growth rate to the above identity, then
\begin{align}
Q \|\partial_3 \eta^{\mathrm{d}}\|_{0}^2+\|  u^{\mathrm{d}}\|^2_0+{P_\vartheta} \|\vartheta^{\mathrm{d}} \|_0^2   \leqslant 2\left( \Lambda^*(\| u^{\mathrm{d}}\|^2_0+P_\vartheta\|  \vartheta^{\mathrm{d}}\|_0^2 ) +Q\|\partial_3 u^{\mm{d}}\|_0^2/ \Lambda^*\right)+ c\delta^3 e^{3\Lambda^* t}.\label{0314sfasfdasf}
\end{align}
For $Q=0$, using  Gronwall's lemma (see Lemma 1.1 in \cite{NASII04}), we get from \eqref{0314sfasfdasf} the error estimate
\begin{align}
\label{201903111545}
Q \|\partial_3 \eta^{\mathrm{d}}\|_{0}^2+\|  u^{\mathrm{d}}\|^2_0+{P_\vartheta} \|\vartheta^{\mathrm{d}} \|_0^2 \lesssim  \delta^3 e^{3\Lambda^* t} .
\end{align}

Obviously, the above idea fails for $Q>0$. To overcome this difficulty, we expect there exists a term $\|\partial_3u^{\mm{d}}\|_0^2$ appearing on the left hand of \eqref{0314sfasfdasf}, and thus rewrite \eqref{0314} as follows:
\begin{align}
&
Q \|\partial_3 \eta^{\mathrm{d}}\|_{0}^2+\|  u^{\mathrm{d}}\|^2_0+{P_\vartheta} \|\vartheta^{\mathrm{d}} \|_0^2 +2 \sqrt{Q}\int_0^t\|\nabla  u^{\mathrm{d}} \|_0^2\mm{d}\tau \nonumber \\
& = 2 \int_0^t  D_{R}(u^{\mm{d}},\vartheta^{\mm{d}},1-\sqrt{Q})\mm{d}\tau +2 \int_0^t {R}_1(\tau)\mm{d}\tau  +  R_2,\label{0314cccccc}
\end{align}
where $Q$ shall further satisfy $Q<1$.
Applying  the method of largest growth rate to the above identity  yields that
\begin{align}
&Q \|\partial_3 \eta^{\mathrm{d}}\|_{0}^2+\|  u^{\mathrm{d}}\|^2_0+{P_\vartheta} \|\vartheta^{\mathrm{d}} \|_0^2+ 2  \sqrt{Q} \int_0^t\|\nabla  u^{\mathrm{d}} \|_0^2\mm{d}\tau \nonumber  \\
 &\leqslant  2\Lambda(1-\sqrt{Q} )\int_0^t(\| u^{\mathrm{d}}\|^2_0+P_\vartheta\|  \vartheta^{\mathrm{d}}\|_0^2 ) \mm{d}y+\frac{2Q}{ \Lambda(1-\sqrt{Q} ) }\int_0^t\|\nabla  u^{\mathrm{d}} \|_0^2\mm{d}\tau+c\delta^3 e^{3\Lambda^* t}.\label{201903101526xxxxx}
\end{align}

An important observation is that the second integral (with respect to time) on the right hand of \eqref{201903101526xxxxx} may be absorbed by the first integral on the left hand,  if $Q$ is properly small. Thus, by careful analysis, we find that, if $Q$ further satisfies
\begin{equation}
\label{201903222102xx} \sqrt{Q}\leqslant \Lambda_0/(1+\sqrt{ \mathcal{R} }) ,
 \end{equation}the estimate \eqref{201903101526xxxxx} reduces to
\begin{align}
&Q \|\partial_3 \eta^{\mathrm{d}}\|_{0}^2+\|  u^{\mathrm{d}}\|^2_0+{P_\vartheta} \|\vartheta^{\mathrm{d}} \|_0^2 \nonumber  \\
 &\leqslant  2\Lambda(1-\sqrt{Q} )\int_0^t(\| u^{\mathrm{d}}\|^2_0+P_\vartheta\|  \vartheta^{\mathrm{d}}\|_0^2 ) \mm{d}y +c\delta^3 e^{3\Lambda^* t}.\label{201903101526}
\end{align}
In addition,  if $Q$ further satisfies
\begin{equation}
\label{201903222102}
Q\leqslant \Lambda_0 /({\mathcal{R} +2\sqrt{ \mathcal{R}}}),
\end{equation}  we have \begin{equation}
   \label{201903101525}
   \Lambda^*\geqslant 2\Lambda(1-\sqrt{Q})/3.
    \end{equation}
Thus we deduce the desired error estimate \eqref{201903111545} from \eqref{201903101526} and \eqref{201903101525} under the instability condition \eqref{201903071510}.

Similarly to Theorem \ref{thm3origxxx}, we easily deduce the the thermal instability of the zero-MB problem in Eulerian coordinates under the instability condition \eqref{201903071510} from Theorem \ref{thmsdfa3} and Remark \ref{201903110946}.
\begin{thm}\label{201806dsaff012301}
Under  the instability condition \eqref{201903071510}, the zero solution of the zero-MB problem is unstable in the Hadamard sense, that is, there are positive constants $C$ and $\delta_0$ such that, for any $\delta\in (0,\delta_0)$, there exists $(v^0,\theta^0,N^0,q ^0 )\in  H^4_{\sigma,h} \times H_{0,h}^4 \times H_{\sigma,h}^4 \times \underline{H}^3_h $ satisfying
\begin{equation}
\| (  v^0,\theta^0, N^0)\|_{4,h}+\|q ^0\|_{3,h}\leqslant C\delta,
\end{equation}
and  a unique classical solution $( v,\theta,N,q )\in C^0(\overline{I_T}, H^4_{\sigma,h} \times  H_{0,h}^4\times H^4_{\sigma,h} \times \underline{H}_h^3)$ to the zero-MB problem with  initial data $ ( v^0,\theta^0, N^0,q ^0)  $,
but the solution satisfies
\begin{align}
&\|\chi(T^\delta)\|_{L^1_h}\geqslant  \varepsilon   \label{201809041319xx}
\end{align}
for some escape time $T^\delta \in I_T$, where $\chi=v_{\mm{h}}$, $v_3$ and $\vartheta$ (we can further take $\chi=N_{\mm{h}}$ and $N_3$ for $Q\neq 0$), and
$L^1_h$ is defined as $L^1$ with $\Omega_h$ in place of $\Omega$. Moreover, the initial data $v^0$, $\theta^0$, $N^0$ and $q^0$ satisfy  necessary compatibility conditions.
\end{thm}

\section{Proof for stability of the transformed MB problem}\label{lowerorder}

To obtain Theorem \ref{thm3}, the key step is to derive the stability estimate \eqref{1.19}.  To this end, let $(\eta,u,\theta)$ be a solution
of the transformed MB problem with a perturbation pressure $q$, such that
\begin{equation}\label{aprpioses}
\sqrt{\mathcal{G}_1(T)+ \sup_{0\leqslant \tau\leqslant T}\mathcal{E}^H(\tau)}\leqslant \delta\in (0,1)\;\; \mbox{ for some  }T>0
\end{equation}
and
\begin{equation}
\label{2019022201518}
\det(I+\nabla \eta^0)=1,
\end{equation}
where $\delta$ is sufficiently small. It should be noted that the smallness depends on the known physical parameters, and  will be repeatedly used in what follows. Moreover, we assume that the solution $(\eta,{u},\vartheta, q)$ possesses proper regularity, so that the procedure
of formal calculations makes sense. Next we mention some estimates of $\mathcal{A}$ and $\mm{div}\eta$, which will be repeatedly used in what follows.

Since $\eta^0$ satisfies volume-preserving condition \eqref{2019022201518}, then we also  have
\begin{equation}
\det(I+\nabla \eta)=1.
\end{equation}
Thus
\begin{equation}
\label{matirxA}
\mathcal{A}=(A^*_{ij})_{3\times3}, \end{equation}
where $A^{*}_{ij}$ is the algebraic complement minor of $(i,j)$-th entry in the matrix $I+\nabla\eta$. Moreover, we can compute out $ \partial_kA^{*}_{ik}=0$, and thus get the relation
\begin{equation}\label{diverelation}
\mm{div}_{\mathcal{A}}u=\partial_l (\mathcal{A}_{kl}u_k)=0,
\end{equation} which will be used in the derivation of temporal derivative estimates.

Under the assumptions of \eqref{aprpioses} and \eqref{2019022201518},
 we can easily derive that the  matrix $\mathcal{A}$ enjoys the following estimates (see \cite[Section 2.1]{JFJSJMFMOSERT}):
\begin{align} & \label{aimdse}
\|\mathcal{A}\|_{L^\infty}   \lesssim 1 ,  \\
&\label{prtislsafdsfs}\| \partial_{t}^i\ml{A}\|_j \lesssim \sum_{m=0}^{i-1} \|\partial_t^m  u\|_{j+1},\\
&\label{prtislsafdsfsfds}\|\tilde{\mathcal{A}}\|_k \lesssim \|  \eta\|_{k+1},
\\
& \label{lemdssdf} \| \nabla f\|_{0}\leqslant \sqrt{2} \| \nabla_{\mathcal{A}}  f\|_{0}\mbox{ for any } f\in H^1.
\end{align}
where  $1\leqslant i\leqslant 4$, $0\leqslant j\leqslant 8-2i$ and $0\leqslant k\leqslant 6$.

Since $\det(I+\nabla \eta)=1$, we have by Sarrus' rule that
\begin{align}
\mathrm{div}\eta=&\partial_1\eta_2\partial_2 \eta_1+ \partial_2 \eta_3\partial_3\eta_2+\partial_ 3\eta_1\partial_1\eta_3
-\partial_1\eta_1\partial_2 \eta_2-
\partial_1 \eta_1\partial_3\eta_3-\partial_2\eta_2\partial_3\eta_3 \nonumber \\
&+ \partial_1\eta_1(\partial_2\eta_3\partial_3\eta_2 -\partial_2\eta_2\partial_3\eta_3)+
\partial_2\eta_1(\partial_1\eta_2\partial_3\eta_3 -\partial_1\eta_3\partial_3\eta_2) \nonumber
\\
&+\partial_3\eta_1(\partial_1\eta_3\partial_2\eta_2 -\partial_1\eta_2\partial_2\eta_3). \label{201903211118}
\end{align}
Thus we further estimate that
\begin{equation}
\label{201903011350}
\|\mm{div}\eta\|_{k}\lesssim \left\{
                       \begin{array}{ll}
                    \|\eta\|_{3}\|\eta\|_{k+1} &\hbox{ for }0\leqslant k\leqslant 5,\\
                \|\eta\|_3\|\eta\|_{k+1}+\|\eta\|_{6}\|\eta\|_{k-2}&\hbox{ for }  k=6.
                       \end{array}
                     \right.
\end{equation}We mention that, by the expression \eqref{201903211118} (see \cite[Section 2.3]{JFJSJMFMOSERT}),\begin{equation}
\label{201903082042}
\mm{div}\eta=\mm{div}\Phi,
\end{equation}
where
\begin{equation}
\label{201903082040}
\Phi(\eta):=\left(\begin{array}{c}
         - \eta_1(\partial_2\eta_2+\partial_3\eta_3 )+
           \eta_1(\partial_2\eta_3 \partial_3\eta_2
- \partial_2\eta_2\partial_3\eta_3) \\
      \eta_1\partial_1\eta_2  - \eta_2\partial_3\eta_3+
        \eta_1( \partial_1\eta_2 \partial_3\eta_3
-\partial_1\eta_3\partial_3\eta_2)\\
          \eta_1\partial_1\eta_3
+\eta_2\partial_2\eta_3+\eta_1(\partial_1\eta_3 \partial_2\eta_2
-\partial_1\eta_2\partial_2\eta_3)
        \end{array}\right) .
        \end{equation}

\subsection{Temporal derivative estimates}

In this subsection, we derive the estimates of temporal derivatives. For this purpose, we apply $\partial_t^j$ to \eqref{0101xxxx} to get
 \begin{equation}\label{01s06p}
\begin{cases}
\partial_t^{j+1}\eta=\partial_t^j u, \\[1mm]
 \partial_t^{j+1} u-\Delta_{\mathcal{A}} \partial_t^ju + \nabla_{\ml{A}} \partial_t^j q
= Q  \partial_{3}^2\partial_t^j\eta+
R \partial_t^{j} \vartheta e_3
 + N^{t,j}_u+N^{t,j}_q,\\[1mm]
P_{\vartheta} \partial_t^{j+1} \vartheta  =  \Delta_{\mathcal{A}}\partial_t^{j }\vartheta+R \partial_t^{j}u_3 +N^{t,j}_\vartheta , \\[1mm]
\div_{\mathcal{A}} \partial_t^j u  =\mm{div} D_u^{t,j},
\end{cases}
\end{equation}
where
\begin{align}
&N^{t,j}_u:= \sum_{0\leqslant m<j,\ 0\leqslant n\leqslant j} \partial_t^{j-m-n }\mathcal{A}_{il}\partial_{l}(\partial_t^{n}\mathcal{A}_{ik}\partial_t^{m}
\partial_ku),\nonumber\\
&N^{t,j}_q:=-\sum_{0\leqslant l< j}(\partial^{j-l}_t\mathcal{A}_{ik}\partial^l_t \partial_k q)_{3\times 1},\nonumber\\
&\label{partialusd} D^{t,j}_u:= \left(-\sum_{0\leqslant l<j}C_j^{j-l}\partial_t^{j-l} \mathcal{A}_{ki}\partial_t^l u_k\right)_{3\times 1},\\
&N^{t,j}_\vartheta:= \sum_{0\leqslant m<j,\ 0\leqslant n\leqslant j} \partial_t^{j-m-n }\mathcal{A}_{il}\partial_{l}(\partial_t^{n}\mathcal{A}_{ik}\partial_t^{m}
\partial_k\vartheta) +  R \partial_t^{j } \left( u  \cdot \nabla_{\tilde{\mathcal{A}}} \zeta_3+ u\cdot   \nabla \eta_3\right),\nonumber\\
& C_j^{j-l}\mbox{ denotes the  number of }(j-l)\mbox{-combinations from a given set }S\mbox{ of }j\mbox{ elements},\nonumber
\end{align}
and we have used the relation \eqref{diverelation} in \eqref{partialusd}. Then from \eqref{01s06p} we can establish the following estimates.
\begin{lem}\label{badiseqin}
Under the assumptions of \eqref{aprpioses} and \eqref{2019022201518},
\begin{align}
&\frac{1}{2}\frac{\mm{d}}{\mm{d}t}\left(E_1(\partial_t^j
\eta,\partial_t^j \vartheta)+\| \partial_t^j u\|_{0}^2\right)
+ \|\nabla\partial_t^j ( u,\vartheta)\|_0^2\nonumber \\
& \leqslant \frac{2 R}{P_\vartheta} \int\nabla \partial_t^{j} \eta_3\cdot \nabla \partial_t^{j }\vartheta  \mm{d}y  +
\left\{
  \begin{array}{ll}
c\sqrt{\mathcal{E}^H} \mathcal{D}^L  & \hbox{ for }  j=1;\\
  c\sqrt{\mathcal{E}^L} \mathcal{D}^H -\partial_t\int  \nabla \partial_t^{j-1}q\cdot D^{t,j}_u\mm{d}y& \hbox{ for }j=2,\ 3,
\end{array}
\right.
 \label{badiseqin30}\\
& \frac{1}{2} \frac{\mm{d}}{\mm{d}t} \| \nabla_{\mathcal{A}}  u_t\|_{0}^2 +\|
 u_{tt}\|^2_{0}\lesssim  \| (\partial_3^2  u, \theta_t)\|_0^2+\sqrt{\mathcal{E}^H}  \mathcal{D}^L,
\label{badiseqin31}
\\
& \frac{1}{2}\frac{\mm{d}}{\mm{d}t} \| \nabla_{\mathcal{A}}  \vartheta_t\|_{0}^2 +\|
\vartheta_{tt}\|^2_{0}\lesssim  \|   u_t\|_0^2+\sqrt{\mathcal{E}^H}  \mathcal{D}^L.
\label{badiseqin30xxxxxx}
\end{align}
\end{lem}
\begin{pf}
Multiplying \eqref{01s06p}$_2$ by $\partial_t^j u$ in $L^2$ and using the integral by parts, we get
\begin{align}
&\frac{1}{2}\frac{\mm{d}}{\mm{d}t}\left(Q \|\partial_{3} \partial_t^j\eta\|^2_{0}+\|\partial_t^j u\|_{0}^2\right) +\| \nabla \partial_t^ju\|^2_{0}\nonumber \\
&=R \int \partial_t^{j}u_3 \partial_t^j \vartheta  \mm{d}y +
\int \partial_t^j q \mm{div}_{\ml{A}}\partial_t^j u \mm{d}y +
\int N^{t,j}_u\cdot \partial_t^j u\mm{d}y\nonumber \\
&\quad+\int N^{t,j}_q \cdot \partial_t^j u\mm{d}y- \int\left(\nabla_{\tilde{\mathcal{A}}}\partial_t^j u:\nabla_{\mathcal{A}}\partial_t^ju+\nabla\partial_t^j u:\nabla_{\tilde{\mathcal{A}}}\partial_t^ju\right)\mm{d}y =:\sum_{k=1}^5I_k.
\label{estimdba}
\end{align}
The integrals $I_2$--$I_4$ can be estimated as follows (see (2.17)--(2.19) and Lemma 3.1 in \cite{JFJSJMFMOSERT}).
\begin{align}
& I_2 \leqslant  \left\{
               \begin{array}{ll}
         c \sqrt{\mathcal{E}^H} \mathcal{D}^L & \hbox{ for } j= 1;\\
        c\sqrt{\mathcal{E}^{L}}\mathcal{D}^H-\partial_t\int  \nabla \partial_t^{j-1}q\cdot D^{t,j}_u\mm{d}y & \hbox{ for } j=2,\  3,
               \end{array}
             \right. \label{201902171425} \\
& I_3,\ I_4 \lesssim \left\{
               \begin{array}{ll}
          \sqrt{\mathcal{E}^H} \mathcal{D}^L & \hbox{ for } j= 1;\\
            \sqrt{\mathcal{E}^L} \mathcal{D}^H & \hbox{ for } j=2,\  3.
               \end{array}
             \right. \label{201902171425x}
\end{align}

Similarly to the estimate of $I_3$, we also have
\begin{align}
& I_5\lesssim  \left\{
               \begin{array}{ll} \sqrt{\mathcal{E}^H} \mathcal{D}^L & \hbox{ for } j= 1;\\
            \sqrt{\mathcal{E}^L} \mathcal{D}^H & \hbox{ for } j=2,\  3,  \end{array}
             \right. \label{201902211620}
\\
&I_6:=\int\partial_t^{j} \nabla_{\tilde{\mathcal{A}}} \eta_3\cdot \nabla_{ {\mathcal{A}}} \partial_t^{j }\vartheta  \mm{d}y+\int\nabla \partial_t^{j} \eta_3\cdot \nabla_{\tilde{\mathcal{A}}} \partial_t^{j }\vartheta  \mm{d}y\nonumber \\
&\qquad +\int \partial_t^{j} \eta_3 N^{t,j}_\vartheta \mm{d}y \lesssim  \left\{
               \begin{array}{ll} \sqrt{\mathcal{E}^H} \mathcal{D}^L & \hbox{ for } j= 1;\\
            \sqrt{\mathcal{E}^L} \mathcal{D}^H & \hbox{ for } j=2,\  3.  \end{array}
             \right. \label{201902211619}
\end{align}
 Making using \eqref{01s06p}$_1$, \eqref{01s06p}$_3$, \eqref{201902211619} and the integral by parts, we have
\begin{align}
I_1
=  &R \left( \frac{\mm{d}}{\mm{d}t}\int\partial_t^{j} \eta_3 \partial_t^j \vartheta \mm{d}y - \frac{1}{P_\vartheta} \int\partial_t^{j} \eta_3  (  \Delta_{\mathcal{A}}\partial_t^{j }\vartheta+
R\partial_t^{j}u_3 +N^{t,j}_\vartheta ) \mm{d}y \right)\nonumber \\
 =&  R \left( \frac{\mm{d}}{\mm{d}t}\int\left(\partial_t^{j} \eta_3 \partial_t^j \vartheta -\frac{R|\partial_t^{j} \eta_3 |^2}{2 P_\vartheta}\right)\mm{d}y   + \frac{1}{P_\vartheta}\int \nabla \partial_t^{j}\eta_3\cdot \nabla \partial_t^{j }\vartheta  \mm{d}y  \right)  + \frac{R I_6}{P_\vartheta}   \nonumber \\
\leqslant & R \left(\frac{\mm{d}}{\mm{d}t}\int\left(\partial_t^{j} \eta_3 \partial_t^j \vartheta -\frac{R|\partial_t^{j} \eta_3 |^2}{ 2 P_\vartheta}\right)\mm{d}y   +\frac{1}{P_\vartheta} \int \nabla \partial_t^{j}\eta_3\cdot \nabla \partial_t^{j }\vartheta  \mm{d}y\right)\nonumber \\
&  + c\left\{
               \begin{array}{ll} \sqrt{\mathcal{E}^H} \mathcal{D}^L & \hbox{ for } j= 1;\\
            \sqrt{\mathcal{E}^L} \mathcal{D}^H & \hbox{ for } j=2,\  3. \end{array}
             \right.
\label{201902171425xx}
\end{align}
Consequently, putting \eqref{201902171425}--\eqref{201902211620} and \eqref{201902171425xx} into \eqref{estimdba} yields  the desired estimate
\begin{align}&
\frac{1}{2}\frac{\mm{d}}{\mm{d}t}\left(Q \|\partial_{3} \partial_t^j
\eta\|^2_{0}+  \frac{R^2}{P_\vartheta} \|\partial_t^{j} \eta_3\|^2_0+\| \partial_t^j u\|_{0}^2 -  2R   \int \partial_t^{j} \eta_3 \partial_t^j \vartheta \mm{d}y\right)
+ \| \nabla \partial_t^ju\|^2_{0}\nonumber \\
& \leqslant \frac{R}{P_\vartheta} \int\nabla \partial_t^{j} \eta_3\cdot \nabla \partial_t^{j }\vartheta  \mm{d}y  +
\left\{
  \begin{array}{ll}
c\sqrt{\mathcal{E}^H} \mathcal{D}^L  & \hbox{ for }  j=1;\\
  c\sqrt{\mathcal{E}^L} \mathcal{D}^H -\partial_t\int  \nabla \partial_t^{j-1}q\cdot D^{t,j}_u\mm{d}y& \hbox{ for }j=2,\ 3.
\end{array}
\right. \nonumber
  \end{align}
Similarly, we can easily derive from \eqref{01s06p}$_3$ that
  \begin{align}& \frac{1}{2}\frac{\mm{d}}{\mm{d}t}
 \left( \frac{ R^2 }{P_\vartheta} \|\partial_t^j \eta_3\|^2+ P_{\vartheta} \| \partial_t^j \vartheta \|_{0}^2- 2R  \int  \partial_t^j \eta_3 \partial_t^j \vartheta \mm{d}y
\right)  + \| \nabla   \partial_t^j\vartheta\|^2_{0}
\nonumber \\
&\leqslant \frac{R}{P_\vartheta} \int \nabla \partial_t^{j}\eta_3\cdot \nabla \partial_t^{j }\vartheta  \mm{d}y+ c \left\{
  \begin{array}{ll}
\sqrt{\mathcal{E}^H} \mathcal{D}^L  & \hbox{ for }  j=  1;\\
\sqrt{\mathcal{E}^L} \mathcal{D}^H & \hbox{ for }j=2,\ 3.
\end{array}
\right.\nonumber
\end{align}
Adding the above two estimates together yields \eqref{badiseqin30}.

Next we turn to the proof of \eqref{badiseqin31}. Multiplying \eqref{01s06p}$_2$ with $j=1$ by $u_{tt}$ in $L^2$,
and using \eqref{01s06p}$_1$ and the integral by parts, we conclude
\begin{align}
  &\frac{1}{2}\frac{\mm{d}}{\mm{d}t}\| \nabla_{\ml{A}} u_{t}\|^2_{0}  +\|u_{tt} \|_{0}^2 =
\int (  Q \partial_{3}^2u\cdot  u_{tt}+ R\theta_t \partial_{t}^2u_3)\mm{d}y\nonumber \\
 & +  \int  \left(   q_{t}  \mm{div}_{\ml{A}} u_{tt}+
  (N^{t,1}_u+ N^{t,1}_q) \cdot u_{tt} + \nabla_{\ml{A}}u_{t}: \nabla_{\ml{A}_t} u_{t}\right)\mm{d}y=:I_7+I_8.   \label{indeforut}
 \end{align}
The integral $I_8$ can be estimated as follows (see (2.26)--(2.29) in \cite{JFJSJMFMOSERT}):
\begin{align}
I_8 \lesssim
 \sqrt{\mathcal{E}^H} \mathcal{D}^L .\nonumber
\end{align}
In addition,
$$I_7\lesssim \|(\partial_3^2 u,\theta_t )\|_0\| u_{tt}  \|_0. $$
Thus, substituting the above two estimates  into \eqref{indeforut}, we immediately get \eqref{badiseqin31}. Similarly, we can easily derive \eqref{badiseqin30xxxxxx} from \eqref{01s06p}$_3$. This completes the proof. \hfill$\Box$
\end{pf}

\subsection{Horizontal spatial estimates}

In this subsection, we establish  the estimates of horizontal spatial derivatives.
For this purpose, we rewrite \eqref{0101xxxx} as the following non-homogeneous linear form:
\begin{equation}\label{s0106pnnnn}
\begin{cases}
 \eta_t=  u, \\[1mm]
 u_{t}+\nabla  q- \Delta  u- Q \partial_{3}^2\eta-  R\vartheta e_3= N^{\mm{h}}_u+N^{\mm{h}}_q,\\[1mm]
P_{\vartheta} \vartheta_t   - \Delta \vartheta - R u_3 = N^{\mm{h}}_\vartheta, \\[1mm]
\div u  =D^{\mm{h}}_u,\\
 (\eta,u,\theta)|_{\partial\Omega}=0,
\end{cases}
\end{equation}
where
$$\begin{aligned}
&N^{\mm{h}}_u:=\partial_l((\tilde{\ml{A}}_{jl}
\tilde{\mathcal{A}}_{jk}+
\tilde{\ml{A}}_{lk}+\tilde{\mathcal{A}}_{kl})\partial_ku_i)_{3\times 1},\
 N^{\mm{h}}_q:=-(\tilde{\mathcal{A}}_{ik}\partial_k q)_{3\times 1},\ D^{\mm{h}}_u:=\tilde{\mathcal{A}}_{lk}\partial_k u_l,\\
& N^{\mm{h}}_\vartheta:=
\partial_l((\tilde{\ml{A}}_{jl}
\tilde{\mathcal{A}}_{jk}+
\tilde{\ml{A}}_{lk}+\tilde{\mathcal{A}}_{kl})\partial_k\vartheta )  + R u\cdot (\nabla_{\tilde{\mathcal{A}}} \zeta_3+\nabla  \eta_3) .\end{aligned}$$
Then we have the following estimates on horizontal spatial derivatives of $(\eta,u,\vartheta)$.
\begin{lem}
\label{ssebadiseqin}
Let $0\leqslant j\leqslant 6$. Under the assumptions of \eqref{aprpioses} and \eqref{2019022201518},
\begin{align}
&\frac{\mm{d}}{\mm{d}t}\left(\int  \partial_{\mm{h}}^j \eta\cdot\partial_{\mm{h}}^j  u\mm{d}y +\frac{1}{2}\|\nabla\partial_{\mm{h}}^j\eta\|_{0}^2\right)
+ Q \|\partial_{3}\partial_{\mm{h}}^j\eta\|_{0}^2\nonumber \\
&  \leqslant \| \partial_{\mm{h}}^j  u\|^2_{0}+\left\{
              \begin{array}{ll}
 c\sqrt{\mathcal{E}^H} \mathcal{D}^L \mbox{ (or }\sqrt{\mathcal{E}^L} \mathcal{D}^H\mbox{)}   +  R\int\partial_{\mm{h}}^j\eta_3\partial_{\mm{h}}^j \vartheta\mm{d}y  & \mbox{ for } 0\leqslant j\leqslant 3 ; \\
  c \sqrt{\mathcal{E}^L}(\|(\eta,u)\|_7^2+\mathcal{D}^H)  + R\int\partial_{\mm{h}}^j\eta_3\partial_{\mm{h}}^j \vartheta\mm{d}y  & \mbox{ for } 4\leqslant  j\leqslant 6,
              \end{array}
            \right. \label{201902171647}
\\
&
\frac{1}{2} \frac{\mm{d}}{\mm{d}t}\left(Q \|\partial_{3}\partial_{\mm{h}}^j\eta\|_{0}^2+\|\partial_{\mm{h}}^j u\|^2_0
 \right)
 +  \|\nabla \partial_{\mm{h}}^j u\|_{0}^2
 \nonumber \\
& \leqslant  R  \int\partial_{\mm{h}}^ju_3\partial_{\mm{h}}^j\vartheta\mm{d}y
+c \left\{
              \begin{array}{ll}
 \sqrt{\mathcal{E}^H} \mathcal{D}^L &\mbox{ for }  j= 3 ; \\
 \sqrt{\mathcal{E}^L}(\|(\eta,u)\|_7^2+\mathcal{D}^H)     & \mbox{ for }   j= 6,
              \end{array}
            \right.
\label{201902201632} \\
&\frac{1}{2}  \frac{\mm{d}}{\mm{d}t} \left( E_1(\partial_{\mm{h}}^j \eta, \partial_{\mm{h}}^j \vartheta)+\|\partial_{\mm{h}}^j u\|^2_0
 \right)+ \|\nabla\partial_{\mm{h}}^j( u,\vartheta) \|_0^2\nonumber \\
&\leqslant  \frac{2 R}{P_\vartheta} \int\nabla \partial_{\mm{h}}^j\eta_3\cdot
\nabla\partial_{\mm{h}}^j \vartheta\mm{d}y+ c \left\{
  \begin{array}{ll}
\sqrt{\mathcal{E}^H} \mathcal{D}^L\mbox{ (or }\sqrt{\mathcal{E}^L} \mathcal{D}^H\mbox{)}    & \mbox{ for } 0\leqslant j\leqslant 3 ;\\
 \sqrt{\mathcal{E}^L}(\|(\eta,u)\|_7^2+\mathcal{D}^H) & \mbox{ for }  4\leqslant  j\leqslant 5.
\end{array}
\right. \label{201902171723}
\end{align}
\end{lem}
\begin{pf}
Applying $\partial_{\mm{h}}^j$ to \eqref{s0106pnnnn}$_2$, and then
multiplying the resulting equation by $\partial^j_{\mm{h}}\eta$ in $L^2$, we get
\begin{align}
&\frac{\mm{d}}{\mm{d}t}\left(\int \partial_{\mm{h}}^j \eta \cdot \partial_{\mm{h}}^j u \mm{d}y+\frac{1}{2}\|\nabla \partial_{\mm{h}}^j \eta\|_{0}^2\right) +Q  \|\partial_{3}\partial_{\mm{h}}^j\eta\|_{0}^2\nonumber \\
&= \| \partial_{\mm{h}}^j  u\|^2_{0}+ R\int\partial_{\mm{h}}^j\eta_3\partial_{\mm{h}}^j \vartheta\mm{d}y+\int \left(\partial_{\mm{h}}^j(N^{\mm{h}}_u+  N^{\mm{h}}_q)\cdot  \partial_{\mm{h}}^j\eta  + \partial_{\mm{h}}^j q\mm{div}\partial^j_{\mm{h}}\eta\right)\mm{d}y .
\label{estimforhoedsds1}
\end{align}
We denote  the last integral  in \eqref{estimforhoedsds1}  by $J_1$, then $J_1$ can be estimated as follows (see (2.33), (2.34), (3.5)--(3.7) in \cite{JFJSJMFMOSERT}).
 \begin{align}
J_1\lesssim  \left\{
  \begin{array}{ll}
 \sqrt{\mathcal{E}^H} \mathcal{D}^L \mbox{ (or }\sqrt{\mathcal{E}^L} \mathcal{D}^H\mbox{)}   & \mbox{ for } 0\leqslant j\leqslant 3;\\
 \sqrt{\mathcal{E}^L}(\|(\eta,u)\|_7^2+\mathcal{D}^H) & \mbox{ for }  4\leqslant j\leqslant 6 .
\end{array}
\right. \label{201902171705} \end{align}
Thus we immediately obtain the desired estimate \eqref{201902171647}.
We mention that the estimate \eqref{02181dsd} shall be used in the derivation of \eqref{201902171705}.

 Applying $\partial_{\mm{h}}^j$ to \eqref{s0106pnnnn}$_2$, and then
multiplying the resulting equality by $\partial^j_{\mm{h}} u$ in  $L^2$, we get
\begin{align}
&\frac{1}{2}\frac{\mm{d}}{\mm{d}t}\left(Q \|\partial_3\partial_{\mm{h}}^j\eta\|_{0}^2+\|\partial_{\mm{h}}^j u\|^2_0\right)+  \|\nabla \partial_{\mm{h}}^j u\|_{0}^2\nonumber \\
&=R \int\partial_{\mm{h}}^ju_3\partial_{\mm{h}}^j\vartheta\mm{d}y +\int \left(\partial_{\mm{h}}^j( N^{\mm{h}}_u+ N^{\mm{h}}_q)\cdot  \partial_{\mm{h}}^j u + \partial_{\mm{h}}^j q\mm{div}\partial^j_{\mm{h}} u\right)\mm{d}y=: J_2+J_3 .
\label{estimforhoe1nwe}
\end{align}
The last integral $J_3$ can be estimated as follows (see (2.37), (2.38) and Lemma 3.3 in \cite{JFJSJMFMOSERT}).
 \begin{align}
J_3 \lesssim  \left\{
  \begin{array}{ll}
 \sqrt{\mathcal{E}^H} \mathcal{D}^L\mbox{ (or }\sqrt{\mathcal{E}^L} \mathcal{D}^H\mbox{)}   & \mbox{ for } 0\leqslant j\leqslant 3 ;\\
  \sqrt{\mathcal{E}^L}(\|(\eta,u)\|_7^2+\mathcal{D}^H) & \mbox{ for }  4\leqslant j\leqslant 6 .\label{201902171836}
\end{array}
\right. \end{align}

 Using \eqref{01s06p}$_1$, \eqref{01s06p}$_3$ and the integral by parts, we have, for $ j\neq 6$,
\begin{align}
J_2 =&  R \left(\frac{\mm{d}}{\mm{d}t}\int\partial_{\mm{h}}^j\eta_3\partial_{\mm{h}}^j\vartheta\mm{d}y
- \frac{1}{P_\vartheta}\int\partial_{\mm{h}}^j\eta_3\partial_{\mm{h}}^j
(   \Delta \vartheta + R u_3 + N^{\mm{h}}_\vartheta )\mm{d}y\right)\nonumber \\
 = & R \left(\frac{\mm{d}}{\mm{d}t}\left(
 \int\left(\partial_{\mm{h}}^j\eta_3\partial_{\mm{h}}^j\vartheta
-\frac{R |\partial_{\mm{h}}^j\eta_3|^2}{2P_\vartheta}\right)\mm{d}y \right)
+ \frac{1}{P_\vartheta}\int(  \nabla \partial_{\mm{h}}^j\eta_3\cdot
\nabla\partial_{\mm{h}}^j \vartheta -\partial_{\mm{h}}^j\eta_3\partial_{\mm{h}}^j N^{\mm{h}}_\vartheta  )\mm{d}y\right). \label{201902171836xx}
\end{align}
In addition, similarly to the estimate of $J_3$, we also estimate that
\begin{align}
&\int \partial_{\mm{h}}^{j} \eta_3 \partial_{\mm{h}}^j N^{\mm{h}}_\vartheta  \mm{d}y \lesssim  \left\{
               \begin{array}{ll} \sqrt{\mathcal{E}^H} \mathcal{D}^L\mbox{ (or }\sqrt{\mathcal{E}^L} \mathcal{D}^H\mbox{)}   & \hbox{ for } 0\leqslant  j\leqslant  3;\\
            \sqrt{\mathcal{E}^L}\|(\eta,u)\|_7^2& \hbox{ for } 4\leqslant  j\leqslant 5. \end{array}
             \right.
\label{201902261641}
\end{align}

Consequently, making use of \eqref{201902171836}--\eqref{201902261641}, we derive from \eqref{estimforhoe1nwe} that
\begin{align}
&\frac{1}{2} \frac{\mm{d}}{\mm{d}t}\left(Q \|\partial_{3}\partial_{\mm{h}}^j\eta\|_{0}^2+\|\partial_{\mm{h}}^j u\|^2_0
 \right)
 +  \|\nabla \partial_{\mm{h}}^j u\|_{0}^2 \nonumber\\
& \leqslant  c \left\{
              \begin{array}{ll}
 \sqrt{\mathcal{E}^H} \mathcal{D}^L\mbox{ (or }\sqrt{\mathcal{E}^L} \mathcal{D}^H\mbox{)}   &\mbox{ for } 0\leqslant j\leqslant 3 ; \\
 \sqrt{\mathcal{E}^L}(\|(\eta,u)\|_7^2+\mathcal{D}^H)     & \mbox{ for }  4\leqslant j\leqslant 6
              \end{array}
            \right. \nonumber
\\
&\quad  +R \left\{
              \begin{array}{ll}
 \partial_t\left(
   \int\partial_{\mm{h}}^j\eta_3\partial_{\mm{h}}^j\vartheta\mm{d}y
- {R} \| \partial_{\mm{h}}^j\eta_3\|_0^2/2P_\vartheta \right)+\int\nabla \partial_{\mm{h}}^j\eta_3\cdot
\nabla\partial_{\mm{h}}^j \vartheta\mm{d}y/  P_\vartheta& \mbox{ for } j\neq 6 ; \\
\int\partial_{\mm{h}}^ju_3\partial_{\mm{h}}^j\vartheta\mm{d}y & \mbox{ for } j= 3,\ 6 ,
              \end{array}
            \right. \nonumber
\end{align}
which yields \eqref{201902201632} for $j=3$ and $6$.
Similarly, we easily derive from \eqref{s0106pnnnn}$_3$ that
\begin{align}
&\frac{1}{2}  \frac{\mm{d}}{\mm{d}t} \left(  \frac{R^2}{P_\vartheta} \| \partial_{\mm{h}}^j\eta_3\|_0^2+P_{\vartheta}  \| \partial_{\mm{h}}^j \vartheta \|_{ 0}^2
- 2R \int\partial_{\mm{h}}^j\eta_3\partial_{\mm{h}}^j\vartheta\mm{d}y
 \right)+ \| \nabla \partial_{\mm{h}}^j\vartheta\|^2_{0}\nonumber \\
&\leqslant
  \frac{R}{P_\vartheta} \int\nabla \partial_{\mm{h}}^j\eta_3\cdot
\nabla\partial_{\mm{h}}^j \vartheta\mm{d}y + c \left\{
  \begin{array}{ll}
\sqrt{\mathcal{E}^H} \mathcal{D}^L\mbox{ (or }\sqrt{\mathcal{E}^L} \mathcal{D}^H\mbox{)}    & \mbox{ for } 0\leqslant j\leqslant 3 ;\\
\sqrt{\mathcal{E}^L} \|(\eta,\vartheta)\|_7^2 & \mbox{ for }  4\leqslant  j\leqslant 5 .
\end{array}
\right.\nonumber
\end{align}
Adding the above two estimates together yields \eqref{201902171723}.
\hfill$\Box$
\end{pf}

\subsection{Stokes estimates and elliptic estimates}

In this subsection, we use the regularity theories of the both of Stokes
and elliptic problems to derive estimates of hybrid derivative of $(\eta,u)$. To this purpose, we rewrite \eqref{s0106pnnnn}$_2$ and \eqref{s0106pnnnn}$_4$ as the following Stokes problem
\begin{equation}
\label{Stokesequson}
\begin{cases}
-\Delta  \omega+\nabla  q =  R \vartheta e_3- Q\Delta_{\mm{h}}  \eta  - u_{t}+   N^{\mm{h}}_u+ N^{\mm{h}}_q,\\
\div \omega= D^{\mm{h}}_u+ Q \mm{div}\eta
\end{cases}
\end{equation}
coupled with boundary-value condition
\begin{equation}
\label{boudnStokesequson}
\omega|_{\partial\Omega}=0,
\end{equation}
where $\omega= Q \eta+ u$.

Now, applying $\partial^{k}_{\mm{h}}$ to \eqref{Stokesequson} and \eqref{boudnStokesequson}, we get
\begin{equation*}
\begin{cases}
-\Delta \partial^{k}_{\mm{h}}\omega+\nabla  \partial^{k}_{\mm{h}}q = \partial^{k}_{\mm{h}}\left(R \vartheta e_3- Q
 \Delta_{\mm{h}}  \eta  - u_{t}+ N^{\mm{h}}_u+N^{\mm{h}}_q\right),\\
\div \partial^{k}_{\mm{h}}\omega =\partial^{k}_{\mm{h}}\left( D^{\mm{h}}_u+Q\mm{div} \eta\right),\\
\partial^{k}_{\mm{h}}\omega|_{\partial\Omega}=0.
\end{cases}
\end{equation*}
Then we can exploit the classical regularity theory of Stokes problem (see \cite[Proposition 2.3]{TRNSETNA}) to deduce that
\begin{equation}   \label{omessetsim}
\|\omega\|_{k,i-k+2}^2+\|\nabla  q\|_{k,i-k}^2 \lesssim \|( \Delta_{\mm{h}}  \eta, u_t,\vartheta)\|_{k,i-k}^2+S_{i}  , \end{equation}
where $ S_{i} :=\|  (N^{\mm{h}}_u, N^{\mm{h}}_q)\|_{i}^2 +\|(D^{\mm{h}}_u ,\mm{div}\eta) \|_{i+1}^2$.

Applying $\partial_t^k$ to  \eqref{s0106pnnnn}$_2$, \eqref{s0106pnnnn}$_4$ and \eqref{s0106pnnnn}$_5$, we see that
$$
\begin{cases}
\nabla  \partial_{t}^k q- \Delta  \partial_{t}^k u=\partial_{t}^k
\left(Q \partial_{3}^2 \eta +R \vartheta e_3-u_t
+ N^{\mm{h}}_u+ N^{\mm{h}}_q\right), \\[1mm]
\div  \partial_{t}^ku  = \partial_{t}^kD^{\mm{h}}_u, \\
 \partial_{t}^ku|_{\partial\Omega}=0 .
\end{cases}$$
Hence, we use again the classical regularity theory of Stokes problem to get
\begin{align}\label{02181}
&  \|\partial_t^k u\|_{i-2k+2}^2+\|\nabla\partial_t^k  q\|_{i-2k}^2\nonumber \\
& \lesssim\|\partial^{k}_t(\partial_{3}^2\eta  ,u_t,\vartheta)\|_{i-2k}^2 +  \|\partial_t^k (  N^{\mm{h}}_u, N^{\mm{h}}_q)\|_{i-2k}^2+\|\partial_t^k D^{\mm{h}}_u\|_{i-2k+1}^2. \end{align}

Similarly, applying $\partial_t^k$ to \eqref{s0106pnnnn}$_3$ and \eqref{s0106pnnnn}$_5$, we see that
$$
\begin{cases}
 -  \Delta  \partial_{t}^k \vartheta = \partial_{t}^k(R  u_3- P_{\vartheta} \vartheta_t   + N^{\mm{h}}_\vartheta),   \\
 \partial_{t}^k\vartheta|_{\partial\Omega}=0 .
\end{cases}
$$
By the classical regularity theory of elliptic problem (see Theorem 5 in \cite[Section 6.3]{ELGP}), we have
\begin{align}
\|\partial_t^k \vartheta\|_{i-2k+2}^2
 \lesssim \|\partial^{k}_t( u_3,\vartheta_t, N^{\mm{h}}_\vartheta)\|_{i-2k}^2 .
\label{201902261711}
\end{align}
Thus we can easily derive the following three lemmas from  \eqref{omessetsim}--\eqref{201902261711}, resp..
%%%%%%%%%%%%%%%%%%%%%%%%%%%%%
\begin{lem}\label{lem:dfifessim}
Under the assumptions of \eqref{aprpioses} and \eqref{2019022201518}, it holds that
\begin{align}
&   \frac{\mm{d}}{\mm{d}t}  \|\eta\|_{i+2,*}^2+ \|  (\eta ,u)\|_{i+2}^2 +\|\nabla q\|_i^2\nonumber  \\
&\lesssim \left\{
   \begin{array}{ll} \|  \eta\|_{ {3},0}^2 +\|  (u_{t}, \vartheta) \|_1^2
+
 {\mathcal{E}^H }{ {\mathcal{D}}^L}, & \hbox{for }i=1; \\
 \|  \eta\|_{6,0}^2 +\|  (u_{t}, \vartheta) \|_4^2
+
 {\mathcal{E}^L}\mathcal{D}^H & \hbox{for }i=4;
 \\
\mathcal{E}^H+\mathcal{D}^H  &  \hbox{for }i=5
   \end{array}
 \right. \label{dfifessim}
\end{align}
on $(0,T]$, where the norm $\|\eta\|_{i+2,*}$ is equivalent to $ \|\eta\|_{i+2}$ for $i=1$, $4$ and $5$.
\end{lem}
\begin{pf}
Noting that by virtue of \eqref{s0106pnnnn}$_1$,
$$
\begin{aligned}
&\|\omega\|_{k,i-k+2}^2 =\|(Q \eta, u)\|_{k,i-k+2}^2
+ {Q} \frac{\mm{d}}{\mm{d}t} \|\eta\|_{k,i-k+2}^2,
\end{aligned}
$$
we deduce from \eqref{omessetsim} that
 \begin{align}
 & {Q}  \frac{\mm{d}}{\mm{d}t} \|\eta\|_{{ 2k,i- 2k+2}}^2
 +  \|   (Q\eta,  u)\|_{ 2k,i- 2k+2  }^2+\|\nabla   q\|_{ 2k,i- 2k}^2 \nonumber \\
 &\lesssim  \| \eta\|_{ 2k+2, i- 2k}^2+\| ( u_t,\theta)\|_{i}^2+S_{i},\label{201903201620} \\
 & {Q}  \frac{\mm{d}}{\mm{d}t} \|\eta\|_{1,2}^2
 +  \|   (Q\eta,  u)\|_{1,2}^2+\|\nabla   q\|_{1,0}^2\lesssim  \| \eta\|_{3,0 }^2+\| ( u_t,\theta)\|_{3}^2+S_3.
 \label{omdm12dfs2}
 \end{align}
Exploiting the recursion formula \eqref{201903201620} from $k=0$ to $[i/2]$,
we see that there are positive constants $c_k$, $0\leqslant k\leqslant [i/2] $, such that
\begin{align}
& \frac{\mm{d}}{\mm{d}t} \sum_{k=0}^{[i/2]} c_k\|\eta\|_{{2k,i-2k+2}}^2
 +  \| (\eta, u)\|_{i+2}^2+\|\nabla   q\|_i^2\nonumber \\
 & \lesssim
 \|  \eta\|_{2+2[i/2],i-2[i/2]}^2+ \|  (u_t,\vartheta)\|_{i}^2+ S_i ,
 \end{align}
 where $[i/2]:=$ the integer part of $i/2$.
In addition,
\begin{equation}
\label{201902202048}
S_i  \lesssim \left\{
                 \begin{array}{ll}
 {\mathcal{E}^H E^L} \lesssim {\mathcal{E}^H  {\mathcal{D}}^L},
                   & \hbox{for }i=1;\\
             \mathcal{E}^L { {E}}^H \lesssim  \mathcal{E}^L \mathcal{D}^H & \hbox{for }i=4;  \\
 \mathcal{E}^H \|\eta\|_7^2 + \mathcal{D}^H & \hbox{for }i=5,
                 \end{array}
               \right.
\end{equation}
where  ${{E}^H}:=\mathcal{E}^H-\|\nabla \eta\|_{6,0}^2$ and $E^L:=\mathcal{E}^L-\|\nabla \eta\|_{3,0}^2$  (see  (2.48), the first estimate before (3.16), and the proof of (3.19) in \cite{JFJSJMFMOSERT} for the derivation of \eqref{201902202048}). Thus \eqref{dfifessim} immediately follows by \eqref{201903201620}--\eqref{201902202048}. \hfill$\Box$
\end{pf}
\begin{lem}\label{lem:201902181540}
Under the assumptions of \eqref{aprpioses} and \eqref{2019022201518}, we have the following estimates:
\begin{align} &  \label{dfifessimlas}
 \| u\|_{3}^2+\|\nabla   q\|_{1}^2 \lesssim  \|(\partial_{3}^2\eta,u_t,\vartheta)\|_1^2+ {\ml{E}^H}E^L, \\
  & \label{dfifeddssimlas}\|u_t\|_{2}^2+\|\nabla   q_t\|_{0}^2   \lesssim \|(\partial_{3}^2u,u_{tt},\vartheta_t)\|_{0}^2+ {\ml{E}^H}\ml{D}^L,\\
\label{omdm122nsdfsf}
&   \sum_{k=0}^2(\| \partial_t^k u\|_{6-2k}^2+ \|\nabla\partial_t^k q\|_{4-2k}^2)
  \lesssim  \|\eta\|_6^2+ \| ( u,u_t, \partial_t^3 u, \vartheta_{tt})\|_0^2+
\|\vartheta_t\|_2^2+\| \vartheta \|_4^2+ {\mathcal{E}^L} {E}^H,\\
&\label{highestdids}
  \sum_{k=1}^2( \| \partial_t^k u\|_{7-2k}^2+\|\nabla\partial_t^k  q\|_{5-2k}^2)
  \lesssim \|u\|_5^2+ \|( u_t, \partial_t^3 u, \vartheta_{tt})\|_1^2+
 \| \vartheta_t\|_3^2+\mathcal{E}^H\mathcal{D}^H,\\
& \label{02181dsd}
 \|u\|_{7}+\|\nabla q\|_{5} \lesssim \|\eta\|_{7}  + \sqrt{\mathcal{D}^H}.
\end{align}
\end{lem}
\begin{pf}
Since the proof is standard, we omit it (or please refer to (2.44), (2.45), (3.9), (3.11) and (3.2) in \cite{JFJSJMFMOSERT} for \eqref{dfifessimlas}--\eqref{02181dsd}, resp.). \hfill$\Box$
\end{pf}

\begin{lem}
\label{201903231042}
Under the assumptions of \eqref{aprpioses} and \eqref{2019022201518}, we have the following estimates:
\begin{align}
&  \| \vartheta\|_{3}^2 \lesssim \| (u_3, \vartheta_t)\|_{1}^2+ {\ml{E}^H}E^L,
\label{201902262135}
\\
 & \|\vartheta\|_4^2+ \| \vartheta_t\|_2^2 \lesssim \|u_3\|_2^2+ \|(u_t, \theta_{tt})\|_0^2 + {\ml{E}^H}\mathcal{D}^L,
\label{omdm122nsdfsfxxx}
\\
&   \sum_{k=0}^2 \| \partial_t^k \vartheta\|_{6-2k}^2
  \lesssim  \| (   u_{tt}, \partial_t^3\vartheta )\|_0^2+
\|u_t\|_2^2+\|u\|_4^2+ {\mathcal{E}^L} {E}^H,\label{highestdidsxxxxxxx}\\
&\label{highestdidsxxxx}
  \sum_{k=0}^2  \| \partial_t^k \vartheta\|_{7-2k}^2
  \lesssim \| (   u_{tt}, \partial_t^3\vartheta )\|_1^2+
\|u_t\|_3^2+\|u\|_5^2+\mathcal{E}^H\mathcal{D}^H .
\end{align}
 \end{lem}
\begin{pf}
The proof is very similar to Lemma \ref{lem:201902181540}, so we omit it. \hfill$\Box$
\end{pf}

\subsection{Modified estimate for $\eta$ and an equivalence estimate}

Since  $R$, $Q$ and $P_\vartheta$ satisfies \eqref{201903032131},  then,  for any $(\varphi,\psi,\phi)\in H^4_\sigma\times H^3_\sigma\times H_0^3$,
\begin{align}
&E_1(\varphi, \phi)\geqslant (1-\Upsilon_1)\left(Q \|\partial_{3}\varphi\|^2_{0}+ {2R^2}P_\vartheta^{-1}\|\varphi_3\|^2_0+
P_{\vartheta} \|\phi \|_{0}^2\right),
\label{201903011036} \\
& \tilde{{\mathcal{D}}}_1(\varphi,\psi,\phi) \geqslant (1-\Upsilon_2) \tilde{{\mathcal{D}}}_{1,1}(\varphi,\psi,\phi).\label{201903011036xx}
\end{align}

Using the integral by parts, Young's inequality  and  Wirtinger's inequality
$$ { \pi^2} \|w\|_0^2\leqslant \|\partial_3 w\|_0^2\mbox{ for any }w\in H_0^1 \mbox{ see
 \cite[Remark 4]{JFJSOMARMA}}, $$
we can check that
$$
\begin{aligned}
 \Gamma( \varphi,\psi,\phi):=& \frac{-\sum_{0\leqslant |\alpha|\leqslant 3}
 \int  \partial^{\alpha}_{\mm{h}}\varphi\cdot \partial^{\alpha}_{\mm{h}} \psi\mm{d}y
 }{\tilde{\mathcal{E}}_{1,1}(\varphi,\psi,\phi)}\\
 \leqslant & \frac{
 2\|\varphi\|_0^2+3\|\varphi\|_{2,0}^2 +\|\varphi\|_{3,0}^2  +2( \|\psi\|_0^2+\|\psi\|_{2,0}^2+\|\psi\|_{3,0}^2)
 }{
\sqrt{2}(  \|\nabla \varphi\|_{\underline{3},0}^2 +
2\|  \psi \|^2_{\underline{2},0}) }\leqslant \frac{1}{\sqrt{2}},
\end{aligned}$$
which implies that
\begin{align}
& \tilde{\mathcal{E}}_1 (\varphi,\psi,\phi) \geqslant(1-1/\sqrt{2})\tilde{\mathcal{E}}_{1,1}(\varphi,\psi,\phi ). \label{201903011036x}
\end{align}

We expect that the above estimates can be applied to  $(\eta,u,\vartheta)$.
However \eqref{201903011036}--\eqref{201903011036x} fail for  $(\eta,u,\vartheta)$, since $(\eta,u)$ does not satisfy $\mm{div}\eta=0$ and $\mm{div}u=0$. Next we shall modify the above three estimates \eqref{201903011036}--\eqref{201903011036xx}.
\begin{lem}
\label{lem:201903101813}
If $R$, $Q$ and $P_\vartheta$ satisfies \eqref{201903032131},  we have the following estimates: for any $(\varphi,\psi,\phi)\in H^{4}_0\times H^3_0\times H_0^3$,
 \begin{align}
& E_1(\varphi,\phi)\geqslant (1-\Upsilon_1)\left(Q \|\partial_{3}\varphi \|^2_{0}+  2  R^2 P_\vartheta^{-1}\|\varphi_3\|^2_0+
P_{\vartheta} \| \phi\|_{0}^2\right)-c\Xi_1(\varphi), \label{201903011352} \\
&\tilde{\mathcal{E}}_1(\varphi,\psi,\phi) \geqslant(1-1/\sqrt{2})\tilde{\mathcal{E}}_{1,1}(\varphi,\psi,\phi) -c\Xi_2(\varphi,\psi)
, \\
& \tilde{{\mathcal{D}}}_1(\varphi,\psi,\phi)   \geqslant (1-\Upsilon_2) \tilde{{\mathcal{D}}}_{1,1}(\varphi,\psi,\phi) -c\Xi_2(\varphi,\psi),\label{201903011352xx}
\end{align}
where
\begin{align}
& \Xi_1(\varphi):= \|\mm{div}\varphi\|_1^2+\|\mm{div}\varphi\|_1\|\nabla \varphi\|_{\underline{1},0} ,\nonumber \\
&\Xi_2(\varphi,\psi):= \|\mm{div}\varphi\|_{3}^2+\|\mm{div}\psi\|_{2}^2+
\|\mm{div}\varphi\|_{3}\|
\nabla \varphi\|_{\underline{3},0}+\|\mm{div}\psi\|_{2}\|
 \psi\|_{\underline{3},0} ,\nonumber
  \end{align}and the positive constant $c$ depends on  the domain, and the parameters  $\Upsilon_1$, $\Upsilon_2$ and $\Upsilon_3$.
\end{lem}
\begin{pf}
By the existence theory of the Stokes problem, there exists $\tilde{\eta}^i$ such that
\begin{equation}
\label{Stokesequsonxx}
\begin{cases}
\nabla  q -\Delta \tilde{\eta}^i= 0,\
\div  \tilde{\eta}^i  = \mm{div}\Psi^i,\\
\tilde{\eta}^i|_{\partial\Omega}=0
\end{cases}
\end{equation}
and, for $0\leqslant j\leqslant 2 $,
\begin{equation} \nonumber
\|\tilde{\eta}^i\|_{j+2}\lesssim \|\mm{div}\Psi^i\|_{j+1},
\end{equation}
where $\Psi^1=\varphi$ and $\Psi^2=\psi$.
Then we can derive from \eqref{201903011036}--\eqref{201903011036x} that
\begin{align}
&E_1(  \varphi-\tilde{\eta}^1 , \phi )\geqslant (1-\Upsilon_1)\left(Q \|\partial_{3} (\varphi-\tilde{\eta}^1)\|^2_{0}+ 2R^2P_\vartheta^{-1}\| \varphi_3-\tilde{\eta}_3^1\|^2_0+
P_{\vartheta} \|\phi \|_{0}^2\right),\nonumber \\
&\tilde{\mathcal{E}}_1 ( \varphi-\tilde{\eta}^1,\psi-\tilde{\eta}^2,\phi) \geqslant(1-1/\sqrt{2})\tilde{\mathcal{E}}_{1,1}( \varphi-\tilde{\eta}^1,\psi-\tilde{\eta}^2,\phi  )
,\nonumber  \\
& \tilde{{\mathcal{D}}}_1( \varphi-\tilde{\eta}^1,\psi-\tilde{\eta}^2,\phi) \geqslant (1-\Upsilon_2) \tilde{{\mathcal{D}}}_{1,1}( \varphi-\tilde{\eta}^1,\psi-\tilde{\eta}^2,\phi  ).\nonumber
\end{align}Thus we can immediately deduce \eqref{201903011352}--\eqref{201903011352xx} from the four estimates above. \hfill $\Box$
\end{pf}

Finally we introduce an equivalence estimate.
\begin{lem}
\label{201903081525}
Under the assumptions of \eqref{aprpioses} and \eqref{2019022201518},
\be \label{esdoval}
\mathcal{E}^H\lesssim  \mathcal{E}:=\|\eta\|_{7}^2 +\|(u,\vartheta)\|_6^2.
\ee
\end{lem}
\begin{pf}
The proof is trivial, please refer to Lemma 3.6 in \cite{JFJSJMFMOSERT}.
 \hfill $\Box$
\end{pf}

\subsection{Energy inequalities}

Next we establish
lower-order energy and higher-order inequalities.
In what follows the letters $c_1^L$, $c_2^L$, $c_1^H$ and $ c_2^H$, will denote generic constants
which may depend on the domain $\Omega$ and the known physical parameters.
\begin{pro}
\label{pro:0301n}
Let $R$, $Q$ and $P_\vartheta$ satisfy \eqref{201903032131}. Then, under the assumptions of \eqref{aprpioses} and \eqref{2019022201518}, there exist functionals  $\tilde{\ml{E}}^L$ and  $\tilde{\ml{E}}^H$
which are equivalent to $\mathcal{E}^L$ and $\mathcal{E}^H$, resp.,
such that
\begin{align}
 \label{emdslds}
& \frac{\mm{d}}{\mm{d}t} \tilde{\mathcal{E}}^{L}+\mathcal{D}^{L}\leqslant  0 ,\\ \label{hgheefdsstim}
& \frac{\mm{d}}{\mm{d}t}\tilde{\mathcal{E}}^{H}+\mathcal{D}^{H}\lesssim \sqrt{\mathcal{E}^L}\|(\eta,u)\|_7^2.
\end{align}
\end{pro}
\begin{pf} (1)
We can derive from \eqref{201902171647} with $0\leqslant j \leqslant 3$, and \eqref{201902171723} with $0\leqslant j \leqslant 2$ that there exists a constant $c_1^L$ such that
\begin{equation} \label{basdicl1s}
\frac{\mm{d}}{\mm{d}t}\tilde{\mathcal{E}}_1(\eta,u,\vartheta) +\tilde{\mathcal{D}}_1(\eta,u,\vartheta) \leqslant c_1^L  \sqrt{\mathcal{E}^H} \mathcal{D}^L .\end{equation}
Then we deduce from \eqref{badiseqin30} for $j=1$, \eqref{badiseqin31}, \eqref{badiseqin30xxxxxx}, \eqref{201902201632} with $ j =3$, \eqref{dfifessim} with $i=1$ and \eqref{basdicl1s} that
 there is constant $c_2^L$,  such that, for any $c_1^l>1$ and any $c_1^s\in (0,1)$,
\begin{align}
\frac{\mm{d}}{\mm{d}t}\tilde{\mathcal{E}}^L +\tilde{\mathcal{D}}^L \leqslant  &
\frac{2R}{P_\vartheta}\int\nabla  u_3\cdot \nabla \vartheta_t  \mm{d}y+
  R\sum_{|\alpha|=3}\int \partial_{\mm{h}}^\alpha u_3
 \partial_{\mm{h}}^\alpha \vartheta\mm{d}y+c_2^L\big(c_1^s(\|  \eta\|_{3,0}^2
 \nonumber \\
& +\|  (u_{t}, \vartheta) \|_{1}^2
+ \| (\partial_3^2  u, \theta_t\|_0^2
)
+ (1+c_1^s+c_1^l ) \sqrt{\mathcal{E}^H} \mathcal{D}^L \big),\label{basdicl1sxx}\end{align}
where
$$ \begin{aligned}
\tilde{\mathcal{E}}^L:=& c_1^l\tilde{\mathcal{E}}_1(\eta,u,\vartheta)
+ (
 Q\|\partial_{3}\eta\|_{3,0}^2+\|u\|_{3,0}^2  +\|  u_t\|_0^2
\\
&+E_1(
u,\vartheta_t)  +   c_1^s(\| \nabla_{\mathcal{A}} ( u_t, \vartheta_t)\|_{0}^2
+2\|\eta\|_{3,*}^2))  /2 ,  \\
\tilde{\mathcal{D}}^L:=& c_1^l\tilde{{\mathcal{D}}}_1(\eta,u,\vartheta )+  \|\nabla (u_t,\vartheta_t)\|_0^2
+ \|\nabla  u\|_{ 3,0}^2 \\
& +
 c_1^s(\|(u,\vartheta)_{tt}\|_0^2 +  \|  (\eta ,u)\|_{3}^2 +\|\nabla q\|_1^2).
\end{aligned}$$

Making use of \eqref{lemdssdf}, Lemma \ref{lem:201903101813},  interpolation inequality (see \cite[5.2 Theorem]{ARAJJFF}), Friedrichs' inequality  (see \cite[Lemma 1.42]{NASII04}), and Wirtinger's inequality,
we can estimate that
\begin{align}
&   c_1^l(\|\nabla \eta\|_{\underline{3},0}^2+\|( u,\vartheta)\|_{\underline{2},0}^2)
  +   c_1^s(\|(u, \vartheta)_t\|_{1}^2
+\|\eta\|_{3}^2)   -c_3^L\|(\eta,u)\|_5\|(\eta,u)\|_3^2
\lesssim  \tilde{\mathcal{E}}^L\lesssim \mathcal{E}^L,\label{201903182050} \\
&\tilde{\mathcal{D}}^L_1:= c_1^l(\|\partial_3\eta\|_{\underline{3},0}^2+\|(u,\vartheta)\|_{\underline{2},1}^2 )+  \|(u,\vartheta)_t\|_1^2
+ \|\nabla  u\|_{ 3,0}^2 \nonumber  \\
&\qquad \quad  +
 c_1^s(\|(u,\vartheta)_{tt}\|_0^2 +  \|  (\eta ,u)\|_{3}^2 +\|\nabla q\|_1^2)  -c_3^L\|(\eta,u)\|_5\|(\eta,u)\|_3^2
\lesssim  \tilde{\mathcal{D}}^L\lesssim \mathcal{D}^L,\label{xx201903182050} \\
&\|u\|_2^2\leqslant \|u\|_0^{2/3}\|u\|_3^{4/3}\lesssim \|u\|_0^{2/3}  \|(\partial_{3}^2\eta,u_t,\vartheta)\|_1^{4/3}+ {\ml{E}^H}E^L, \label{xx201903182050xx}
\end{align}
where the constant $c_3^L$ depends on $c_1^l$.
By \eqref{xx201903182050} and \eqref{xx201903182050xx}, we immediately see from \eqref{basdicl1sxx} that,  for properly large $c_l^l$ and properly small $c_1^s$,
\begin{align}
\frac{\mm{d}}{\mm{d}t}\tilde{\mathcal{E}}^L +c_4^L (\tilde{\mathcal{D}}^L_1+ c_3^L\|\eta\|_5\|\eta\|_3^2) \lesssim   &  \sqrt{\mathcal{E}^H} \mathcal{D}^L .\label{basdicl1sxxsdafa}\end{align}

Making use of \eqref{dfifessimlas}, \eqref{dfifeddssimlas},
\eqref{201902262135} and \eqref{omdm122nsdfsfxxx},
we have
\begin{align}
& \|(u,\vartheta)\|_3^2+\|\nabla q\|_1^2 \lesssim \|(u,\vartheta)\|_0^2+\|(u,\vartheta)_t\|_1^2+\|\eta\|_3^2+   {\ml{E}^H}E^L,\label{201903201819}\\
&\|(u, \vartheta)_t\|_2^2+\|\vartheta\|_4^2+\|\nabla q_t\|_0^2 \lesssim \|u\|_2^2+ \|u_t, (u,\theta)_{tt}\|_0^2 + {\ml{E}^H}\mathcal{D}^L . \label{201903182050xxxx}
\end{align}
Then we can derive from \eqref{201903182050}--\eqref{xx201903182050}, \eqref{201903201819} and \eqref{201903182050xxxx} that, there exists $\delta^L<1$, such that, for any sufficiently small $\delta\in (0,\delta^L)$,
 \begin{align}
 \mathcal{E}^L\mbox{ and }\mathcal{D}^L\mbox{ are equivalent to }\tilde{\mathcal{E}}^L\mbox{ and }
(\tilde{\mathcal{D}}^L_1+ c_3^L\|\eta\|_5\|\eta\|_3^2)\mbox{ resp.},\label{201902282100}
 \end{align}
 where the equivalence coefficients can be independent of $\delta$.
Consequently we derive from \eqref{basdicl1sxx} and \eqref{201902282100} that
\begin{equation}\label{defindgfdfeofengylasdf}
 \frac{\mm{d}}{\mm{d}t}\tilde{\mathcal{E}}^L +c_5^L {\mathcal{D}}^L \leqslant0,
\end{equation}
which yields \eqref{emdslds} by re-defining $\tilde{\mathcal{E}}^L := \tilde{\mathcal{E}}^L/c_5^L$.

(2) We turn to derive the higher-order inequality \eqref{hgheefdsstim}. We can derive from \eqref{201902171647} with $3\leqslant j \leqslant 6$,  and \eqref{201902171723} with $3\leqslant j \leqslant 5$,
\begin{align}
& \sum_{|\alpha|=3}\left(\frac{\mm{d}}{\mm{d}t}
\tilde{\mathcal{E}}_1(\partial^{\alpha}_{\mm{h}}\eta,\partial^{\alpha}_{\mm{h}}u,
\partial^{\alpha}_{\mm{h}}\vartheta) +
\tilde{\mathcal{D}}_1(\partial^{\alpha}_{\mm{h}}\eta,\partial^{\alpha}_{\mm{h}}u,
\partial^{\alpha}_{\mm{h}}\vartheta )\right)\nonumber  \\
&\leqslant c_1^H  \sqrt{\mathcal{E}^L} (\|(\eta,u)\|_7^2+\mathcal{D}^H).\label{basdicl1syy}\end{align}
Then we deduce from \eqref{badiseqin30} for $2\leqslant j\leqslant 3$, \eqref{201902201632} with $ j =6$, \eqref{dfifessim} for $i=4$ and \eqref{basdicl1syy} that
there is constant $c_2^H$,  such that, for any $c_2^l>1$,
\begin{align}
\frac{\mm{d}}{\mm{d}t}\tilde{\mathcal{E}}^H_2 +\tilde{\mathcal{D}}^H_2 \leqslant & \frac{2 R}{P_\vartheta}\sum_{2\leqslant j \leqslant 3}c_{i}^l \int\nabla \partial_t^{j} \eta_3\cdot \nabla \partial_t^{j }\vartheta  \mm{d}y+  R\sum_{|\alpha|=6}\int \partial_{\mm{h}}^\alpha u_3
 \partial_{\mm{h}}^\alpha \vartheta\mm{d}y\nonumber \\
 &+ c_2^H\left(\|  \eta\|_{6,0}^2 +\|  (u_{t}, \vartheta) \|_{4}^2
+ (1+ c_2^l ) \sqrt{\mathcal{E}^L}(\|(\eta,u)\|_7^2+\mathcal{D}^H) \right),
 \label{basdicl1ssdfafda}\end{align}
where
$$ \begin{aligned}
\tilde{\mathcal{E}}^H_2:=&c_2^l \sum_{|\alpha|=3}
\tilde{\mathcal{E}}_1(\partial^{\alpha}_{\mm{h}}\eta,\partial^{\alpha}_{\mm{h}}u,
\partial^{\alpha}_{\mm{h}}\vartheta)
+\frac{1}{2}\left(Q \|\partial_{3} \eta\|_{6,0}^2+\| u\|^2_{6,0}
+ 2\|\eta\|_{6,*}^2\right)  \\
& + \frac{1}{2}\sum_{2\leqslant j \leqslant 3}c_{i}^l\left( \| \partial_t^{j} u\|_{0}^2+E_1(
\partial_t^{j-1}u,\partial_t^{j}\vartheta)   +2 \int  \nabla \partial_t^{j-1}q\cdot D^{t,j}_u\mm{d}y \right),\\
\tilde{{\mathcal{D}}}^H_2:=&c_2^l  \sum_{|\alpha|=3} \tilde{{\mathcal{D}}}_1(\partial^{\alpha}_{\mm{h}}\eta,\partial^{\alpha}_{\mm{h}}u,
\partial^{\alpha}_{\mm{h}}\vartheta)  +  \|\nabla  u\|_{6,0}^2
\\ &+
   \|  (\eta ,u)\|_{6}^2 +\|\nabla q\|_4^2  + \sum_{2\leqslant j \leqslant 3}c_{i}^l \|\nabla \partial_t^j( u, \vartheta)\|_0^2\mbox{ and }
c_3^l=1.
\end{aligned}$$

Noting that
\begin{align}
&\sum_{k=2}^3 \int \nabla \partial_t^{k-1}q\cdot D^{t,k}_u\mm{d}y
\lesssim  (\mathcal{E}^H)^{3/2},\ \Xi_1(u_t)+\Xi_1(u_{tt})\lesssim  (\mathcal{E}^H)^{3/2}, \nonumber\\
& \sum_{|\alpha|=3}\Xi_2(\partial^{\alpha}_{\mm{h}}\eta,\partial^{\alpha}_{\mm{h}}u)\lesssim \|\eta\|_7\mathcal{E}^H\mbox{ or }\|\eta\|_3\|\eta\|_7^2+\sqrt{\mathcal{E}^H}\mathcal{D}^H,\nonumber
\end{align}
thus, similarly to \eqref{201903182050} and \eqref{xx201903182050}, we easily derive that, for any sufficiently large $c_2^l$,
\begin{align}
&    c_2^l\sum_{|\alpha|=3}( \|\nabla\partial^\alpha_{\mm{h}} \eta\|_{\underline{3},0}^2+\|\partial^\alpha_{\mm{h}} ( u,\vartheta)\|_{\underline{2},0}^2)
  +  \|\partial_3\eta\|_{6,0}^2+\|u\|_{6,0}^2+  \|\eta\|_{6}^2  \nonumber  \\
  & \sum_{2\leqslant j \leqslant 3}c_{i}^l  \| (\partial_t^{j} u,
\partial_3\partial_t^{j-1}u,\partial_t^{j}\vartheta)\|_0^2  -c_3^H \mathcal{E}^H(\sqrt{\mathcal{E}^H}+\|\eta\|_7)
\lesssim  \tilde{\mathcal{E}}^H\lesssim \mathcal{E}^H,\label{201903182050xx} \\
&\tilde{\mathcal{D}}^H_3:= c_2^l\sum_{|\alpha|=3}(\|\partial_3\partial^\alpha_{\mm{h}} \eta\|_{\underline{3},0}^2+
\|\partial^\alpha_{\mm{h}} (u,\vartheta)\|_{\underline{2},1}^2 )+    \|\nabla  u\|_{6,0}^2  +
   \|  (\eta ,u)\|_{6}^2 +\|\nabla q\|_4^2  \nonumber
\\ &
\qquad  \quad+ \sum_{2\leqslant j \leqslant 3}c_{i}^l \| \partial_t^j( u, \vartheta)\|_1^2  -c_3^H(\|\eta\|_3\|\eta\|_7^2+\sqrt{\mathcal{E}^H}\mathcal{D}^H)
\lesssim  \tilde{\mathcal{D}}^H\lesssim \mathcal{D}^H,\label{xx2019014493182050xx}
\end{align}
where the constant $c_3^H$ depends on $c_2^l$.
In addition,
making use of \eqref{omdm122nsdfsf}, \eqref{highestdids},  \eqref{highestdidsxxxxxxx} and \eqref{highestdidsxxxx} and the interpolation inequality,
\begin{align}
& \sum_{k=0}^2(\| \partial_t^k u\|_{6-2k}^2+ \|\nabla\partial_t^k q\|_{4-2k}^2
+ \| \partial_t^k \vartheta\|_{6-2k}^2)\nonumber
  \\
  &\lesssim  \|\eta\|_6^2+\  \| (u,u_t,u_{tt}, \partial_t^3 u,   \partial_t^3\vartheta)\|_0^2 + {\mathcal{E}^L} {E}^H,\label{201903011519} \\
&  \sum_{k=1}^2( \| \partial_t^k u\|_{7-2k}^2+\|\nabla\partial_t^k  q\|_{5-2k}^2)+  \sum_{k=0}^2  \| \partial_t^k \vartheta\|_{7-2k}^2\nonumber
  \\
  &\lesssim  \|u\|_5^2+ \|( u_t,  u_{tt},\partial_t^3 u, \vartheta_{tt}, \partial_t^3 \vartheta)\|_1^2+\mathcal{E}^H\mathcal{D}^H,\label{201903011519dfsafsd}\\
&\|u_t\|_4^2\leqslant \|u_t\|_2^{2/3}\|u_t\|_5^{4/3}\lesssim \|u_t\|_2^{2/3} (\|u\|_5+ \|( u_t,  \partial_t^3 u, \vartheta_{t},  \vartheta_{tt}, \partial_t^3 \vartheta)\|_1)^{4/3}+\mathcal{E}^H\mathcal{D}^H.\label{xx201903182050xxxx}
\end{align}

Similarly to \eqref{basdicl1sxxsdafa}, we immediately derive from \eqref{basdicl1syy},  \eqref{xx2019014493182050xx}  and \eqref{xx201903182050xxxx} that,  for properly large $c_2^l$,
\begin{align}
\frac{\mm{d}}{\mm{d}t}\tilde{\mathcal{E}}^H +c_4^H\tilde{\mathcal{D}}^H \lesssim   & \sqrt{\mathcal{E}^L} \|(\eta,u)\|_7^2+\sqrt{\mathcal{E}^H} \mathcal{D}^H,\label{basdicl1sxxsdafasdafa}\end{align}
where $\tilde{\mathcal{E}}^H:=(c_2^l)^2\tilde{\mathcal{E}}^L+ \tilde{\mathcal{E}}^H_2$ and $\tilde{\mathcal{D}}^H:=  (c_2^l)^2\mathcal{D}^L+\tilde{\mathcal{D}}^H_3+
c_3^H(\|\eta\|_3\|\eta\|_7^2+\sqrt{\mathcal{E}^H}\mathcal{D}^H)  $.
In addition, we also derive from \eqref{201903182050xx}--\eqref{201903011519dfsafsd} that
$\mathcal{E}^H$ and $\mathcal{D}^H$ are equivalent to $\tilde{\mathcal{E}}^H$ and $\tilde{\mathcal{D}}^H$ resp..
Consequently we further obtain \eqref{hgheefdsstim} from \eqref{basdicl1sxxsdafasdafa}. \hfill$\Box$
\end{pf}

\subsection{Stability estimate}\label{sec:02}

Now we are in the position to derive the \emph{a priori} stability estimate \eqref{1.19}. We begin with estimating the terms
$\mathcal{G}_1,\cdots,\mathcal{G}_4$.

Using \eqref{dfifessim} with $i=5$, and recalling the equivalence of $\|\eta\|_{7,*}^2$ and $\|\eta\|_7^2$, we deduce that
\begin{align}
\|\eta\|_7^2\lesssim& \|\eta^0\|_7^2e^{-t}+\int_0^te^{-(t-\tau)}( \mathcal{E}^H(\tau)+\mathcal{D}^H(\tau))\mm{d}\tau\nonumber \\
\lesssim & \|\eta^0\|_7^2e^{-t}+\sup_{0\leqslant \tau\leqslant t}\mathcal{E}^H(\tau)\int_0^te^{-(t-\tau)}\mm{d}\tau
+\int_0^t\mathcal{D}^H(\tau)\mm{d}\tau  \nonumber \\
\lesssim &\|\eta^0\|_7^2e^{-t}+\mathcal{G}_3(t), \nonumber
\end{align}
which yields
\begin{equation}
\label{etasef}
\mathcal{G}_1(t)\lesssim \|\eta^0\|_7^2 +\mathcal{G}_3(t).
\end{equation}

Multiplying \eqref{dfifessim} with $i=5$ by $(1+t)^{-3/2}$, we get
\begin{equation*} \frac{\mm{d}}{\mm{d}t}\frac{
\|\eta\|_{7,*}^2}{(1+t)^{3/2}}+
\frac{3}{2}\frac{\|\eta\|_{7,*}^2}{(1+t)^{5/2}}+\frac{\|(\eta,u)\|_{7}^2}{(1+t
)^{3/2}}\lesssim \frac{ \mathcal{E}^H
+\mathcal{D}^H }{(1+t)^{3/2}},\end{equation*}
which implies that
\begin{equation}
\label{etasef1}\mathcal{G}_2(t)\lesssim  \|\eta^0\|_7^2+\mathcal{G}_3(t).\end{equation}

An integration of \eqref{hgheefdsstim} with respect to $t$ gives
$$
\mathcal{G}_3(t)\lesssim \mathcal{E}^H(0)+
\int_0^t\sqrt{\mathcal{E}^L(\tau)}\|(\eta,u)(\tau)\|_{7}^2\mm{d}\tau.
$$
Let $$\mathcal{G}_5(t):=\mathcal{G}_1(t)+ \sup_{0\leqslant \tau\leqslant t}\mathcal{E}^H(\tau)+\mathcal{G}_4(t).$$
From now on, we further assume
$\sqrt{\mathcal{G}_5(T)}\leqslant \delta$ which is a stronger requirement than \eqref{aprpioses}.
Thus, we make use of \eqref{etasef1} to find that
$$\begin{aligned}\mathcal{G}_3(t) \lesssim &\mathcal{E}^H(0)+
\int_0^t {\delta}(1+\tau)^{-3/2}\|(\eta,u)(\tau)\|_{7}^2\mm{d}\tau\\
\lesssim &\mathcal{E}^H(0)+  {\delta}\left( \|\eta^0\|_7^2+\mathcal{G}_3(t)\right),
\end{aligned}$$
which implies
\begin{equation}
\label{G3testim}
\mathcal{G}_3(t) \lesssim \|\eta^0\|_7^2+\mathcal{E}^H(0).
\end{equation}

Finally, we show the time decay behaviour of $\mathcal{G}_4(t)$. Noting that $\mathcal{E}^L$ can be controlled by $\mathcal{D}^L$
except the term $\|\eta\|_{4,0}$ in $\mathcal{E}^L$. To deal with $\|\eta\|_{4,0}$, we use interpolation inequality to get
$$ \|\eta\|_{4,0}\lesssim \|\eta\|_{3,0}^{\frac{3}{4}}\|\eta\|_{7 ,0}^{\frac{1}{4}}.$$
On the other hand, we combine \eqref{etasef} with \eqref{G3testim} to get
$${\mathcal{E}}^{L}+\|\eta\|_{7 ,0}^2\lesssim \tilde{\mathcal{E}}^{L}+\|\eta\|_{7 ,0}^2\lesssim \|\eta^0\|_7^2+\mathcal{E}^H(0) .$$
Thus,
$$\tilde{\mathcal{E}}^L\lesssim{\mathcal{E}}^L\lesssim (\mathcal{D}^L)^{\frac{3}{4}}({\mathcal{E}}^{L}+\|\eta\|_{7 ,0}^2)^{\frac{1}{4}}
\lesssim (\mathcal{D}^L)^{\frac{3}{4}} (\|\eta^0\|_7^2+\mathcal{E}^H(0))^{\frac{1}{4}}. $$
Putting the above estimate into the lower-order energy inequality \eqref{emdslds}, we obtain
$$ \frac{\mm{d}}{\mm{d}t} \tilde{\mathcal{E}}^{L}+
 \frac{ (\tilde{\mathcal{E}}^{L})^{\frac{4}{3}}}{\mathcal{I}_0 ^{1/{3}}}\lesssim 0,$$
which yields
$${\mathcal{E}}^{L}\lesssim \tilde{\mathcal{E}}^{L}\lesssim \frac{\mathcal{I}_0}{\left((\mathcal{I}_0/\mathcal{E}^{L}(0))^{1/3}+  t/3\right)^3}
\lesssim  \frac{ \|\eta^0\|_7^2+ \mathcal{E}^H(0)}{ 1 +t^3} $$
with $\mathcal{I}_0:=c( \mathcal{E}^H(0)+ \|\eta^0\|_7^2)$ for some positive constant $c$. Therefore,
\begin{equation}
\label{etasef12}\mathcal{G}_4(t)\lesssim \|\eta^0\|_7^2 +  \mathcal{E}^H(0). \end{equation}

Now we sum up the estimates \eqref{etasef}--\eqref{etasef12} to conclude that
\begin{equation*}
\mathcal{G}(t): =\sum_{k=1}^4\mathcal{G}_k(t)\lesssim  \|\eta^0\|_7^2+\mathcal{E}^H(0)\lesssim \|\eta^0\|_7^2+\|(u^0, \vartheta^0)\|_6^2,
\end{equation*}
where \eqref{esdoval} has been also used.
Consequently, we have proved the following \emph{a priori} stability estimate.
\begin{pro}  \label{125pro:0401}
Let $(\eta,u,\vartheta)$ be a solution of the transformed MB problem \eqref{0101xxxx}--\eqref{defineditconxx}
with an associated perturbation pressure $q$. Then there is a sufficiently small $\delta_1^{\mm{s}}$, such that $(\eta,u,\vartheta,q)$ enjoys
the following \emph{a priori} stability estimate
\begin{equation}\label{decadsyeste}
\mathcal{G}(T_1)\lesssim  \|\eta^0\|_7^2+\|(u^0, \vartheta^0)\|_6^2,
\end{equation}
provided that $\sqrt{\mathcal{G}_5(T_1)}\leqslant \delta_1^{\mm{s}}$ for some $T_1>0$.
\end{pro}
Proposition \ref{125pro:0401}, together with the following local existence of unique solution of transformed MB problem, immediately yields Theorem \ref{thm3}.
\begin{pro} \label{pro:0401n}
There is a sufficiently small $\delta_2^{\mm{s}}$, such that for any given initial data $(\eta_0,u_0,\vartheta_0)\in H^{7,1}_0\times H^6_0\times H^6_0$ satisfying
\begin{equation*}
\sqrt{\|\eta^0\|_7^2+\|(u^0,\vartheta^0)\|_6^2}\leqslant \delta_2^{\mm{s}}
\end{equation*}
and the compatibility conditions (i.e., $\partial_t^i u(y,0)|_{\partial\Omega}=0$ and $\partial_t^i \vartheta(y,0)|_{\partial\Omega}=0$, $i=1,2$),
there exist a $T_2:=T_2(\delta_2^{\mm{s}})>0$ (depends on $\delta_2^{\mm{s}}$, the domain $\Omega$ and known physical parameters),
and a unique classical solution $(\eta, u,\vartheta)\in C^0([0,T_2],H^{7,1}_0\times H^6_0
\times H^6_0)$ to the transformed MB problem
\eqref{0101xxxx}--\eqref{defineditconxx} with an associated perturbation pressure $q$. Moreover,
$(\partial_t^j u ,\partial_t^j\vartheta)\in C^0([0,T_2],H^{6-2j})$ for $1\leqslant j\leqslant 3$, $\nabla q\in C^0([0,T_2],H^4)$, $\mathcal{E}^H(0)\lesssim \|\eta^0\|_7^2+\|(u^0,\vartheta^0)\|_6^2$, and
$$ \sup_{0\leqslant \tau\leqslant T_2}\left(\|\eta(\tau)\|_7^2+\mathcal{E}^H(\tau)\right)+\int_0^{T_2}\left( \mathcal{D}^H(\tau)+\|u(\tau)\|_7^2+\| \nabla q(\tau)\|_5^2\right) \mm{d}\tau<\infty ,$$
and $\mathcal{G}(t)$ is continuous on $[0,T_2]$.
\end{pro}
%\begin{rem} In view of the definition of $\mathcal{G}(t)$, we have \begin{equation}\label{Gestsim}\mathcal{G}(0)\leqslant  \|\eta_0\|^2_7+2\mathcal{E}^H(0)\leqslant  C_2\left(\|\eta_0\|_7^2+\|(u_0, \vartheta_0)\|_6^2\right), \end{equation} where $C_2\geqslant 1$ denotes a constant depending on the domain $\Omega$ and other physical parameters in the transformed MB equations.\end{rem}
\begin{pf} The transformed MB problem
is very similar to the surface wave problem (1.4) in \cite{GYTILW1}. Moreover, the current problem is indeed simpler than
the surface wave problem due to the non-slip boundary condition $(u,\vartheta)|_{\partial\Omega}=0$. Using the standard method in \cite{GYTILW1},
one can easily establish Proposition \ref{pro:0401n}, hence we omit its proof here.
\hfill $\Box$
\end{pf}

\subsection{Analysis of stability condition}\label{201903102308}

This section is devoted to verifying that the stability criterion can be satisfied, if $Q$ satisfies \eqref{201903111923}.

By the definition of $\Upsilon_1$, and Wirtinger's and Young's inequalities, we have
$$\Upsilon_1\leqslant \sup_{(\varphi,\phi)\in H^1_\sigma\times L^2}\frac{4R\|\varphi_3\|_0\|\phi\|_0}{(Q\pi^2 +2R^2P_\vartheta^{-1})
\|\varphi_3\|_0^2+ P_\vartheta\|\phi\|_0^2}
\leqslant \frac{2R}{\sqrt{  {2  R^2}+{Q\pi^2}{P_\vartheta}    }}.$$
Obviously, $\Upsilon_1< 1$, if $Q$ satisfies
$$\frac{2R^2}{\pi^2P_\vartheta}< Q .$$

By Wirtinger's inequality,
$$
\begin{aligned}
\Upsilon_2\leqslant &\sup_{(\varphi,\psi,\phi)\in H^4_\sigma\times H_\sigma^3\times H^3_0}
\frac{ \|\psi\|_{\underline{3},0}^2 + (R+4R  /{P_\vartheta})\|\nabla \varphi_3\|_{\underline{2},0}\|\nabla \phi\|_{\underline{2},0}
}{Q\|\partial_{3}\varphi\|_{\underline{3},0}^2
 +  2 \|\nabla ( \psi,\phi)\|_{\underline{2},0}^2}\\
<& \sup_{(\varphi,\psi,\phi)\in H^4_\sigma\times H_\sigma^3\times H^3_0} \frac{ \|\nabla \psi\|_{\underline{2},0}^2 +
Q\|\partial_3 \varphi_3\|_{\underline{3},0}^2  +
(R +4R/{P_\vartheta})^2   \|\nabla \phi \|_{\underline{2},0}^2/4Q
}
{
Q\|\partial_{3} \varphi\|_{\underline{3},0}^2
 +2  \|\nabla ( \psi,\phi)\|_{\underline{2},0}^2} .
\end{aligned}
$$
Hence we immediately observe that $\Upsilon_2<1$, if $$ Q> {R^2(1+4{P_\vartheta}^{-1})^2  }/{8 } . $$

Since
$$\frac{2R^2}{\pi^2P_\vartheta}<\frac{2R^2}{P_\vartheta}\leqslant \frac{R^2(1+4{P_\vartheta}^{-1})^2   }{8},$$
we immediately see that  the stability criterion can be satisfied, if $Q$ satisfies \eqref{201903111923}.

\section{Proof for instability of the transformed MB problem}\label{sec:201903102311}

This section is devoted to the proof of instability of transformed MB problem in Theorem \ref{thmsdfa3}. Next we will complete the proof by four subsections.

\subsection{Linear instability}\label{201903082005}

To begin with, we exploit modified variational method of PDE as in \cite{JFJSO2014} to prove the existence of unstable solutions of the linearized zero-MB problem with parameter $\tau$.
\begin{pro}\label{thm:0201201622}
Let  $R>R_0$ and $\tau\in (0,1]$.
If $Q$ satisfies \eqref{201903052059},
 then zero solution
 is unstable to the linearized MB problem with parameter $\tau$:
 \begin{equation}   \left\{\begin{array}{ll}
\eta_t=u &\mbox{ in } \Omega ,\\[1mm]
  u_t-\tau \Delta u+ \nabla q=Q\partial_3^2\eta + R\vartheta e_3&\mbox{ in }  \Omega ,\\[1mm]
  P_{\vartheta} \vartheta_t   - \Delta \vartheta = R u_3& \mbox{ in }  \Omega ,\\
\mathrm{div} {u} =0& \mbox{ in }  \Omega ,\\
 (\eta,u,\vartheta)=0 &\mbox{ on }\partial\Omega ,\\
 (\eta,u,\vartheta)|_t= (\eta^0,u^0,\vartheta^0)  &\mbox{ on }\Omega .\end{array}\right.\label{201803101401}\end{equation}
 That is, there is an unstable solution
$(\eta, u,\vartheta,  q):=e^{\Lambda t}(w/\Lambda,w,\theta ,\beta)$
 to  the above problem, where
 \begin{equation}
 \label{201903092042}
 (w,\theta,\beta )\in H^\infty_\sigma\times  H^\infty_0\times \underline{H}^\infty
 \end{equation}
 solves  the boundary value problem:
 \begin{equation}
 \label{201903081551}
\left\{    \begin{array}{l}
   \Lambda  w =\tau \Delta w
-  \nabla \beta +R\theta e_3 +   Q
 \partial_3^2 w/\Lambda,\\
     \Lambda P_{\vartheta} \theta = \Delta  \theta+ R w_3 ,\\
   \mathrm{div}w=0, \\
 (w, \vartheta)|_{\partial\Omega}= {0}  \end{array}  \right.
                       \end{equation}
  with a
growth rate  $\Lambda(\tau)>0$ satisfying
\begin{equation}\label{Lambdard}
\begin{aligned}
\Lambda (\tau)= \sup_{(\varpi,\phi)\in\mathcal{A} }F( \varpi,\phi,\Lambda (\tau),\tau). \end{aligned}\end{equation}
Moreover,
\begin{align}
\label{201602081445MH}
&\|w_{\mm{h}}\|_{L^1}\|w_3\|_{L^1}  \|\partial_3w_{\mm{h}}\|_{L^1}\|\partial_3w_3\|_{L^1}
\|\vartheta \|_{L^1}\|\partial_3\vartheta\|_{L^1}>0 . \end{align}
\end{pro}
\begin{rem}
We call $ \Lambda_0$ the largest growth of thermal instability of the linearized Rayleigh--B\'enard problem  \begin{equation} \label{201903092026}  \left\{\begin{array}{ll}
  u_t-\Delta u+ \nabla q=  R\vartheta e_3, \ \mathrm{div} {u} =0&\mbox{ in }  \Omega ,\\[1mm]
  P_{\vartheta} \vartheta_t   - \Delta \vartheta = R u_3& \mbox{ in }  \Omega ,\\
 ( u,\vartheta)|_{t=0}= ( u^0,\vartheta^0) &\mbox{ in } \Omega ,\\
 ( u,\vartheta)=0 &\mbox{ on }\partial\Omega ,
\end{array}\right.\end{equation}
since any classical solution of \eqref{201903092026} enjoys the estimate
$$ \|(u,\vartheta)\|_0^2\lesssim\|(u^0,\vartheta^0)\|_0^2e^{\Lambda_0 t}, $$
and there exists a classical solution $(u,\vartheta,q)=(\tilde{w},\tilde{\theta},\tilde{\beta})e^{\Lambda_0t}$ to \eqref{201903092026} with $(\tilde{w},\tilde{\theta},\tilde{\beta})$ satisfies the regularity
 \eqref{201903092042} and the property \eqref{Lambdard} with $(\tau,Q)=(1,0)$.
\end{rem}
\begin{rem}
\label{201903102310}
We further remark the upper and lower bounds for $\Lambda_0$.
Noting that there exists  $(w^1,\theta^1)\in \mathcal{A}$ such that
$\Lambda_0=D_R( w^1,\theta^1,1) >0$.
Obviously $w_3^1\neq 0$, thus $w_{\mm{h}}^1\neq 0$ due to $w^1|_{\partial\Omega}=0$ and $\mm{div}w^1=0$.
By Young's inequality and the fact $\|w_3^1\|_0<\|w^1\|_0$, we get an upper bound for $\Lambda_0$:
\begin{equation}
\label{201903101512}
\Lambda_0\leqslant R(\|w_3^1\|_0^2+P_\vartheta\|\theta^1\|_0^2)/\sqrt{P_\vartheta}<R/\sqrt{P_\vartheta}.
\end{equation}

By the definition of $R_0$, for any $(u,\vartheta)\in \mathcal{A}$,
$$D_{R_0}( w,\theta,1):=  2 R_0\int u_3\vartheta\mm{d}y-\|\nabla (u,\vartheta)\|_0^2 \leqslant 0,$$
thus we get another upper bound for $\Lambda_0$:
\begin{equation*}
\Lambda_0 \leqslant 2(R-R_0)  \int w_3^1\theta^1\mathrm{d}y.
 \end{equation*}
 In particular, we see that $\Lambda_0$ as $R\to R_0$.

In addition, by the definition of $R_0$, there exists a  $(u^2,\vartheta^2)\in \mathcal{A}$ such that $D_{R_0}( w^2,\theta^2,1)=0$,
thus we get a lower bound for $\Lambda_0$:
\begin{equation*}
2(R-R_0)\int u_3^2\vartheta^2\mm{d}y\leqslant
 \sup_{(\varpi,\phi)\in\mathcal{A}} D_R(\varpi,\phi,1)=\Lambda_0,
 \end{equation*}
 which yields that
 \begin{equation}
\label{201903092105}
2(R-R_0)\xi\leqslant  \Lambda_0.
 \end{equation}
\end{rem}
\begin{pf}
Next we divide the proof of Proposition \ref{thm:0201201622} into four steps.

(1) \emph{Existence of weak solutions to the modified problem:}
 \begin{equation}
\label{2016040614130843xd}
\left\{    \begin{array}{ll}
\alpha(s,\tau) w= \tau\Delta w
-  \nabla \beta +R\theta e_3 +   Q
 \partial_3^2 w/s ,\\
    \alpha(s,\tau) P_{\vartheta} \theta  = \Delta  \theta+ R w_3 ,\\
   \mathrm{div}w=0, \\
  (w,\theta)|_{\partial\Omega} = {0},  \end{array} \right.
                       \end{equation}
\emph{where $\tau\in (0,1]$ and $s>0$  are given parameters.}

To prove the existence of weak solutions of the above problem, we consider the variational problem of the functional $F(\varpi,\phi,s ,\tau)$:
\begin{equation}
\label{201903111430}
\alpha(s,\tau): =\sup_{(\varpi,\phi)\in\mathcal{A}} F(\varpi,\phi,s,\tau ).
\end{equation}
For simplicity, next we temporarily denote $\alpha(s,\tau)$ and $F(\varpi,\phi,s,\tau)$
by $\alpha$  and $F(\varpi,\phi)$, resp..

Using the condition $J(w^n,\theta^n)=1$, by Young's inequality, we easily see that
\begin{equation}
\alpha\leqslant R/\sqrt{P_\vartheta}.
\label{201903052117}
\end{equation}
Hence, $F(\varpi,\phi)$ has a maximizing sequence $\{(w^n,\theta^n)\}_{n=1}^\infty\subset \mathcal{A}$, which satisfies $$\alpha=\lim_{n\to\infty} F(w^n,\theta^n).$$
Moreover,  $\|(w^n,\theta^n)\|_1 \leqslant c$ for some constant $c$, which is independent of $n$. Thus  there exists a subsequence, still labeled
 by $(w^n,\theta^n)$, and a function $(w,\theta)\in \mathcal{A}$, such that
$$\begin{aligned}
& (w^n,\theta^n)\rightharpoonup (w,\theta)\mbox{ in }H^1_\sigma\times H^1_0\mbox{ and } (w^n,\theta^n)\to (w,\theta)\mbox{ in }L^2.
\end{aligned}$$
Exploiting the above convergence results, and the lower semicontinuity of weak convergence, we find that
\begin{equation*}
\alpha =\limsup_{n\to \infty}F(w^n,\phi^n)
\leqslant F(w,\theta) \leqslant \alpha .\end{equation*}
Hence $ (w,\theta)$
 is a maximum point of the functional $F(\varpi,\phi)$ with respect to $(\varpi,\phi)\in \mathcal{A}$.

Obviously, $w$ constructed above is also  the maximum point of the functional $F(\varpi,\phi )/  J(\varpi,\phi)$ with respect to
$(\varpi,\phi)\in H_{\sigma}^1\times H_0^1$. Moreover, $\alpha= F(w,\theta)/J(w,\theta)$. Thus, for any given $(\varphi,\phi) \in H^1_{\sigma} $, the point $t=0$ is the maximum point of the functional
$$I(t):=F(w+t\varpi ,\theta+t \phi )- \alpha \int J(w+t\varpi, \theta+t\phi)\mm{d}y\in C^1(\mathbb{R}).$$
Then, by evaluating $I'(0)=0$, we obtain the weak form:
\begin{align}
& R\int (w_3 \phi+\theta\varpi_3)\mm{d}y
-\tau\int \nabla w: \nabla\varpi\mm{d}y   -\int \nabla \theta: \nabla \phi\mm{d}y  -
\frac{Q}{s}\int \partial_3 w\cdot \partial_3 \varpi \mm{d}y
\nonumber \\
 & = \alpha \int ( w\cdot \varpi +P_\vartheta \theta \phi ) \mm{d}y.
\label{201805161108}
\end{align}
This means that $(w,\theta)$ is a weak solution of the modified problem \eqref{2016040614130843xd}.

(2) \emph{Improvement of regularity of the weak solution $(w,\theta)$.}

To begin with, we establish the following preliminary conclusion:

\emph{For any $i\geqslant 0$, we have}
\begin{equation}
(w,\theta) \in H^{1,i}_\sigma \times H^{1,i}_0
\label{201806181000}
\end{equation}
\emph{and}
\begin{align}
&R\int (\partial_{\mm{h}}^{i} w_3\phi+\partial_{\mm{h}}^{i}\theta\varpi_3)\mm{d}y-\tau\int \nabla \partial_{\mm{h}}^i w :\nabla\varpi \mm{d}y -
 \int  \nabla \partial_{\mm{h}}^{i} \theta:\nabla \phi \mm{d}y
\nonumber \\
& -
\frac{Q}{s}\int \partial_3 \partial_{\mm{h}}^i w\cdot \partial_3\varpi\mm{d}y=  \alpha \int (\partial_{\mm{h}}^i w\cdot   \varpi+P_\vartheta\partial_{\mm{h}}^i \theta\phi ) \mm{d}y.
\label{201805161108sdfs}
\end{align}

Obviously, by induction, the above assertion reduces to verification of the following recurrence relation:

\emph{For given $i\geqslant 0$, if $(w, \theta) \in H^{1,i}_{\sigma}\times H^{1,i}_{0}$ satisfies \eqref{201805161108sdfs}
for any $(\varpi, \phi) \in H^1_{\sigma } \times H^1_{0} $,  then}
\begin{equation}
\label{201806181412}
(w, \theta) \in H^{1,i+1}_{\sigma }\times H^{1,i+1}_{0 },
\end{equation}\emph{ and
 $(w,\theta)$ satisfies}
\begin{align}
&R\int (\partial_{\mm{h}}^{i+1}w_3\phi+\partial_{\mm{h}}^{i+1}\theta\varpi_3)\mm{d}y - {\tau} \int  \nabla \partial_{\mm{h}}^{i+1} w :\nabla  \varpi \mm{d}y - \int  \nabla \partial_{\mm{h}}^{i+1}\theta\cdot\nabla \phi \mm{d}y\nonumber \\
& -\frac{Q}{s}
\int \partial_3 \partial_{\mm{h}}^{i+1} w\cdot \partial_3\varpi \mm{d}y   = \alpha  \int  (\partial_{\mm{h}}^{i+1} w\cdot   \varpi + P_\vartheta \partial_{\mm{h}}^{i+1} \theta \phi )\mm{d}y.\label{201805161108sdfssafas}
\end{align}
Next we verify the above recurrence relation by the method of difference quotients as in \cite{JFJSOMITNWZ}.

We define the difference quotient $D_j^h f:=(f(y+h e_j)-f(y))/h$ for $j=1$ and $2$, and $D^h_{\mm{h}}f:=(D_1^hf,D_2^hf)$, where $h$ is a constant. Assume that $(w,\theta) \in H^{1,i}_{\sigma }\times H^{1,i}_{0}$ satisfies \eqref{201805161108sdfs}, then we can deduce from \eqref{201805161108sdfs}
that
\begin{align}
& R\int ( \partial_{\mm{h}}^i w_3 D_j^h\phi+\partial_{\mm{h}}^i \theta D_j^h \varpi_3)\mm{d}y-
\tau\int \nabla \partial_{\mm{h}}^i w : \nabla D_j^h \varpi
 \mm{d}y- \int \nabla \partial_{\mm{h}}^i \theta\cdot \nabla D_j^h \phi
 \mm{d}y\nonumber \\
& - \frac{Q}{s}\int \partial_3\partial_{\mm{h}}^i w \cdot \partial_3 D_j^h \varpi  \mm{d}y
= \alpha  \int(
  \partial_{\mm{h}}^i w\cdot  D_j^h \varpi+ P_\vartheta
  \partial_{\mm{h}}^i \theta  D_j^h \phi )  \mm{d}y \nonumber
\end{align}
and
\begin{align}
& R\int (\partial_{\mm{h}}^i w_3 D_j^{-h} D_j^h \partial_{\mm{h}}^i\phi+\partial_{\mm{h}}^i\theta D_j^{-h} D_j^h \partial_{\mm{h}}^i\varpi_3)\mm{d}y
-\tau\int \nabla \partial_{\mm{h}}^i w : \nabla D_j^{-h} D_j^h \partial_{\mm{h}}^i w
 \mm{d}y\nonumber \\
& - \int \nabla \partial_{\mm{h}}^i \theta \cdot \nabla D_j^{-h} D_j^h \partial_{\mm{h}}^i \phi
 \mm{d}y-
\frac{Q}{s}\int \partial_3 \partial_{\mm{h}}^i w\cdot \partial_3 D_j^{-h} D_j^h \partial_{\mm{h}}^i \varpi \mm{d}y \nonumber \\
&=  \alpha  \int (
 \partial_{\mm{h}}^i w\cdot D_j^{-h}  D_j^h  \partial_{\mm{h}}^i \varpi+ P_\vartheta \partial_{\mm{h}}^i \theta   D_j^{-h}  D_j^h  \partial_{\mm{h}}^i \phi ) \mm{d}y, \nonumber
\end{align}
which yield
\begin{align}
&R\int ( D_j^{-h} \partial_{\mm{h}}^i w_3\phi + D_j^{-h} \partial_{\mm{h}}^i\theta\varpi_3 )\mm{d}y
-\tau\int \nabla D_j^{-h} \partial_{\mm{h}}^i w :\nabla  \varpi
 \mm{d}y - \int \nabla D_j^{-h}\partial_{\mm{h}}^i \theta\cdot \nabla \phi
 \mm{d}y \nonumber \\
& -\frac{Q}{s}\int \partial_3  D_j^{-h}\partial_{\mm{h}}^i w \cdot \partial_3 \varpi\mm{d}y  = \alpha \int (
  D_j^{-h}\partial_{\mm{h}}^i w\cdot \varpi+  P_\vartheta D_j^{-h}\partial_{\mm{h}}^i \theta \phi)\mm{d}y ,
\label{20180612713459}
\end{align}
and
\begin{align}
&\|\nabla D_j^h\partial_{\mm{h}}^i(\sqrt{\tau} w,\theta)\|^2_0  + Q \|\partial_3  D_j^h\partial_{\mm{h}}^i w \|^2_0/s\leqslant c \|
D_j^h\partial_{\mm{h}}^i (w ,\theta)\|^2_0, \label{201806161512}
\end{align}
resp., where the constant $c$ is independent $h$.

By Friedrichs'  inequality, we
deduce from \eqref{201806161512} that
\begin{align}
  \|D^h_{\mm{h}}  \partial_{\mm{h}}^i (w,\theta)\|^2_1
\lesssim  \|\nabla D^h_{\mm{h}} \partial_{\mm{h}}^i (w,\theta) \|^2_0\lesssim \| \nabla_{\mm{h}} \partial_{\mm{h}}^i ( w,\theta) \|^2_0\lesssim 1.\nonumber
\end{align}
Thus
there is a subsequence of $\{-h\}_{h\in \mathbb{R}}$, still denoted by $-h$, such that
\begin{equation}  \label{201806127345}
 D^{-h}_{\mm{h}} \partial_{\mm{h}}^i  (w,\theta) \rightharpoonup \nabla_{\mm{h}}\partial_{\mm{h}}^i (w,\theta) \mbox{ in }H^1_\sigma\times H^1_0,
 \quad  D^{-h}_{\mm{h}} \partial_{\mm{h}}^i  (w,\theta) \to \nabla_{\mm{h}}\partial_{\mm{h}}^i (w,\theta) \mbox{ in }L^2.
\end{equation}
Employing the regularity of $(w,\theta)$ in \eqref{201806127345} and the fact $(w,\theta)\in H_{\sigma}^{1,i}\cap H_0^{1,i}$, we get \eqref{201806181412}.
In addition, exploiting the limiting results in \eqref{201806127345}, we deduce \eqref{201805161108sdfssafas} from \eqref{20180612713459}.
This completes the proof of the recurrence relation, and thus \eqref{201806181000} and \eqref{201805161108sdfs} hold.

With \eqref{201806181000} in hand, we can consider a Stokes problem:
 \begin{equation}\label{2016040614130843x}      \left\{  \begin{array}{l}
 \nabla \beta^k -(\tau+Q/s)\Delta\omega^k =\partial_{\mm{h}}^k\mathcal{L}^1 ,\\
- \Delta  \sigma^k=\partial_{\mm{h}}^k\mathcal{L}^2 ,\\
\mm{div}\omega^k=0 ,
\\ (\omega^k,\sigma^k)|_{\partial\Omega} =0  ,\end{array}\right.
\end{equation}
where $k\geqslant 0$ is a given integer, and $\mathcal{L}^1:= R \theta e_3-\alpha w - Q\Delta_{\mm{h}} w/s$ and $\mathcal{L}^2:= R  w_3-  \alpha P_{\vartheta} \theta $.
Recalling the regularity \eqref{201806181000} of $w$, we see that $\partial_{\mm{h}}^k(\mathcal{L}^1,\mathcal{L}^2)\in H^1$. By the existence theories of the both of Stokes and elliptic problems, there exists a unique strong solution $(\omega^k,\sigma^k,\beta^k)\in H^3_\sigma\times H_0^3\times\underline{H}^2$ to the above problem \eqref{2016040614130843x}.

Multiply \eqref{2016040614130843x}$_1$ and \eqref{2016040614130843x}$_2$ by $\varpi\in H^1_{\sigma}$ and $\phi\in H^1_{0}$ in $L^2$, reps., and then adding the two identities, we get
\begin{align}
 &R\int (\partial_{\mm{h}}^{i} w_3\phi+\partial_{\mm{h}}^{i}\theta\varpi_3)\mm{d}y- \left(\tau+\frac{Q}{s}\right)\int \nabla \omega^k:\nabla \varpi \mm{d}y
-\int \nabla \sigma^k:\nabla \phi\mm{d}y\nonumber  \\
   &= \alpha\int (\partial_{\mm{h}}^k w\cdot \varpi+P_\vartheta\partial_{\mm{h}}^k \theta \phi )\mm{d}y
+\frac{Q}{s}\int \partial_{\mm{h}}^k  \Delta_{\mm{h}} w\cdot \varpi\mm{d}y.
\label{201808072057}
\end{align}
Thus, subtracting \eqref{201808072057} from  \eqref{201805161108sdfs}, we obtain
$$\int\left(\left(\tau+\frac{Q}{s}\right) \nabla(\partial_{\mm{h}}^k  w-\omega^k):\nabla \varpi
+ \nabla(\partial_{\mm{h}}^k  \theta-\sigma^k)\cdot \nabla \phi\right)\mm{d}y=0.  $$
Taking $(\varpi,\phi):= (\partial_{\mm{h}}^k w-\omega^k, \partial_{\mm{h}}^k \theta-\sigma^k) \in H^1_{\sigma}\times H_0^1$ in the above identity, and then using the Friedrichs'  inequality,
we find that $(\omega^k,\sigma^k)=\partial_{\mm{h}}^k (w,\theta)$. Thus, we immediately conclude that
\begin{equation}
\label{201806181507}
\partial_{\mm{h}}^k ( w,\theta)\in H^{3}\mbox{ for any }k\geqslant 0,
\end{equation}
which implies $\partial_{\mm{h}}^k(\mathcal{L}^1,\mathcal{L}^2)\in H^3$ for any $k\geqslant 0$.
Hence, applying the regularity theory of Stokes problem to \eqref{2016040614130843x}, we get
\begin{equation}
\label{201806181507fdsdsgsdfgsg}
\partial_{\mm{h}}^k  (w,\theta)\in H^{5} \mbox{ for any }k\geqslant 0,
\end{equation}
Obviously, by induction we can easily follow the improving regularity method from \eqref{201806181507} to \eqref{201806181507fdsdsgsdfgsg} to
deduce that $(w,\theta)\in H^\infty_0$. In addition, we have $\beta:=\beta^0\in \underline{H}^\infty$. Moreover, $\beta^k$ in \eqref{2016040614130843x}
is equal to $\partial^k_{\mm{h}}\beta$. Hence we  see that $(w,\theta,\beta)$ constructed above is indeed a classical solution to the modified problem \eqref{2016040614130843xd}.

(3) \emph{Some properties of the function $\alpha(s,\tau)$ on $\mathbb{R}_+$, where $\tau$ is given:}
\begin{align}
\label{201702081047}&
\alpha(s,\tau) \in C^{0,1}_{\mm{loc}}(\mathbb{R}_+),\\
&\label{201702081122n}
\Lambda_0/2< \alpha(s,\tau)\leqslant R/\sqrt{P_\upsilon}  \mbox{ \emph{ on  interval }}[\Lambda_0/2,\infty),
\\
\label{201702081046}
&\alpha(s_1,\tau) <\alpha(s_2,\tau)\mbox{ \emph{ for any }}\Lambda_0/2<s_1<s_2,\mbox{ if }Q\neq 0.
\end{align}
For simplicity, next we temporarily denote $\alpha(s,\tau)$ and $F(\varpi,\phi,s,\tau)$
by $\alpha(s)$  and $F(\varpi,\phi,s)$, resp.. It is easy to see that $\alpha(s)$ is a positive constant for $Q=0$, since $R>R_0$.

We begin with verification of \eqref{201702081047}.
 Choosing a bounded interval $[c_1,c_2]\subset \mathbb{R}_+$, then for any $s\in [c_1,c_2]$,
there exists  $(w^{s}, \theta^{s})$ satisfying $\alpha(s)=F(w^{s}, \theta^{s}, s) $. So, by the monotonicity
\begin{align}\alpha(s_2) \geqslant  \alpha(s_1)\mbox{ for any }s_2\geqslant s_1>0,\nonumber
\end{align}  we have
$$\alpha(c_1)\leqslant \alpha(s) \leqslant \alpha( 2s)- Q\|\partial_3 w^s\|_0^2/2s \leqslant \alpha(2c_2)- Q\|\partial_3 w^s\|_0^2/2c_2 ,$$
which yields
$$Q\|\partial_3 w^s\|_0^2  \leqslant  2c_2(\alpha(2c_2)-\alpha(c_1))=:c_1^2 K\mbox{ for any }s\in [c_1,c_2].$$
Thus for any $s_1$, $s_2\in [c_1,c_2]$,
$$
\begin{aligned}
\alpha(s_1)-\alpha(s_2)\leqslant & F(w^{s_1},\theta^{s_1},s_1)-F(w^{s_1},\theta^{s_1},s_2)\leqslant K|  {s_2}-s_1|
\end{aligned}$$
 and $$\alpha(s_2)-\alpha(s_1)\leqslant K|  {s_2}-s_1|,$$
which immediately imply $|\alpha(s_1)-\alpha(s_2)|\leqslant K| {s_2}-s_1|$. Hence, \eqref{201702081047} holds.

Now we turn to prove \eqref{201702081122n}. For given $s\in [\Lambda_0/2,\infty)$, it is easy to see that
\begin{align}
&\sup_{(\varpi,\phi)\in \mathcal{A}}\left(
D_R(\varpi,\phi,1)-2Q\|\partial_3 \varpi\|_0^2/\Lambda_0  \right)\nonumber \\
& \leqslant \sup_{(\varpi,\phi)\in \mathcal{A}}\left(
D_R(\varpi,\phi,\tau)-2Q\|\partial_3 \varpi\|_0^2/\Lambda_0  \right)=\alpha( \Lambda_0/2 )\leqslant \alpha(s).\label{201903132109}
\end{align}
Moreover, by step (1), there exists $(w^0,\theta^0)\in\mathcal{A}$ such that
\begin{equation}
\label{201903101039}
u^0\neq0,\ \theta^0\neq0 \mbox{ and }\Lambda_0=D_R(w^0,\theta^0,1)>0.
\end{equation}

Recalling the definition of $\Lambda_0$, and using \eqref{201903132109} and \eqref{201903101039}, we have, for $Q\neq0$,
\begin{align}
(1+2Q/\Lambda_0)\Lambda_0<&
D_R(w^0,\theta^0,1)-2Q\|\partial_3 w^0\|_0^2/\Lambda   \nonumber \\
  & +4QR \Lambda_0^{-1} \int w_3^0\theta^0\mm{d}y  \leqslant  \alpha(s)  +2Q R /\Lambda_0 \sqrt{P_\vartheta},
\end{align}
which yields that
$$ \alpha(s)>  \Lambda_0+2Q(1-R/\Lambda_0 \sqrt{P_\vartheta}),  $$
which, together with \eqref{201903052059} and \eqref{201903052117}, yields \eqref{201702081122n} for $Q\neq0$. For $Q=0$, \eqref{201702081122n} obviously holds; moreover, $\Lambda_0<\alpha(s)$ for any given $\tau\in (0,1)$, and $\Lambda_0=\alpha(s)$ for $\tau=1$.

Finally, we verify \eqref{201702081046}.
For given $s_2>s_1\geqslant \Lambda_0/2$, there exists $(w^{s_1}, \theta^{s_1})\in \mathcal{A}$, such that
$$\alpha(s_1)  = F(w^{s_1},\theta^{s_1},s_1)>\Lambda_0/2>0\mbox{ and }w^{s_1}\neq 0.$$ Therefore, by Wirtinger's inequality,  we further have
$$ \alpha(s_1)= F(w^{s_1},\theta^{s_1 },s_1) \leqslant  \alpha(s_2) + Q\left(\frac{ 1}{s_2}- \frac{1}{s_1}\right) \| \partial_3 w^{s_1}\|_0^2<  \alpha(s_2), $$
 which implies \eqref{201702081046}.

(4) \emph{Construction of an interval for fixed point}:

 Now exploiting \eqref{201702081047}--\eqref{201702081046},
we find by a fixed-point argument on $[\Lambda_0/2,\infty)$ that there is a unique $\Lambda(\tau)\in (\Lambda_0/2,\infty)$ satisfying
\begin{equation*}  \Lambda(\tau)= \alpha(\Lambda(\tau)) =
  \sup_{\varpi\in\mathcal{A}}F(\varpi,\phi, \Lambda(\tau),\tau )\in (\Lambda_0/2,\infty).
\end{equation*}
Hence, there is a classical solution $(w,\phi)\in(\mathcal{A}\cap H^\infty)\times H_0^\infty$ to the boundary-value problem \eqref{201903081551}
with $\Lambda(\tau)$ constructed above and with $\beta\in \underline{H}^\infty $.
Moreover,
\begin{equation}
\label{growthnn} \Lambda(\tau)= F(w,\phi, \Lambda(\tau),\tau)>0.
\end{equation}
In addition, \eqref{201602081445MH}  directly follows \eqref{growthnn},  the fact
$(w,\vartheta)\in H_\sigma^1\times H_0^1$, and  \eqref{201903081551}$_2$.
This completes the proof of Proposition \ref{thm:0201201622}.
\hfill $\Box$
\end{pf}

Next we shall establish upper and lower bounds for $\Lambda(1)$, which will be used in the derivation of error estimates.
\begin{lem}
\label{thm:0201201622xx}
Let $\Lambda^*:=\Lambda(1)$ be constructed in Proposition \ref{thm:0201201622} and $Q$ further satisfies $Q\in [0,1)$,   then
\begin{equation}
\label{201903101106}
\Lambda_0-\sqrt{ Q\mathcal{R}}\leqslant \Lambda^*\leqslant  \Lambda(1-\sqrt{Q})\leqslant  (1-\sqrt{Q} )\Lambda_0 +R\sqrt{Q/P_\vartheta}.
\end{equation}
 \end{lem}
 \begin{pf} Obviously, for any $0<\tau_2\leqslant \tau_1\leqslant 1$,
 \begin{equation}
\Lambda_0/2<  \alpha(\tau_1,s)\leqslant\alpha(\tau_2,s)\mbox{ for any }s\geqslant \Lambda_0/2. \nonumber
 \end{equation}
By the construction of \eqref{growthnn}, we easily see
 \begin{equation*}
\Lambda_0/2<  \Lambda(\tau_1)\leqslant \Lambda(\tau_2) ,
  \end{equation*}
  which yields that
\begin{align}
\label{201903221400}
\Lambda_0/2< \Lambda^*\leqslant  \Lambda(1-\sqrt{Q}).
\end{align}

Noting that there exists $(w^0,\phi^0)\in \mathcal{A}$ such that
$$
\Lambda(1-\sqrt{Q})=D_R(w^0,\phi^0, 1-\sqrt{Q} )-Q\|\partial_3 w^0\|_0^2/\Lambda(1-\sqrt{Q}), $$
thus we have
\begin{align}
\Lambda(1-\sqrt{Q})\leqslant &(1-\sqrt{Q})D_R(w^0,\phi^0, 1)+2R \sqrt{Q}\int u_3\vartheta\mathrm{d}y\nonumber \\
\leqslant& (1-\sqrt{Q} )\Lambda_0 +R\sqrt{Q/P_\vartheta}. \label{201903221402}
\end{align}

Using \eqref{201903101039}, we have
\begin{align}
(1+Q/\Lambda^* )\Lambda_0= &(1+Q /\Lambda^{*}) D_R(w^0,\theta^0,1) \nonumber  \\
\leqslant & F(w^0,\theta^0,\Lambda^*,1)+ \frac{2QR}{\Lambda^{*} }\int w_3^0\theta^0\mathrm{d}y
 \leqslant \Lambda^* +\frac{QR}{\Lambda^{*}\sqrt{P_\vartheta}}.
\label{201903221354}
\end{align}
Noting that $\Lambda_0/2< \Lambda^*$, we deduce from \eqref{201903221354} that
\begin{align}
\Lambda^*\geqslant \frac{\Lambda_0+\sqrt{\Lambda_0^2-4Q\mathcal{R}}}{2}\geqslant \Lambda_0-\sqrt{ Q\mathcal{R}} .
\label{201903221403}
\end{align}
Putting \eqref{201903221400}, \eqref{201903221402}
and \eqref{201903221403} together
yields \eqref{201903101106}.
 \hfill $\Box$
 \end{pf}

\subsection{Gronwall-type energy inequality of nonlinear solutions}

This section is devoted to establishing Gronwall-type energy inequality for solutions of the transformed MB problem.
Let $(\eta,u)$ be a solution
of the transformed RT problem, such that
\begin{equation}\label{aprpiosesnew}
 \sup_{0\leqslant t < T}\sqrt{\|\eta(t)\|_5^2+\|(u,\vartheta)(t)\|_4^2}\leqslant \delta \in (0,1)\;\;\mbox{ for some  }T>0
\end{equation}
where $\delta$ is sufficiently small.

Noting that, under the assumptions of  \eqref{2019022201518} and \eqref{aprpiosesnew},
\eqref{aimdse}--\eqref{lemdssdf} and \eqref{201903011350} also holds for  $1\leqslant i\leqslant 3$, $0\leqslant j\leqslant 6-2i$, $0\leqslant k\leqslant 4$.
Thus, following the argument of \eqref{badiseqin30}, \eqref{201902171647}, \eqref{201902201632}, \ref{dfifessim}, and Lemmas \ref{lem:201902181540}, \eqref{201903231042} and \ref{201903081525} with slight modification (or referring to Lemmas 3.3--3.7 in \cite{JFZYYOAR}),  we easily establish the following five lemmas:
\begin{lem}\label{201612132242nn}
Under the assumptions of  \eqref{2019022201518} and \eqref{aprpiosesnew},
                        for $0\leqslant i\leqslant 4$,
\begin{align}
&
\frac{\mm{d}}{\mm{d}t}\int  \left(  \partial_\mm{h}^i \eta \cdot  \partial_\mm{h}^i u  + \frac{ 1}{2}|\nabla \partial_\mm{h}^i \eta |^2\right) \mm{d}y+
Q \|\partial_3\partial_{\mm{h}}^i \eta\|_0^2 \lesssim
\|(  u,\vartheta)\|^2_{i,0}
 + \sqrt{\mathfrak{E}} \mathfrak{D} ,
\label{ssebdaiseqinM0846} \\
&
 \label{201702061418}
 \frac{\mm{d}}{\mm{d}t}(Q\|\partial_3\partial_{\mm{h}}^i \eta\|_0^2+\| \partial_\mm{h}^i ( u,\vartheta)\|^2_0
   )
+ c\|\partial_\mm{h}^i  ( u, \vartheta) \|_{1}^2 \lesssim
\| (u_3,\vartheta)\|_{i,0}^2+ \sqrt{\mathfrak{E} }\mathfrak{D} . \end{align}
\end{lem}

\begin{lem}\label{201612132242nxsfssdfs}Under the assumptions of  \eqref{2019022201518} and \eqref{aprpiosesnew},
\begin{align}
&  \frac{\mm{d}}{\mm{d}t}( \|    u_{t}\|_{0}^2
+E_1(u,\vartheta_t) )
+ c\| (u,\vartheta)_t\|_{1}^2 \lesssim  \|u_3\|_1^2+ \sqrt{\mathfrak{E} } \mathfrak{D} , \label{Lem:030dsfafds1m0dfsf832}\\
&  \frac{\mm{d}}{\mm{d}t}\left(\| u_{tt}\|_{0}^2+E_1(u_t,\vartheta_{tt})
+2\int \nabla q_t\cdot D^{t,2}_{u}\mm{d}y)
\right)
+ c\| (u,\vartheta)_{tt}\|_{1}^2 \lesssim  \|u_t\|_{1}^2+ \sqrt{\mathfrak{E} } \mathfrak{D} . \label{Lem:030dsfafds1m0dfsf832xx}
 \end{align}
\end{lem}
\begin{lem}\label{201612132242nx}
Under the assumptions of  \eqref{2019022201518} and \eqref{aprpiosesnew},
\begin{align}
  &\|u\|_4+\|\nabla q\|_{2}+\|u_t\|_2+\|\nabla q_t\|_{0} \nonumber \\
  &\lesssim \|\eta\|_4+\|\vartheta\|_2+\|u\|_0+ \|(u_{tt},\vartheta_t)\|_0  ,\label{2017020614181721}\\
&\|   u\|_5+ \|\nabla q\|_{3}   \lesssim
 \|\eta\|_{5} + \|\vartheta\|_3 + \|u_{t}\|_3,\label{201903192036xx}\\
&\|  u_t\|_3+\|\nabla q_t\|_1 \lesssim
 \|u\|_3  +\|(u_{tt},\vartheta_t)\|_1+\|u\|_3 \|\nabla q\|_1 ,\label{201903192036}
\\  & \|\vartheta\|_{4}+\| \vartheta_t\|_{2}
 \lesssim  \|u\|_2+  \|(u_t,\vartheta_{tt})\|_0 ,\label{2017020614181721xxxxxxx}\\
&  \| \vartheta \|_5 + \| \vartheta_t \|_3  \lesssim  \| u \|_3+ \|(u_t,\vartheta_{tt})\|_1 .\label{2017020614181721xxxxx}
\end{align}
\end{lem}
\begin{lem}
\label{lem:dfifessim2057}
Under the assumptions of  \eqref{2019022201518} and \eqref{aprpiosesnew}, there exists a functional  $ \|\eta\|_{5,*}$, which is equivalent to $\|\eta\|_5$ and satisfies
 \begin{align}
& \frac{\mm{d}}{\mm{d}t}
 {\|\eta\|}_{5,*}^2+
\left\|  \left( \eta,u\right)\right\|_5^2+
\|\nabla q\|_3^2 \lesssim  \| \eta \|_{\underline{5},0}^2+\|(u_t,\vartheta)\|_3^2.\label{201812301849} \end{align}
\end{lem}

\begin{lem}
\label{lem:dfifessim2057sadfafdas}
Under the assumptions of  \eqref{2019022201518} and \eqref{aprpiosesnew}, we have
\begin{equation} \mathfrak{E}\mbox{ is equivalent to }\|\eta\|_5^2+\|(u,\vartheta)\|_4^2 . \label{201702071610nb}
\end{equation}
\end{lem}
%%%%%%%%%%%%%%%%%%%%%%%%%%%%%%

With the five lemmas in hand, next we derive  \emph{a prior} Gronwall-type energy inequality for the transformed MB problem.
\begin{lem}  \label{pro:0301n0845}
There exist  constants ${\delta}_1\in (0,1)$ and ${C}_1>0$ such that, for any $\delta\leqslant {\delta}_1$,  then, under the assumption of \eqref{aprpiosesnew}, $(\eta,u,\vartheta)$ satisfies the  Gronwall-type energy inequality,
for any $t\in I_T$,
\begin{align}
\tilde{\mathfrak{E}}(t)   \leqslant \Lambda^*\int_0^t
\tilde{\mathfrak{E}}(\tau) \mm{d}\tau +  {C}_1\left(\|\eta^0\|_5^2+ \|(u^0,\vartheta^0)\|_4^2 +\int_0^t
\|( u,\vartheta)\|_{0}^2\mm{d}\tau\right)
\label{2016121521430850}
\end{align}and the equivalent estimate
\begin{align}
\label{2018008121027}
\mathfrak{E}(t)\leqslant {C}_1\tilde{\mathfrak{E}}(t)\lesssim \mathfrak{E}(t),
\end{align}
see Lemma \ref{thm:0201201622xx} for the definition of $\Lambda^*$.
\end{lem}
\begin{pf}
We can derive from Lemmas \ref{201612132242nn}--\ref{201612132242nxsfssdfs} and \ref{lem:dfifessim2057}, and \eqref{201903192036} that, for any sufficiently large constant  $c_3^l>1$ and for any sufficiently small constant $c_3^s\in (0,1)$,
 \begin{align}
 \label{201702061553}
\frac{\mm{d}}{\mm{d}t} \tilde{\mathfrak{E}}+c \tilde{\mathfrak{O}}\leqslant & c_1^I  (
(c_3^l)^2  ( \|(u,\vartheta)\|_{4}^2 +
 \sqrt{{\mathfrak{E}}} {\mathfrak{D}}  )+ c_3^s(\|  \eta\|_{\underline{5},0}^2 +\|u\|_3^2\|\nabla q\|_1^2)) ,
 \end{align}
  where
 $$\begin{aligned}
 \tilde{\mathfrak{E}}:=& \sum_{|\alpha|\leqslant 4}
   \partial_\mm{h}^{\alpha} \eta \cdot  \partial_\mm{h}^{\alpha} u    +
\frac{1}{2}\|\nabla \eta\|_{\underline{4},0}^2
+(c_3^l)^2(Q\|\partial_3  \eta\|_{\underline{4},0}^2  +\| (u,\vartheta)\|^2_{\underline{4},0})+c_3^s\|\eta\|_{5,*}^2\nonumber
  \\
 & +c_3^l(\|  u_{t}\|_{0}^2+E_1(u,\vartheta_t)) + \|  u_{tt}\|_{0}^2
+2\int \nabla q_t\cdot  D^{t,2}_{u}\mm{d}y  +E_1(u_t,\vartheta_{tt}) ,   \\
 \tilde{\mathfrak{O}}:= &\|\partial_3\eta\|_{\underline{4},0}^2+  (c_3)^2   \| (u ,\vartheta)\|_{\underline{4},1}^2 +  c_3^l \|(u,\vartheta)_t\|_{1}^2+\|(u,\vartheta)_{tt}\|_{1}^2+c_3^s( \left\|  \left( \eta,u \right)\right\|_5^2+\| \nabla q \|_3^2).
\end{aligned}$$

Similarly to \eqref{201902282100}, using Lemma \ref{201612132242nx},
  there exists constants $c$, $c_3^l$, $c^s_3:=\Lambda^*/2 c_1^I $ and $\delta_1$  such that, for any $\delta\leqslant \delta_1$,
\begin{align}
&\tilde{\mathfrak{E}},\ \mathfrak{E}\mbox{ and } \|\eta\|_5^2+\|(u,\vartheta)\|_4^2 \mbox{ are equivalent},\label{201808safdas02asdfsadf1219xxafasf}\\
\label{201808safdas02asdfsadf1219xx}
& \|\nabla \eta\|_{\underline{4},0}^2 /2 \leqslant  \tilde{\mathfrak{E}},\\
\label{201808061629xx}&
\tilde{\mathfrak{O}}\mbox{ is equivalent to }\mathfrak{D},
\end{align}
where the equivalence coefficients in \eqref{201808safdas02asdfsadf1219xxafasf} and \eqref{201808061629xx} can be independent of $\delta$.
 Using interpolation inequality, we immediately deduce   the desired conclusions  \eqref{2016121521430850} and \eqref{2018008121027} from \eqref{201702061553}--\eqref{201808061629xx}.
 \hfill $\Box$\end{pf}

Similarly to Proposition \ref{pro:0401n}, we also prove the existence of a unique local-in-time classical solution for the transformed MB problem under the initial-value condition
\begin{equation}
\label{201903081829}
(\eta^0,u^0,\vartheta^0)\in H^{5,1}_0\times H_\sigma^4\times H_0^4 .
\end{equation}
Consequently we have the following conclusions.
\begin{pro} \label{pro:0401nxd}
  \begin{enumerate}
    \item[(1)]  Let $(\eta^0,u^0,\vartheta^0)$ satisfy \eqref{201903081829} and  $\zeta^0:=\eta^0+y$.
           There exists a sufficiently small ${\delta}_2\in (0,1)$, such that, if $(\eta^0,u^0,\vartheta^0)$ satisfying
\begin{align}
\label{201801032021}
& \sqrt{\|\eta^0\|_5^2+\|(u^0,\vartheta^0)\|_4^2}<  {\delta}_2
\end{align}
and necessary compatibility conditions
\begin{alignat}{2}
& \mm{div}_{\mathcal{A}^0}u^0=0& &\ \mbox{ in }\Omega, \label{201808040977}\\
& u_t|_{t=0}=0\mbox{ and }\vartheta_t|_{t=0}=0 & &\ \mbox{ in }\Omega ,\label{201808040977xx}
\end{alignat}
then there is a  local existence time $T^{\max}>0$ (depending on ${\delta}_2$, the domain  and the known parameters), and a unique local-in-time classical  solution
$(\eta, (u,\vartheta),q)\in C^0([0,T^{\max}) ,  H^{5,1}_{0}\times H_0^4 \times \underline{H}^3) $ to the transformed MB problem, where the solution enjoys the regularity \begin{align}
& ((u,\vartheta)_t, (u,\vartheta)_{tt}, q_t  ) \in C^0(\overline{I_{ T^{\max}}}, H^2_0 \times L^2\times \underline{H}^1), \label{201811142106}\\
& ((u, \vartheta), (u,\vartheta)_t, (u,\vartheta)_{tt}, q, q_t ) \in L^2(I_{ T^{\max}}, H^5_0\times  H^3_0\times H^1_0\times \underline{H}^4\times \underline{H}^2). \label{201811142107}
\end{align}
    \item[(2)]
 In addition, if the solution $(\eta,u,\vartheta)$  further satisfies
$$\sup_{t\in [0,T)}\sqrt{\|\eta(t)\|_5^2 +\|(u,\vartheta)(t)\|_4^2}\leqslant {\delta}_1\mbox{ for some }T<T^{\max},$$
 then  $(\eta, u,\vartheta)$ enjoys the Gronwall-type energy inequality \eqref{2016121521430850} and
 the equivalent estimate \eqref{2018008121027}.
  \end{enumerate}
\end{pro}
\begin{rem}
For any given initial data $(\eta^0, ({u}^0, \vartheta^0))\in H^{5,1}_0\times H^4_0 $ satisfying  \eqref{201801032021}--\eqref{201808040977} with sufficiently small $\delta_2$, there exists  a unique local-in-time strong solution $(\eta, ({u},\vartheta),q)\in C^0([0,T),H^{3,1}_{0}$ $\times H^2_0\times \underline{H}^1)$. Moreover, the initial date of  $q$ is a weak solution to
\begin{equation}
\left\{
\begin{array}{ll}
\Delta_{\mathcal{A}^0} q^0  =\mm{div}_{\mathcal{A}^0} \left(( \Delta_{\mathcal{A}^0}u^0+Q\partial_3^2\eta^0+R \vartheta^0e_3)   - u^0 \cdot \nabla_{\mathcal{A}^0} u^0  \right)  &\mbox{ in }\Omega,\\
 \nabla_{\mathcal{A}^0}q^0  \cdot   e_3 =\left( \Delta_{\mathcal{A}^0}u^0+ Q \partial_3^2\eta^0 \right) \cdot   e_3 &\mbox{ on }\partial\Omega.
\end{array}
\right. \label{201808040916}
\end{equation}
If the condition \eqref{201808040977xx} is further satisfied, i.e., $(\eta^0,u^0,\vartheta^0,q^0)$ satisfies
\begin{equation}
\label{201812271841x}
\nabla_{\mathcal{A}^0} q^0- \Delta_{\mathcal{A}^0} u^0 -Q\partial_3^2\eta^0=0\mbox{ and }
 \Delta_{\mathcal{A}^0}\vartheta^0=0   \mbox{ on }\partial\Omega,
\end{equation}
then we can improve the regularity of $(\eta,u,\vartheta,q)$ so that it is a classical solution for sufficiently small ${\delta}_2$.

In addition, since there exists a constant ${\delta}_3$ such that, for any $\|w\|_3\leqslant {\delta}_3$
$$\|\nabla f\|_0 \leqslant c\|\nabla_{\mathcal{A}} f\|_0 \mbox{ for any }f\in H^1,$$
where $\mathcal{A}^{\mm{T}}:=(\nabla w)^{-1}$. Thus we easily see that, for any
$\|\eta^0\|_3\leqslant {\delta}_3$, the solution  $q^0\in \underline{H}^1$ to \eqref{201808040916} is unique  for given $(\eta^0,u^0,\vartheta^0)$.
\end{rem}

\begin{rem}
\label{201811031411}
For any classical solution $(\eta,u,\vartheta,q)$ constructed by Proposition \ref{pro:0401nxd}, we take  $(\eta,u,\vartheta$, $q)|_{t=t_0}$ as a new initial  datum,  where $t_0\in I_{T^{\max}}$. Then the new initial data can define a unique local-in-time classical solution $(\tilde{\eta},\tilde{u},\tilde{\vartheta},\tilde{q})$ constructed by Proposition
\ref{pro:0401nxd}; moreover the initial datum of $\tilde{q}$ is equal to $q|_{t=t_0}$, if $\|\eta\|_3 \leqslant \delta_3 $ in $I_{T^{\max}}$.
\end{rem}

%%%%%%%%%%%%%%%x%%%%%%%%%%%%%%%%%%%%%%%%%%%%%%%%%%%%%%%%%%%%%%%%%

\subsection{Construction of initial data for  nonlinear solutions}

For any given $\delta>0$,  let
\begin{equation}\label{0501}
\left(\eta^\mm{a}, u^\mm{a}, \vartheta^\mm{a} ,q^{\mm{a}}\right)=\delta e^{{\Lambda^* t}} (\tilde{\eta}^0,\tilde{u}^0, \vartheta^0 ,\tilde{q}^0),
\end{equation}
where $(\tilde{\eta}^0,\tilde{u}^0, \tilde{\vartheta}^0, \tilde{q}^0):=(w/\Lambda^*,w,\theta,\beta)$, and $(w,\theta,\beta)\in  H^5_\sigma \times H^5_0 \times  \underline{H}^4$ comes from  Proposition \ref{thm:0201201622} with $\tau=1$.
Then $\left(\eta^\mm{a}, u^\mm{a}, \vartheta^\mm{a}, q^\mm{a}\right)$  is a solution to the linearized MB problem, and enjoys the estimate, for any $i \geqslant 0$,
 \begin{equation}
\label{appesimtsofu1857}
\|\partial_t^i(\eta^{\mm{a}},u^{\mm{a}},\vartheta^{\mm{a}} )\|_{5}+\|\partial_t^i q^{\mm{a}}\|_{4}\leqslant c(i) \delta e^{\Lambda^* t}.
 \end{equation}
 Moreover, by \eqref{201602081445MH},
 $$
m_0 :=
    \min  \{ \|\tilde{\eta}^0_{\mm{h}}\|_{L^1},\|\tilde{\eta}^0_3\|_{L^1},\|\tilde{u}^0_{\mm{h}}\|_{L^1},\|\tilde{u}^0_3\|_{L^1}, \|\partial_3 \tilde{\eta}^0_{\mm{h}}\|_{L^1}, \|\partial_3 \tilde{\eta}^0_3\|_{L^1} , \|\vartheta\|_{L^1}\}>0.$$

Next we shall modify the initial data of the linear solutions.
\begin{pro}
\label{lem:modfied}
Let $(\tilde{\eta}^0,\tilde{u}^0,\tilde{\vartheta}^0,\tilde{q}^0)$
be the same as in \eqref{0501}, then there exists a constant ${\delta}_4$, such that for any $\delta\in(0,{\delta}_4)$, there exists $(\eta^{\mm{r}},(u^{\mm{r}},\vartheta^{\mm{r}}),q^{\mm{r}})\in H^5_0\times H^4_0\times \underline{H}^3$ enjoying the following properties:
\begin{enumerate}[\quad (1)]
                \item The modified initial data
  \begin{equation}\label{mmmode04091215}(  {\eta}_0^\delta,{u}_0^\delta,\vartheta_0^\delta,q_0^\delta): =\delta
   (\tilde{\eta}^0,\tilde{u}^0,\tilde{\vartheta}^0, \tilde{q}^0) + \delta^2 (\eta^\mm{r},   u^\mm{r}, \vartheta^\mm{r},q^{\mm{r}})
\end{equation}
belongs to $H^{5,1}_{0}\times H^4_0\times H^4_0\times \underline{H}^3$, and satisfies
\begin{align}
\det \nabla ( {\eta}_0^\delta+y)=1,\nonumber
\end{align}
the compatibility conditions  \eqref{201808040977}, \eqref{201808040916}$_1$ and \eqref{201812271841x} with $(\eta_0^\delta,u_0^\delta,\vartheta_0^\delta,q_0^\delta)$ in place of $(\eta^0,u^0$, $\vartheta^0, q^0)$.
\item Uniform estimate:
\begin{equation}
\label{201702091755}
 \sqrt{\|\eta^{\mm{r}}\|_5^2+ \|(u^{\mm{r}},\vartheta^{\mm{r}})\|_4^2+\|q^{\mm{r}}\|_{3}^2 }\leqslant {C}_2,
 \end{equation}
where the constant ${C}_2\geqslant 1$ depends on the domain,   and the known parameters, but is independent of $\delta$.
 \end{enumerate}
\end{pro}
\begin{pf}
Recalling the construction of $(\tilde{\eta}^0,\tilde{u}^0,\tilde{\vartheta}^0,\tilde{q}^0)$, we can see that
$(\tilde{\eta}^0,\tilde{u}^0,\tilde{\vartheta}^0,\tilde{q}^0)$ satisfies
\begin{equation}
\label{201705141510}
\left\{
\begin{array}{ll}
 \mm{div}\tilde{\eta}^0=\mm{div}\tilde{u}^0=0   &\mbox{ in }\Omega,\\
 \Lambda  \tilde{u}^0-\Delta \tilde{u}^0 +\nabla\tilde{q}^0=Q\partial_3^2\tilde{\eta}^0 +R \vartheta e_3  & \mbox{ in }\Omega,\\
    \Lambda P_{\vartheta} \tilde{\theta}^0  - \Delta  \tilde{\theta}^0=R \tilde{u}_3^0   & \mbox{ in }\Omega,\\
 (\tilde{\eta}^0,\tilde{u}^0,\tilde{\vartheta}^0) = {0} & \mbox{ on }\partial\Omega.
 \end{array}\right.
\end{equation}
If $(\eta^\mm{r},(u^\mm{r},\vartheta^\mm{r}),q^\mm{r})\in H^5_0\times H^4_0\times \underline{H}^3$ satisfies
\begin{equation}
\label{201702061320}
\left\{\begin{array}{ll}
\mm{div}\eta^{\mm{r}}=O(\eta^{\mm{r}}),\ \mm{div}u^\mm{r} = \mm{div}((\tilde{\mathcal{A}}_0^\delta)^{\mm{T}} u_0^\delta)/\delta^2 &\mbox{ in }\Omega,\\
 \nabla q^{\mm{r}}- \Delta u^{\mm{r}}= Q\partial_3^2\eta^{\mm{r}}+  R\vartheta^{\mm{r}}_3e_3 + \Upsilon  & \\
+ \left( (\nabla_{\tilde{\mathcal{A}}_0^\delta}\mm{div}_{\mathcal{A}_0^\delta}u^\delta_0+
 \nabla\mm{div}_{\tilde{\mathcal{A}}_0^\delta}u^\delta_0) -   u^\delta_0 \cdot \nabla_{\mathcal{A}^0_\delta} u^\delta_0
  -\nabla_{\mathcal{A}_0^\delta}q^\delta_0
  \right)/\delta^2  &\mbox{ in }\Omega,  \\
 - \Delta \vartheta^{\mm{r}}=    \left( (\nabla_{\tilde{\mathcal{A}}_0^\delta}\mm{div}_{\mathcal{A}_0^\delta}\vartheta^\delta_0+
 \nabla\mm{div}_{\tilde{\mathcal{A}}_0^\delta}\vartheta^\delta_0)
  \right)/\delta^2   &\mbox{ in }\Omega,  \\
     \mm{div}_{\mathcal{A}_0^\delta}\Upsilon =- \delta \Lambda \mm{div}_{\tilde{\mathcal{A}}_0^\delta}\tilde{u}^0  &   \mbox{ in }\Omega,\\
(\eta^\mm{r}, u^\mm{r},\vartheta^\mm{r},\Upsilon  )=0&\mbox{ on }\partial\Omega,
\end{array}\right.\end{equation}
 where $({\eta} ^\delta_0,{u} ^\delta_0, \vartheta^\delta_0,q^\delta_0)$ is given in the mode \eqref{mmmode04091215},
  $\zeta^\delta_0:=\eta^\delta_0 +y$, $\mathcal{A} ^\delta_0:= (\nabla\zeta^{\delta}_0)^{-\mm{T}}$, $\tilde{\mathcal{A}}^{\delta}_0 := \mathcal{A} ^\delta_0-I$ and $O(\eta^{\mm{r}}):=\delta^{-2}\mm{div}
\Phi(\delta\tilde{\eta}+\delta^2\eta^{\mm{r}})$,
see \eqref{201903082040} for the definition of $\Phi$, then, by \eqref{201903082042} with
$\eta^\delta_0$ in place of $\eta$ and \eqref{201705141510}, it is easy to check that $({\eta} ^\delta,({u} ^\delta,\vartheta^\delta),q^\delta)$  belongs to $H^{5}_{0}\times H^4_0  \times \underline{H}^3$, and satisfies
\eqref{201808040977}, \eqref{201808040916}$_1$ and \eqref{201812271841x}  with $(\eta ^\delta,u ^\delta,\vartheta^\delta ,q^\delta)$ in place of $(\eta^0 ,u^0 ,\vartheta^0 ,q^0 )$.

However, making use of the existence theory of Stoke problem, iterative technique  and compactness convergence method as in \cite[Lemma 4.2]{JFJSOMARMA}, we easily establish the following conclusions:
for sufficiently small $\delta_4$, such that, for any $\delta \leqslant \delta_4$,
\begin{itemize}
  \item there exists  $(\eta^{\mm{r}}, \varpi)\in H^5_0\times \underline{H}^4$ such that
\begin{equation}
\label{20170514xx}
 \left\{\begin{array}{l}
\nabla \varpi- \Delta \eta^{\mm{r}} =0,\quad  \mm{div}\eta^{\mm{r}} =O(\eta^{\mm{r}}) ,\\
\eta^{\mm{r}}|_{\partial\Omega} =0
\end{array}\right.  \end{equation}
and $\|\eta^{\mm{r}}\|_5\leqslant c$.
  \item
 there exists $(\Upsilon ,q)\in H_0^2\times \underline{H}^1$ such that
\begin{equation}
\label{201702051905xfwea}
 \left\{\begin{array}{l}
\nabla q-  \Delta \Upsilon   =0, \
 \mm{div}\Upsilon =\mm{div}_{\tilde{\mathcal{A}}_0^\delta}(\Upsilon - \delta \Lambda \tilde{u}^0) , \\
\Upsilon |_{\partial\Omega}=0
\end{array}\right.\end{equation}
and
\begin{equation}
\|\Upsilon\|_2+\|q\|_{1}\lesssim  \delta^2,
\label{xx}
\end{equation}
where $\eta^{\mm{r}}$ is constructed by \eqref{20170514xx}.
  \item
 there exists  $(u^{\mm{r}}, q^{\mm{r}})\in H^4_0 \times \underline{H}^3$ such that
\begin{equation}
\label{201702061320xxx}
\left\{\begin{array}{ll}
 \nabla q^{\mm{r}}- \Delta u^{\mm{r}}= Q\partial_3^2\eta^{\mm{r}}+R\vartheta^{\mm{r}}_3e_3+ \Upsilon  & \\
 + \left( (\nabla_{\tilde{\mathcal{A}}_0^\delta}\mm{div}_{\mathcal{A}_0^\delta}u^\delta_0+
 \nabla\mm{div}_{\tilde{\mathcal{A}}_0^\delta}u^\delta_0) -   u^\delta_0 \cdot \nabla_{\mathcal{A}^0_\delta} u^\delta_0
  -\nabla_{\mathcal{A}_0^\delta}q^\delta_0
  \right)/\delta^2   ,\\
\mm{div}u^\mm{r} =\mm{div}((\tilde{\mathcal{A}}_0^\delta)^{\mm{T}} u_0^\delta)/\delta^2  ,\\
    u^\mm{r}|_{\partial\Omega } =0&
\end{array}\right.\end{equation}
and $\|u^{\mm{r}}\|_4\leqslant c$,
where  $\eta^{\mm{r}}$ and $\Upsilon $ are constructed in \eqref{20170514xx} and \eqref{201702051905xfwea}.
\item
   there exists  $ \vartheta^{\mm{r}} \in H^4_0 $ such that
\begin{equation}
\label{201702061320xx}
\left\{\begin{array}{ll}
- \Delta \vartheta^{\mm{r}}=    \left( (\nabla_{\tilde{\mathcal{A}}_0^\delta}\mm{div}_{\mathcal{A}_0^\delta}\vartheta^\delta_0+
 \nabla\mm{div}_{\tilde{\mathcal{A}}_0^\delta}\vartheta^\delta_0)
  \right)/\delta^2 , \\
    u^\mm{r}|_{\partial\Omega } =0&
\end{array}\right.\end{equation}
and $\|\vartheta^{\mm{r}}\|_4\leqslant c$, where  $\eta^{\mm{r}}$ is constructed in \eqref{20170514xx}.
\end{itemize}
It should be noted that all constants $c$ above are independent of $\delta$.

Consequently, we immediately see that the construction of $(\eta^\mm{r},u^\mm{r},\vartheta^\mm{r},q^\mm{r} )$ constructed above satisfies \eqref{201702091755} and \eqref{201702061320}. This completes the proof of Proposition \ref{lem:modfied}.
\hfill $\Box$
\end{pf}

\subsection{Error estimates and existence of escape times}

Let
\begin{align}
&{C}_3:= \sqrt{\|\tilde{\eta}^0\|_5^2+ \|(\tilde{u}^0,\tilde{\vartheta}^0 )\|_4^2 } +{C}_2 \geqslant 1,\label{201903101530}\\
 & \delta<{\delta}_0:= \min\{{\delta}_1,{\delta}_2,\delta_3, {C}_3{\delta}_4\}/2{C}_3<1,\nonumber
 \end{align}
and $(\eta_0^\delta,u_0^\delta,\vartheta_0^\delta,q_0^\delta)$ be constructed by Proposition \ref{lem:modfied}.

 Noting that
\begin{align}
  \sqrt{\|\eta_0^\delta\|_5^2
+\|(u_0^\delta, \vartheta_0^\delta)\|_4^2 }\leqslant {C}_3\delta<2{C}_3{\delta}_0 \leqslant {\delta}_2 ,
\label{201809121553xx}
\end{align}
by the first assertion in Proposition \ref{pro:0401nxd}, there exists a (nonlinear) solution $(\eta,u ,\vartheta,q)$ of the transformed MB problem defined on some time interval $I_{T^{\max}}$ with initial value $(  {\eta}_0^\delta ,{u}_0^\delta,\vartheta_0^\delta )$.

Let $\epsilon_0\in (0,1)$ be a constant, which will be defined in \eqref{defined}.
 We define
 \begin{align}\label{timesxx}
& T^\delta:=(\Lambda^*)^{-1}\mm{ln}({\epsilon_0}/{\delta})>0,\quad\mbox{i.e.,}\;
 \delta e^{\Lambda^* T^\delta }=\epsilon_0,\\
&T^*:=\sup\left\{t\in I_{T^{\max}} \left|~\sqrt{{\|\eta (\tau)\|_5^2+
\| (u , \vartheta)(\tau)\|_4^2 }}\right.\right. \leqslant 2{C}_3{\delta}_0\mbox{ for any }\tau\in [0,t)\bigg\},\nonumber\\
&  T^{**}:=\sup\left\{t\in  I_{T^{\max}} \left|~\left\|(\eta,u,\vartheta)(\tau)\right\|_{0}\leqslant 2 {C}_3\delta e^{\Lambda \tau}\mbox{ for any }\tau\in [0,t)\right\}.\right.\nonumber
 \end{align}
 thus $T^*>0$ by \eqref{201809121553xx} and the first assertion in Proposition \ref{pro:0401nxd}.  Similarly, we also have  $T^{**}>0$.
Moreover, by the first assertion in Proposition \ref{pro:0401nxd} and Remark \ref{201811031411},
we can easily see that
 \begin{align}\label{0502n1xx}
&\sqrt{\|\eta (T^*)\|_5^2+
\| (u, \vartheta) (T^*)\|_4^2 }=2{C}_3 {\delta}_0,\mbox{ if }T^*<\infty , \\
& \left\|(\eta,u,\vartheta) (T^{**}) \right\|_0
=2 {C}_3\delta e^{\Lambda T^{**}},\mbox{ if }T^{**}<T^{\max}.\label{0502n111}
\end{align}

We denote ${T}^{\min}:= \min\{T^\delta ,T^*,T^{**} \}$. Noting that $(\eta,u,\vartheta)$ satisfies
$$
 \sup_{0\leqslant t< T^{\min}}\sqrt{\|\eta(t)\|_5^2 +\|(u,\vartheta)(t)\|_4^2 }\leqslant {\delta}_1,$$
thus,  for any $t\in I_{T^{\min}}$, $(\eta,u,\vartheta)$ enjoys  \eqref{2016121521430850} and \eqref{2018008121027} with $(\eta^\delta_0,u^\delta_0,\vartheta^\delta_0)$ in place of $(\eta^0,u^0,\vartheta^0)$ by the second conclusion in Proposition \ref{pro:0401nxd}.
Noting that $\|(\eta,u,\vartheta)(t)\|_0\leqslant 2{C}_3 \delta e^{\Lambda^* t}$ in $I_{T^{\min}}$, then we deduce from the estimates \eqref{2016121521430850} and \eqref{201809121553xx} that, for all $t\in I_{T^{\min}}$,
\begin{align}
  \tilde{\mathfrak{E}}(t)  \leqslant c \delta^2e^{2\Lambda^* t}+\Lambda^*\int_0^t \tilde{\mathfrak{E}}(\tau)  \mm{d}\tau.
\label{0503xx} \end{align}
Applying the Gronwall's lemma to the above estimate  yields
$\tilde{\mathfrak{E}}(t)\lesssim (\delta e^{\Lambda^* t})^2$, which, together with \eqref{2018008121027}, further implies
\begin{align}
& \mathfrak{E}(t) \leqslant ({C}_4\delta e^{\Lambda^* t})^2\leqslant {C}_4^2\epsilon_0^2\mbox{ in } I_{T^{\min}}.
\label{201702092114xx}
\end{align}

Next we estimate the error between the (nonlinear) solution $(\eta,u,\vartheta )$ and the linear solution $(\eta^\mm{a},u^{\mm{a}},\vartheta^{\mm{a}})$ provided by \eqref{0501}.
\begin{pro}\label{lem:0401xxxxxx}
 If $Q$ satisfies \eqref{201903071510},  then there exists a constant ${C}_5$ such that,  for any $t\in I_{T^{\min}}$,
\begin{align}
\label{ereroe}
&Q\|\partial_3\eta^{\mm{d}}\|_{\mathfrak{X}}+\|(\eta^{\mm{d}}, u^{\mm{d}}, \vartheta^{\mm{d}})\|_{\mathfrak{X}}  \leqslant C_5\sqrt{\delta^3e^{3\Lambda^* t}},
\end{align}
where $(\eta^{\mathrm{d}}, u^{\mathrm{d}}, \vartheta^{\mathrm{d}}):=(\eta, u,\vartheta)-(\eta^\mm{a},u^{\mm{a}},\vartheta^{\mm{a}})$, $ \mathfrak{X}=L^{1}$ or $L^2$, and ${C}_5$ is independent  of  $\delta$,  and $T^{\min}$.
\end{pro}
\begin{pf}
Subtracting the both transformed MB problem  and the linearized MB problem (i.e., \eqref{201803101401} with $\tau=1$), we get
\begin{equation}\label{201702052209}\left\{\begin{array}{l}
\eta_t^{\mathrm{d}}=u^{\mathrm{d}}  ,\\[1mm]
  u_{t}^{\mathrm{d}}- \Delta u^{\mm{d}}+\nabla q^{\mm{d}} -Q \partial_3^2 \eta^{\mm{d}}- R \vartheta^{\mm{d}}e_3
= N^{\mm{h}}_u+N^{\mm{h}}_q, \\
  P_\vartheta \vartheta_{t}^{\mathrm{d}}- \Delta \vartheta^{\mm{d}}  - R u^{\mm{d}}_3=
 N^{\mm{h}}_\vartheta, \\
\div u^{\mm{d}}=D_{\mm{u}}^h ,\\
 (\eta^{\mathrm{d}},u^{\mathrm{d}},\vartheta^{\mathrm{d}})|_{t=0}=\delta^2 (\eta^{\mm{r}},u^{\mm{r}},\vartheta^{\mm{r}}) ,\\
 (\eta^{\mathrm{d}},u^{\mathrm{d}},\vartheta^{\mathrm{d}})|_{\partial\Omega}=0  .\end{array}\right.\end{equation}
Similarly to \eqref{estimforhoe1nwe}, multiplying \eqref{201702052209}$_2$ and \eqref{201702052209}$_3 $ by $u^{\mathrm{d}}$ and $\vartheta^{\mathrm{d}}$ in $L^2$, resp., and then adding the two resulting identities together, we have
\begin{align}
&\frac{1}{2}\frac{\mm{d}}{\mm{d}t}\left(Q \|\partial_3 \eta^{\mathrm{d}}\|_{0}^2 +\|  u^{\mathrm{d}}\|^2_0+{P_\vartheta} \|\vartheta^{\mathrm{d}}\|_0^2
\right)-D_R(u^{\mathrm{d}},\vartheta^{\mathrm{d}},1) =R_1(t),\nonumber
\end{align}
where
\begin{align}
\label{201903111523}
R_1(t):=  \int \left( ( N^{\mm{h}}_u+ N^{\mm{h}}_q)\cdot    u^{\mathrm{d}} +   q^{\mathrm{d}}\mm{div} u^{\mathrm{d}}+ N^{\mm{h}}_\vartheta   \vartheta^{\mathrm{d}} \right)\mm{d}y.
\end{align}
Integrating the above identity in time from $0$ to $t$ yields \eqref{0314}.

By the existence theory of Stokes problem, there exists $(\tilde{u},\varpi)\in H_0^2\times  \underline{H}^1$ such that, for given $(\eta,u)$,
\begin{equation*}
 \left\{\begin{array}{ll}
\nabla \varpi-\Delta \tilde{u}=0,\quad
 \mm{div}\tilde{u} =-  \div_{\tilde{\mathcal{A}}}
  u  &\mbox{ in }\Omega,  \\
\tilde{u}   =0&\mbox{ on } \partial\Omega.
\end{array}\right.\end{equation*}
 Moreover, by \eqref{201702092114xx},
\begin{equation}
\label{201705012219xx}
\|\tilde{u}\|_2\lesssim \|\mm{div}_{\tilde{\mathcal{A}}}u\|_1 \lesssim \delta^2 e^{2\Lambda^* t} \lesssim \sqrt{\delta^3 e^{3\Lambda^* t}}\mbox{ in } I_{T^{\min}}.
\end{equation}
It is easy to see that $v^{\mm{d}}:=u^{\mm{d}}-
\tilde{u}\in H_\sigma^2$.
We derive from \eqref{Lambdard} that
  \begin{equation*}\begin{aligned}
  D_{R}(v^{\mm{d}},\vartheta^{\mm{d}},1-\sqrt{Q})  \leqslant  \Lambda(1-\sqrt{Q})(\| v^{\mathrm{d}}\|^2_0+P_\vartheta\|  \vartheta^{\mathrm{d}}\|_0^2 ) +Q\|\partial_3 v^{\mm{d}}\|_0^2/ \Lambda(1-\sqrt{Q}).
\end{aligned}\end{equation*}
Exploiting  \eqref{0501}, \eqref{201702092114xx} and \eqref{201705012219xx}, immediately implies
\begin{align}
& D_{R}(u^{\mm{d}},\vartheta^{\mm{d}},1-\sqrt{Q}) \nonumber   \\
&\leqslant  \Lambda(1-\sqrt{Q})(\| u^{\mathrm{d}}\|^2_0+P_\vartheta\|  \vartheta^{\mathrm{d}}\|_0^2 ) +Q\|\partial_3 u^{\mm{d}}\|_0^2/ \Lambda(1-\sqrt{Q})+c\delta^3 e^{3\Lambda^* t}.\nonumber
\end{align}

Similarly to \eqref{201902171836},  we  easily estimate that
\begin{align}
\int_0^t {R}_1(\tau)\mm{d}\tau  \lesssim   \delta^3 e^{3\Lambda^* t}.\nonumber
\end{align}
By \eqref{201702052209}$_5$, we have
\begin{equation}
R_2\lesssim \delta^4\|(\eta^{\mm{r}},u^{\mm{r}},\vartheta^{\mm{r}})\|_2^2\lesssim  \delta^3 e^{3\Lambda^* t}. \nonumber
\label{2018110702038}
\end{equation}
Noting that $Q$ satisfies \eqref{201903222102xx}, then,  by \eqref{201903101106}, $$\sqrt{Q}\leqslant  \Lambda(1-\sqrt{Q}).$$
Thanks to the four estimates above, we can derive \eqref{201903101526} from \eqref{0314}.

Exploiting  \eqref{201903101106}, we have
$$  \Lambda(1-\sqrt{Q})-\Lambda^*\leqslant   \sqrt{Q}\left( \mathcal{R} +\sqrt{ \mathcal{R}}\right).$$
Noting that \eqref{201903222102}, we can use \eqref{201903101106} and  the above relation to deduce that
$$\begin{aligned}
\Lambda^* \geqslant  \Lambda_0-\sqrt{ Q\mathcal{R}}
\geqslant    \sqrt{Q}\left( \mathcal{R} +\sqrt{ \mathcal{R} }\right) \geqslant 2(\Lambda(1-\sqrt{Q})-\Lambda^*),
\end{aligned}$$
which yields \eqref{201903101525}.
 Thus, exploiting Gronwall's lemma and  \eqref{201903101525}, we derive from \eqref{201903101526} that
\begin{align}
Q \|\partial_3 \eta^{\mathrm{d}}\|_{0}^2+\| ( u^{\mathrm{d}},\vartheta^{\mathrm{d}}) \|_0^2  \leqslant c\delta^3 e^{3\Lambda^* t}.\label{0314xxx}
\end{align}

We turn to derive the error estimate for $\eta^{\mathrm{d}}$.
It follows from \eqref{201702052209}$_1$ that
\begin{equation*}\begin{aligned}
 \frac{\mm{d}}{\mm{d}t}\|\eta^{\mathrm{d}}\|_0^2
\lesssim   \| u^{\mathrm{d}}\|_0 \|\eta^{\mathrm{d}}\|_0.
\end{aligned}\end{equation*}
Therefore, using \eqref{0314xxx} and  the initial condition ``$\eta^{\mm{d}}|_{t=0}=\delta^2 \eta^{\mm{r}}$",  it follows that
\begin{equation}\begin{aligned}\label{erroresimts}
 \|\eta^{\mathrm{d}}\|_0\lesssim  &  \int_0^t \|u^{\mathrm{d}}\|_0 \mm{d}\tau
+\delta^2\|\eta^{\mm{r}}\|_0\lesssim  \sqrt{ \delta^3e^{3\Lambda^* t}}.
\end{aligned}\end{equation}
Noting that $L^2\hookrightarrow L^1$, then we can derive
 \eqref{ereroe}  from \eqref{0314xxx} and  \eqref{erroresimts}.
This completes the proof of Proposition \ref{lem:0401xxxxxx}. \hfill$\Box$
\end{pf}

Now we define that
 \begin{equation}\label{defined}
\epsilon_0:=\min\left\{
 \frac{C_3\delta_0}{  C_4},
\frac{C_3^2}{4 C_5^2},
\frac{m_0^2}{4C_5^2}   \right\}>0,
 \end{equation}Consequently, we further have the  relation
\begin{equation}
\label{201702092227xx}
T^\delta =T^{\min}\neq T^*\mbox{ or }T^{**},
\end{equation}
which can be showed by contradiction as follows:

If $T^{\min}=T^*$, then   $T^*<\infty$. Noting that $\epsilon_0\leqslant C_3\delta_0/C_4$, thus we deduce from \eqref{201702092114xx} that
\begin{equation*}\begin{aligned}
\sqrt{\|\eta\|_{5}^2+ \|(u,\vartheta) (T^*)\|_4^2}\leqslant C_3\delta_0< 2 C_3 \delta_0,
 \end{aligned} \end{equation*}
which contradicts \eqref{0502n1xx}. Hence, $T^{\min}\neq T^*$.

If $T^{\min}=T^{**}$, then   $T^{**}<T^*\leqslant T^{\max}$.
 Noting that $\sqrt{\epsilon_0}\leqslant C_3/2C_5$, we can deduce from  \eqref{0501},  \eqref{201903101530},
\eqref{timesxx} and \eqref{ereroe}   that \begin{equation*}\begin{aligned}
\| (\eta,u,\vartheta)  (T^{**})\|_{0}
&\leqslant \|( \eta^\mm{a}, u^\mm{a},\vartheta^\mm{a} )(T^{**})\|_{0}+\|( \eta^\mm{d}, u^\mm{d}, \vartheta^\mm{d} ) (T^{**})\|_{0}\\
&\leqslant  \delta e^{{\Lambda^* T^{**}}}(C_3+ C_5\sqrt{\delta e^{\Lambda T^{**}}})\\
&\leqslant  \delta e^{{\Lambda^* T^{**}}}(C_3+ C_5\sqrt{\epsilon_0})
\leqslant 3C_3 \delta e^{\Lambda^* T^{**}}/2<2C_3\delta e^{\Lambda^* T^{**}},
 \end{aligned} \end{equation*}
which contradicts \eqref{0502n111}. Hence, $T^{\min}\neq T^{**}$.
We immediately see that \eqref{201702092227xx} holds.
Moreover, by the relation \eqref{201702092227xx} and the definition of $T^{\min}$, we see that $T^\delta<T^*\leqslant T^{\max}$, and thus $(\eta,(u,\vartheta),q)\in C^0(\overline{I_T},H^{5,1}_1\times  H^4_0\times \underline{H}^3)$ for some $T\in (T^\delta,T^{\max})$.

Noting that $\sqrt{\epsilon_0}\leqslant m_0/2C_5$ and \eqref{ereroe} holds for $t=T^\delta$, thus, making use of \eqref{0501}, \eqref{timesxx}  and \eqref{ereroe}, we can easily deduce the following instability relation:
$$
 \begin{aligned}
\|  \chi (T^\delta)
\|_{{L^1}}\geqslant &
\| \chi^{\mm{a}}   (T^\delta) \|_{{L^1}}
-
\|  \chi^{\mm{d}} (T^\delta)
\|_{{L^1}} \\
>&   \delta e^{\Lambda^* T^\delta }( \| \tilde{\chi}^{0} \|_{L^1}- C_5\sqrt{\delta e^{\Lambda^* T^\delta }}) \geqslant m_0 \epsilon_0 /2,\\
\end{aligned}
$$
where $\chi=\eta_{\mm{h}}$, $\eta_3$, $u_{\mm{h}}$, $u_3$, and $\vartheta$ (we can further take $\chi=\partial_3\eta_{\mm{h}}$ and  $\partial_3\eta_3$ for $Q\neq 0$).
This completes the proof of  Theorem \ref{thmsdfa3} by taking $\epsilon:= m_0\epsilon_0 /2$.

\vspace{4mm} \noindent\textbf{Acknowledgements.}
The research of Fei Jiang   was supported by NSFC (Grant No. 11671086)
 and the research of Song Jiang by NSFC (Grant Nos. 11631008 and 11371065).

\renewcommand\refname{References}
\renewenvironment{thebibliography}[1]{ %
\section*{\refname}
\list{{\arabic{enumi}}}{\def\makelabel##1{\hss{##1}}\topsep=0mm
\parsep=0mm
\partopsep=0mm\itemsep=0mm
\labelsep=1ex\itemindent=0mm
\settowidth\labelwidth{\small[#1]}%
\leftmargin\labelwidth \advance\leftmargin\labelsep
\advance\leftmargin -\itemindent
\usecounter{enumi}}\small
\def\newblock{\ }
\sloppy\clubpenalty4000\widowpenalty4000
\sfcode`\.=1000\relax}{\endlist}
\bibliographystyle{model1b-num-names}

\end{document}